\documentclass{amsart}
\usepackage{color}

\usepackage{amsmath} 
\usepackage{amssymb}
\usepackage{mathrsfs}

\newtheorem{theorem}{Theorem}[section] 
\newtheorem{claim}[theorem]{Claim}
\newtheorem{subc}[theorem]{Subclaim}
\newtheorem{cld}[theorem]{Claim/Definition}
\newtheorem{lemma}[theorem]{Lemma} 
 
\newtheorem{conclusion}[theorem]{Conclusion}
\newtheorem{observation}[theorem]{Observation} 
 
\newtheorem{convention}[theorem]{Convention}

\theoremstyle{definition}
\newtheorem{definition}[theorem]{Definition}
\newtheorem{dc}[theorem]{Definition/Claim}

\newtheorem{problem}[theorem]{Problem}
\newtheorem{discussion}[theorem]{Discussion}

\theoremstyle{remark}
\newtheorem{remark}[theorem]{Remark}

\newtheorem{notation}[theorem]{Notation}

\newcommand{\wk}{{\rm wk}}
\newcommand{\Dp}{{\rm Dp}}

\newcommand{\Lim}{{\rm Lim}}

\newcommand{\otp}{{\rm otp}}
\newcommand{\CON}{{\rm CON}}
\newcommand{\FS}{{\rm FS}}
\newcommand{\FSI}{{\rm FSI}}
\newcommand{\GCH}{ffff{\rm GCH}}
\newcommand{\ZFC}{{\rm ZFC}}
\newcommand{\tr}{{\rm tr}}

\newcommand{\dep}{{\rm dp}}

\newcommand{\Ord}{{\rm Ord}}

\newcommand{\id}{{\rm id}}

\newcommand{\Min}{{\rm Min}}
\newcommand{\Dom}{{\rm Dom}}
\newcommand{\Rang}{{\rm Rang}}

\newcommand{\ic}{{\rm ic}}
\newcommand{\spl}{{\rm spl}}
\newcommand{\exn}{{\phi}}
\newcommand{\exm}{{\epsilon}}
\newcommand{\fsi}{{\rm fsi}}
\newcommand{\cln}{{\rm cln}}
\newcommand{\extheta}{{\partial}}
\newcommand{\exl}{{\varkappa}}
\newcommand{\snn}{{\sn}}

\newcommand{\lqq}{{{`}{`} }}

\newcommand{\rest}{{\restriction}}
\newcommand{\dom}{{\rm dom}}

\newcommand{\wilog}{{\rm without loss of generality}}

\newcommand{\then}{{\underline{then}}}
\newcommand{\when}{{\underline{when}}}
\newcommand{\Then}{{\underline{Then}}}

\newcommand{\If}{{\underline{if}}}
\newcommand{\Iff}{{\underline{iff}}}
\newcommand{\mn}{{\medskip\noindent}}
\newcommand{\sn}{{\smallskip\noindent}}

\newcommand{\cA}{{\mathscr A}}

\newcommand{\gb}{{\mathfrak b}}
\newcommand{\ga}{{\mathfrak a}}

\newcommand{\cD}{{\mathscr D}}
\newcommand{\cJ}{{\mathscr J}}
\newcommand{\gd}{{\mathfrak d\/}} 
\newcommand{\gC}{{\mathfrak C\/}} 

\newcommand{\cH}{{\mathscr H}}

\newcommand{\cI}{{\mathscr I}}

\newcommand{\bbL}{{\mathbb L}}

\newcommand{\cW}{{\mathscr W}}

\newcommand{\bbP}{{\mathbb P}}
\newcommand{\cP}{{\mathscr P}}

\newcommand{\gK}{{\mathfrak K}}
\newcommand{\gu}{{\mathfrak u}}
\newcommand{\gs}{{\mathfrak s}}
\newcommand{\bbQ}{{\mathbb Q}}
\newcommand{\bbR}{{\mathbb R}}

\newcommand{\cT}{{\mathscr T}}
\newcommand{\gt}{{\mathfrak t}}

\newcommand{\cf}{{\rm cf}}
\newcommand{\pr}{{\rm pr}}

\newcount\skewfactor
\def\mathunderaccent#1#2 {\let\theaccent#1\skewfactor#2
\mathpalette\putaccentunder}
\def\putaccentunder#1#2{\oalign{$#1#2$\crcr\hidewidth
\vbox to.2ex{\hbox{$#1\skew\skewfactor\theaccent{}$}\vss}\hidewidth}}
\def\name{\mathunderaccent\tilde-3 }

\newenvironment{PROOF}[2][\proofname.]
   {\begin{proof}[#1]\renewcommand{\qedsymbol}{\textsquare$_{\rm #2}$}}
   {\end{proof}}

\usepackage[hidelinks]{hyperref}
\begin{document}
\makeatletter\def\shfiuwefootnote{\gdef\@thefnmark{}\@footnotetext}\makeatother\shfiuwefootnote{Version 2021-08-08. See \url{https://shelah.logic.at/papers/700a/} for possible updates.}

\title {Are $\ga$ and $\gd$ your cup of tea;\\ revisited \\ 
Sh:700a}

\author {Saharon Shelah}
\address{Einstein Institute of Mathematics\\
Edmond J. Safra Campus, Givat Ram\\
The Hebrew University of Jerusalem\\
Jerusalem, 91904, Israel\\
 and \\
 Department of Mathematics\\
 Hill Center - Busch Campus \\ 
 Rutgers, The State University of New Jersey \\
 110 Frelinghuysen Road \\
 Piscataway, NJ 08854-8019 USA}
\email{shelah@math.huji.ac.il}
\urladdr{http://shelah.logic.at}
\thanks{The author thanks Alice Leonhardt for the beautiful typing.
 First Typed - 98/Dec/16, Public.as  700.
 This is paper 700a, created cannibalizing paper 700} 
 
% Previous version - 2015/Feb/5
  
   % 2020-12-04 15:16 

% Converted to Latex February 2015

\date{2021-02-10}
% 2021-02-10 07:35 --  2021-02-10 07:52 
% 7-7:40 horadeti + 5akalnu
% 2021-02-10 04:51 --  2021-02-10 06:20   5axar kak - bdiqah 
% 2021-02-09 15:26 - klum/ hari5jon nimxaq
% 2021-02-02 11:57 --
% 2021-02-02 11:57 {2021-01-18}
% 2021-01-18 13:05---2021-01-18 17:36 
% % 2021-01-18 12:07 --% 2021-01-18 12:11 
% 2020-12-27 15:57 siyamti 
% 2020-12-27 12:51 -  % 2020-12-27 14:16 
% 2020-12-27 11:28 --2020-12-27 12:22 
% 2020-12-27 11:00 ---
% 2020-12-27 11:00 {2020-12-25a}
% 2020-12-25 13:06 --2020-12-25 14:07 --latex--2020-12-25 14:49  vatah le-sevip 3
% 2020-12-25 09:57 --2020-12-25 10:53 --latex 2020-12-25 11:01 
% 2020-12-24a}
% 2020-12-24 13:12 --  {2020-12-06}
% 2020-12-06 07:49 --% 2020-12-06 09:02 
% 2020-12-05 07:51 -- 2020-12-05 08:56 
% 2020-12-04 17:27 ---% 2020-12-04 18:24 
% 2020-12-04 15:05 ---2020-12-04 15:54   Ilijas Farah 
% 2020-12-04 08:36--2020-12-04 10:25  {2020-12-03}
% 2020-12-03 15:53 -% 2020-12-03 18:42 
--% 2020-12-01 13:05 }  ula xaci javah  qodem -- 2020-12-01 13:35 
% 2020-11-30 19:53 - bevecem kol hayom
% 2020-11-30 12:39   
% 2020-11-30 11:11 --
% 2020-11-29 18:37 --19:02~~~  
% {February 20, 2015}

\begin{abstract} 
This was  % 2020-12-04 15:06 is a % 2020-12-01 13:36 was 
non-essentially revised in late 2020. 
First point is noting that 
the proof of \cite[Th.4.3]{Sh:700} which says  that 
the proof giving the % 2020-12-04 15:07 
consistently  
$ \mathfrak{b} = \mathfrak{d} = \mathfrak{u}   < \mathfrak{a} $ gives
that also $ \mathfrak{s} = \mathfrak{d} $. The proof use
a measurable cardinal and a c.c.c. forcing
so it give large $ \mathfrak{d} $ and assume a large cardinal.

Second point is adding to the results of \S2,\S3 which say that 
(in \S3  % 2020-12-04 15:08 
with no large cardinals) we can force 
$ {\aleph_1} < \mathfrak{b} =   \mathfrak{d} < \mathfrak{a}$. % 2020-12-03 16:00 
We like to have 
$  {\aleph_1} <  \mathfrak{s} \le  \mathfrak{b} = \mathfrak{d} < \mathfrak{a} $. 
For  % 2020-12-01 13:38 
this we allow in \S2,3 the sets $ K_t $ to be
uncountable; this require non-essential  changes.  
In particular,  we replace usually $ {\aleph_0}, {\aleph_1} $ % 2020-12-04 15:10 
by $ \sigma , \partial $.  % 2020-12-04 15:11 
Naturally we can deal with $ \mathfrak{i} $ and similar invariants.

Third we proof read the work again.
To get $ \mathfrak{s} $ we could have retain the countability 
of the member of the $ I_t$-s but the parameters  would change
with $ A \in I_t$, well for a cofinal set of them;
but the present seem simpler.

We intend to continue in \cite{Sh:F2009}.

\underline{Original abstract}  
 We show that consistently, every MAD family has cardinality
strictly bigger than the dominating number, that is ${\ga} > {\gd}$,
thus solving one of the oldest problems on cardinal invariants of the
continuum.  The method is a contribution to the theory of iterated
forcing for making the continuum large. 
\end{abstract}

\maketitle
\numberwithin{equation}{section}
\setcounter{section}{-1}
\newpage

\centerline {Annotated Content} 
\bigskip

\noindent
\S0 \quad Introduction
% \mn

[Was not changed in 2020]  % 2020-12-24 13:14  % 2020-11-30 11:19  18:56 ]
\bigskip

\noindent
\S1 \quad $\CON({\ga} > {\gd})$

\mn
\begin{enumerate}
\item[${{}}$]  [We prove the consistency of the inequality ${\ga} <
{\gd}$, relying on the theory of CS iteration of nep forcing (from
\cite{Sh:630}, this proof is a concise version).
% 2020-12-04 15:17 2020-11-30 
(2020) Was not changed]  % 2020-12-01 13:40  though]
\end{enumerate}
\bigskip

\noindent
\S2 \quad On $\CON({\ga} > {\gd})$ revisited with FS, ideal 
memory of non-well ordered length
\mn
\begin{enumerate}
\item[${{}}$]  
[We use itaration of c.c.c. forcing along a non-well orderd  linear order
 with non-transitive memory.  % 2020-12-04 15:20 
Does not depend on \S1 but use a measurable $\kappa$.
We define ``FSI-template", a depth on the subsets on which
we shall do induction; we are interested just in the cases where the depth is
$< \infty$.  Now the iteration is defined and its properties are proved
simultaneously by induction on the depth.  After we have understood such
iterations sufficiently well, we proceed to prove the consistency in details.

% 2020-12-04 15:20 2020-11-30 
(2020) The change is that we do not require $ K_t $ 
(and the members of $ I_t $) to be countable, 
this require non-essential changes. We also add the 
 promised result].  % 2020-12-24 13:15 
\end{enumerate}
\bigskip

\noindent
\S3 $\quad$ Eliminating the measurable
\mn
\begin{enumerate}
\item[${{}}$]   [In \S2, for checking the criterion which appears there for
having ``${\ga}$ large", we have used ultra-power by some 
$\kappa$-complete ultrafilter.
Here we construct templates of cardinality, e.g. $\aleph_3$ which
satisfy the criterion; by constructing them 
such that any sequence of $\omega$-tuples of
appropriate length has a (big) sub-sequence which is ``convergent" so
some complete $\kappa$-complete filter behave for appropriate
$\kappa$-sequence of names of reals as if it is an ultrafilter and as
if the sequence has appropriate limit.

% 2020-12-01 13:43 
% 2020-12-04 15:21 2020-11-30 
% 2020-12-04 15:20 2020-11-30 
(2020) We add the elimination of the measurable 
also from the result with $ \mathfrak{s} $.]
\end{enumerate}
\bigskip

\noindent
\S4 \quad  On related cardinal invariants
\mn
\begin{enumerate}
\item[${{}}$]  [We prove e.g. the consistency of ${\gu} < {\ga}$, 
starting with a measurable cardinal.  % 2020-12-03 16:00 
Here the forcing notions are not so definable, so this gives a third
proof of the main theorem (but the points which repeat \S3 are not repeated).

% 2020-12-01 13:43 
% 2020-12-04 15:21 2020-11-30
(2020) The addition is noting that the proof give 
also $ \mathfrak{s} = \mathfrak{d} $ in the consistency, 
again not relying on \S2.] 
\end{enumerate}
\newpage

\section {Introduction} 

We deal with the theory of iteration of 
of c.c.c. % 2020-12-03 16:01 
forcing notions for the
continuum and prove $\CON({\ga} > {\gd})$ and related results.  We
present it in several perspectives; so \S2 + \S3 does not depend on
\S1; and \S4 does not depend on \S1, \S2, \S3.  
In \S2 we introduce and investigate iterations which are of finite
support but with so called ideal, weakly 
transitive memory and linear, non well ordered
length and prove $\CON({\ga} > {\gd})$ using a measurable.
In \S4 we answer also related questions (${\gu} < {\ga}$); in \S3, 
relying on \S2 we eliminate the use of a measurable, and in \S1 we 
rely heavily on \cite{Sh:630}.

Very basically, the difference between ${\ga}$ on the one hand 
and ${\gb},{\gd}$ on the other hand which we use is that ${\ga}$ 
speaks on a set, whereas ${\gb}$ is witnessed by a sequence and 
${\gd}$ by a quite directed family; it essentially deals with cofinality;
so every unbounded subsequence is a witness as well, i.e.
the relevant relation is transitive; when ${\gb} =
{\gd}$ things are smooth, otherwise the situation is still similar.
This manifests itself by using ultrapowers for
some $\kappa$-complete ultrafilter (in model theoretic outlook), and by
using ``convergent sequence" (see \cite{Sh:300}
and later \cite[\S2]{Sh:300a}, \cite{Sh:300b}), or the existence of
Av, the average,
from % 2020-12-03 16:04 in 
\cite{Sh:c}) in \S2, \S3, respectively.
The meaning of ``model theoretic outlook", is that by experience set
theorists starting to hear an explanation of the forcing tend to think
of an elementary embedding
$\mathbf j:\mathbf V \rightarrow M$ and then the limit practically does not
make sense (though of course we can translate).  Note that ultrapowers 
by e.g. an ultrafilter on $\kappa$, preserve any witness for a
cofinality of a linear order being $\ge \kappa^+$ (or the 
cofinality of a $\kappa^+$-directed partial order), 
as the set of old elements is cofinal and a cofinal subset of a cofinal
subset is a cofinal subset.   On the other hand, the ultrapower always
``increase" any set of cardinality at least $\kappa$, the 
completeness of the ultrafilter.
\bigskip

\centerline {$* \qquad * \qquad *$}
\bigskip

This (is ${\ga} \le {\gd}$?) is one of 
the oldest problems and well known on cardinal invariants of the continuum
(see [vD]
and Roitman \cite{Mi91}).
It was mostly thought (certainly by me) that 
consistently ${\ga} > {\gd}$ and that the natural way to proceed is by 
CS iteration $\langle \bbP_i,\name{\bbQ}_i:i < \omega_2 \rangle$ of proper
${}^\omega \omega$-bounding forcing notions, starting with $\mathbf V
\models \GCH$, and $|\bbP_i| = \aleph_1$ for $i < \omega_2$ and 
$\name{\bbQ}_i$ ``deal" with one MAD family ${\cA}_i \in \mathbf V^{\bbP_i},
{\cA}_i \subseteq [\omega]^{\aleph_0}$, adding an infinite 
subset of $\omega$ almost disjoint to every $A \in {\cA}_i$.  
The needed iteration theorem holds by \cite[Ch.V,\S4]{Sh:f}, saying that in 
$\mathbf V^{\bbP_{\omega_2}},{\gd} = {\gb} = \aleph_1$ and no 
cardinal is collapsed, \underline{but} 
the single step forcing is not known to exist.  This has been
explained in details in \cite{Sh:666}.

We do not go in this way but in a totally different 
direction involving making the continuum large, so we still do not know 
\begin{problem}
\label{0.1}
Is $\ZFC + 2^{\aleph_0} = \aleph_2 + {\ga} > {\gd}$ consistent?
 
To clarify our idea, let $D$ be a normal ultrafilter on $\kappa$, 
a measurable cardinal and consider a c.c.c. forcing notion $\bbP$ 
and assume we have 
\mn
\begin{enumerate}
\item[$(a)$]   a sequence $\name{\bar f} = \langle \name f_\alpha:
\alpha < \kappa^+ \rangle$ of $\bbP$-names such that $\Vdash_{\bbP} 
``\langle \name f_\alpha:\alpha < \kappa^+ \rangle$ is 
$<^*$-increasing cofinal in ${}^\omega \omega$" 
(so $\name{\bar f}$ exemplifies $\Vdash_{\bbP} ``{\gb} = {\gd} = \kappa^+"$)
\sn
\item[$(b)$]    a sequence $\langle \name A_\alpha:\alpha < 
\alpha^* \rangle$ of $\bbP$-names such that
\newline
$\Vdash_{\bbP} ``\{\name A_\alpha:\alpha < \alpha^*\}$ is
MAD that is $\alpha \ne \beta \Rightarrow 
\name A_\alpha \cap \name A_\beta$ is finite and
$\name A_\alpha \in [\omega]^{\aleph_0}"$.
\end{enumerate}
\mn
Now $\bbP_1 = \bbP^\kappa/D$ also is a c.c.c. forcing notion by 
{\L}o\'s theorem for $\bbL_{\kappa,\kappa}$; 
let $\mathbf j:\bbP \rightarrow \bbP_1$ be the canonical
embedding; moreover, under the canonical identification we have
$\bbP \prec_{\bbL_{\kappa,\kappa}} \bbP_1$.  
So also $\Vdash_{\bbP_1} ``\name f_\alpha \in {}^\omega \omega"$, 
recalling that $\name f_\alpha$ actually consists of $\omega$ maximal
antichains of $\bbP$ (or think of $(\cH(\chi),\in)^\kappa/D,\chi$ large
enough).  Similarly $\Vdash_{\bbP_1} ``\name f_\alpha <^*
\name f_\beta$ if $\alpha < \beta < \kappa^+"$.

Now, if $\Vdash_{\bbP_1} ``\name g \in {}^\omega \omega"$,
then $\name g = \langle \name g_\varepsilon:
\varepsilon < \kappa \rangle/D,\Vdash_{\bbP}
``\name g_\varepsilon \in {}^\omega \omega"$ so for some
$\alpha < \kappa^+$ we have $\Vdash_{\bbP} 
``\name g_\varepsilon <^* \name f_\alpha$
for $\varepsilon < \kappa"$ hence by {\L}o\'s theorem
$\Vdash_{\bbP_1} ``\name g <^* { \name{ f}}_\alpha"$ (so before the
identification this means $\Vdash_{\bbP_1} ``\name g <^* 
\mathbf j(f_\alpha)"$), so $\langle \name f_\alpha:
\alpha < \kappa^+ \rangle$ exemplifies also
$\Vdash_{\bbP_1} ``\gb = \gd = \kappa^+"$.

On the other hand $\langle \name A_\alpha:\alpha <
\alpha^* \rangle$ cannot exemplify that ${\ga} \le \kappa^+$ in
$\mathbf V^{\bbP_1}$ because $\alpha^* \ge \kappa^+$ 
(as $\ZFC \models \gb \le {\ga}$) so $\langle \name A_\alpha:
\alpha < \kappa \rangle/D$ exemplifies that $\Vdash_{\bbP_1} 
``\{\name A_\alpha:\alpha < \alpha^*\}$ is not MAD".

Our original idea here is to start with a $\FS$ 
iteration $\bar{\bbQ}^0 = \langle \bbP^0_i,
\name{\bbQ}^0_i:i < \kappa^+ \rangle$ of nep c.c.c.
forcing notions, $\name{\bbQ}^0_i$ adding a dominating real, (e.g.
by % 2021-01-18 
dominating real = Hechler forcing), for $\kappa$ a measurable 
cardinal and let $D$ be a $\kappa$-complete uniform ultrafilter 
on $\kappa$ and $\chi \gg \kappa$.  
Then let $L_0 = \kappa^+,\bar{\bbQ}^1 = \langle \bbP^1_i,\bbQ^1_i:
i \in L_1 \rangle$ be $\bar{\bbQ}^0$ as 
interpreted in $({\cH}(\chi),\in,<^*_\chi)^\kappa/D$, it
looks like $\bar{\bbQ}^0$ replacing $\kappa^+$ by
$(\kappa^+)^\kappa/D$.  We look at $\Lim(\bar{\bbQ}^0) =
\bigcup  % 2021-01-18 12:10 \limits_{i}
\{ \ 
\bbP^0  % 2021-01-18 12:09 
_i  i < \kappa ^+ \} $    % 2021-01-18 12:11 % 2021-01-18 13:05 
as a subforcing of $\Lim(\bar{\bbQ}^1)$
identifying $\name{\bbQ}_i$ with $\name{\bbQ}_{\mathbf j_0(i)},
\mathbf j_0$ the canonical elementary embedding of
$\kappa^+$ into $(\kappa^+)^\kappa/D$ (no Mostowski collapse!).  
We continue to define $\bar{\bbQ}^n$ and then $\bar{\bbQ}^\omega$ 
as the following limit: for the original\footnote{
which mean not the ones added by taking ultrapowers
}
$i \in \kappa^+$, we use 
the definition, otherwise we use direct limit (``founding fathers 
privilege" you may say).  
So $\bbP^i = \Lim(\bar{\bbQ}^i)$ is 
$\lessdot$-increasing, continuous when $\cf(i) > \aleph_0$; so now we
have a kind of iteration with so called ideal, weakly 
transitive memory and a not well founded base.  We continue
$\kappa^{++}$ times.  Now in $\mathbf V^{\Lim(\bar{\bbQ}^{\kappa^{++}})}$, the 
original $\kappa^+$ generic reals exemplify ${\gb} = {\gd} =
\kappa^+$, so we know that ${\ga} \ge \kappa^+$.  To finish assume
$p \Vdash ``\{\name A_\gamma:\gamma < \kappa^+\}
\subseteq [\omega]^{\aleph_0}$ is a MAD family".  Each name
$\name A_\gamma$ is a ``countable object" 
and so depends on countably many conditions, so all of 
them are in $\Lim(\bar{\bbQ}^i)$ for some $i < \kappa^{++}$.  
In the next stage, $\bar{\bbQ}^{i+1},\langle \name A_\gamma:
\gamma < \kappa \rangle/D$ is a
name of an infinite subset of $\omega$ almost disjoint to
$\name A_\beta$ for each $\beta < \kappa^+$, contradiction.

All this is a reasonable scheme.  This is done in \S1 but rely on
``nep forcing" from \cite{Sh:630}.  But a self contained another 
approach is in \S2,\S3, where the meaning of the iteration is 
more on the surface (and also, in \S3, help to eliminate the 
use of large cardinals).  In \S4 we deal with
the case of an additional cardinal invariant, ${\gu}$.

Note that just using FS iteration on a non well-ordered linear 
order $L$ (instead of an ordinal) is impossible by a theorem 
of Hjorth.  On nonlinear orders for iterations (history and 
background) see \cite{Sh:670}.  
On iteration with non-transitive memory 
see \cite{Sh:592}, \cite{Sh:619} and in particular \cite[\S3]{Sh:619}.
\end{problem}
\bigskip

Continuing this work J. Brendle has proved the consistency
of $\cf({\ga})     % 2020-12-03 16:25 
= {\aleph_0}$, (note that in \ref{abd.3} we have assumed 
$\lambda = \lambda^{\aleph_0}$ in $\mathbf V$ hence $\cf(\lambda) >
{\aleph_0}$ even in $\mathbf V^{\bbP}$).

I thank Heike Mildenberger and Juris Steprans for their helpful
comments.  After publication this was revised simplifying \S2.

\begin{notation} \label{z2}
1)  $ \mathbb{P}, \mathbb{Q} $ denote forcing notions

\noindent 
2) Let $ \mathbb{P} \subseteq \mathbb{Q} $  means that for 
$ p, q \in \mathbb{P} $ we have 
$ p < _ \mathbb{P} q$ iff $ p < _ \mathbb{Q} q$ 

\noindent 
3) let $ \mathbb{P} \subseteq _{ic} \mathbb{Q} $
iff 
$ \mathbb{P} \subseteq \mathbb{Q} $  and for every 
$ p,q \in \mathbb{P}$  we have
$ p, q $ are compatible in $ \mathbb{P}$  iff they are
compatible in $ \mathbb{Q} $

\noindent 
4) Let $ \mathbb{P} \lessdot \mathbb{Q}$  iff 
$ \mathbb{P} \subseteq _{ic} \mathbb{Q} $
and every maximal anti-chain of $ \mathbb{P} $
is a maximal anti-chain  of $ \mathbb{Q} $
 
\end{notation}

\begin{convention} \label{z5}
1)
When using $ \mathfrak{t} , (\mathfrak{t} , \bar{ K } ) $
we mean as in  Def \ref{ad.1}.

\noindent 
2) When using $ (\mathfrak{t} , \bar{ K } , \bar{ u } , \Lim( \bar{ \mathbb{Q} }))  $
we mean as in \ref{ad.5}

\noindent 
3) We may write $ I_t $ instead $ I^ \mathfrak{t} _ t $ or 
$ I^ \mathbf{q} _t$  when $ \mathfrak{t}, % 2021-01-18 13:09 /
\mathbf{q} $ is clear 
from the contecxt.  % 2020-12-04 15:23 

\noindent 
4) Dealing with e.g. $ \mathfrak{t} ^\zeta $ we may write
$ \mathfrak{t} [\zeta ]$ in subscript and superscripts.  % 2020-12-04 15:24 
\end{convention} 
\newpage

\section {On Con$({\ga} > {\gd})$} 

In this section, we look at it in the context of \cite{Sh:630} and we 
use a measurable  cardinal.
\bigskip

\begin{definition}
\label{da.1}
1) Given sets $A_\ell$ of ordinals
for $\ell < n$, we say ${\cT}$ is an $(A_0,\dotsc,A_{n-1})$-tree if
${\cT} = \bigcup\limits_{k < \omega} {\cT}_k$ where ${\cT}_k \subseteq 
\{(\eta_0,\dotsc,\eta_\ell,\dotsc,\eta_{n-1}):
\eta_\ell \in {}^k(A_\ell)$ for $\ell < n\}$ and ${\cT}$
is ordered by $\bar \eta \le_{\cT} \bar \nu \Leftrightarrow 
\bigwedge\limits_{\ell < n} \eta_\ell \trianglelefteq \nu_\ell$ and we let 
$\bar \eta \upharpoonleft k_1 =: \langle \eta_\ell \restriction 
k_1:\ell < n \rangle$ and demand $\bar \eta \in {\cT}_k 
\wedge  k_1 < k \Rightarrow \bar \eta \upharpoonleft k_1 \in
{\cT}_{k_1}$.  We call ${\cT}$ locally countable if $k \in [1,\omega) 
\wedge  \bar \eta \in {\cT}_k \Rightarrow |\{\bar \nu \in {\cT}_{k+1}:
\bar \eta \le_{\cT} \bar \nu\}| \le \aleph_0$.
Let $\lim(\cT) = \{ \langle \eta_\ell:\ell < n \rangle:\eta_\ell \in
{}^\omega(A_\ell)$ for $\ell < n $ % 2020-12-24 13:21  k$ 
and $m < \omega \Rightarrow 
\langle \eta_\ell \restriction m:\ell < n \rangle \in {\cT}\}$.  

Lastly, for $n_1 \le n$ we let prj lim$_{n_1}({\cT}) = 
\{\langle \eta_\ell:\ell < n_1 \rangle:$ for some 
$\eta_{n_1},\dotsc,\eta_{n-1}$ we have $\langle \eta_\ell:\ell < n 
\rangle \in \lim({\cT})\}$; and if $n_1$ is omitted we mean $n_1 = n-1$.

\noindent
2)

\begin{equation*}
\begin{array}{clcr}
\gK = \{\bar{\cT}:&\text{for some sets } A,B \text{ of ordinals we have} \\
  &(i) \quad \bar{\cT} = ({\cT}_1,{\cT}_2), \\
  &(ii) \quad {\cT}_1 \text{ is a locally countable } (A,B)\text{-tree},\\
  &(iii) \quad {\cT}_2 \text{ is a locally countable } 
(A,A,B)\text{-tree, and} \\
  &(iv) \quad \bbQ_{\bar{\cT}} =: (\text{prj lim}(\cT_1), \text{ prj lim}
({\cT}_2)) \text{ is a c.c.c. forcing notion} \\
  &\qquad \text{absolute   % 2021-01-18 13:09 ly 
  under c.c.c. forcing notions (see below)}\}
\end{array}
\end{equation*}

\mn
2A) We say that $\bbQ_{\bar{\cT}}$ is 
c.c.c. absolute   % 2021-01-18 13:10 ly 
for c.c.c. forcing if:
for c.c.c. forcing notions $\bbP \lessdot \bbR \text{ we have } 
\bbP * \name{\bbQ}_{\bar{\cT}} \lessdot \bbR *
\bbQ_{\name{\bar \cT}}$ (though not necessarily 
$\bbQ^{\mathbf V^{\bbP}}_{\bar{\cT}} \lessdot 
\bbQ^{\mathbf V^{\bbR}}_{\bar{\cT}}$ in 
$\mathbf V^{\bbR}$) so membership, order, 
non-order, compatibility, noncompatibility and being predense over $p$
in the universe $\mathbf V^{\bbP}$, are 
preserved in passing to $\mathbf V^{\bbR}$, note that 
predense sets belong to $\mathbf V^{\bbP}$ 
(the $\bbQ_{\bar{\cT}}$'s are snep, from \cite{Sh:630} with 
slight restriction).  Similarly we define ``$\bbQ_{\bar{\cT}_1} 
\lessdot \bbQ_{\bar{\cT}_2}$ 
absolute % 2021-01-18 13:12 ly 
under c.c.c. forcing" (compare with
\ref{ad.5}, clause (A)(a)(iii) in the definition).

\noindent
3) For a set or class A of ordinals, ${\gK}^\kappa_A$ 
is the family of $\bar{\cT} \in {\gK}$ which
are a pair of objects, the first an $(A,B)$-tree and the second an
$(A,A,B)$-tree % 2020-12-24 13:22 s 
for some $B$ such 
that $|{\cT}_1| \le \kappa,|{\cT}_2| \le \kappa$.  For a
cardinal $\kappa$ and a pairing function pr with inverses 
$\pr_1,\pr_2$, let ${\gK}^\kappa_{\pr_1,\gamma} = 
{\gK}^\kappa_{\{\alpha:\pr_1(\alpha)=\gamma\}}$ and
$\gK^\kappa_{\pr_1,<\gamma} = 
\gK^\kappa_{\{\alpha:\pr_1(\alpha) < \gamma\}}$.  
Let $|\bar{\cT}| = |{\cT}_1| + |{\cT}_2|$. 

\noindent
4) Let $\bar{\cT},\bar{\cT}' \in {\gK}$, we say $\mathbf f$ is 
an isomorphism from $\bar{\cT}$ onto $\bar{\cT}'$ \when \,
$\mathbf f = (f_1,f_2)$ and for $m = 1,2$ we have: $f_m$ is a one-to-one
function from ${\cT}_m$ onto ${\cT}'_m$ 
preserving the level (in the respective trees), preserving the relations 
$x = y \upharpoonleft k,x \ne y \upharpoonleft k$
and if $f_2((\eta_1,\eta_2,\eta_3)) 
= (\eta'_1,\eta'_2,\eta'_3),f_1((\nu_1,
\nu_2)) = (\nu'_1,\nu'_2)$ then $[\eta_1 
= \nu_1 \Leftrightarrow \eta'_1 = \nu'_1],
[\eta_2 = \nu_1 \Leftrightarrow \eta'_2 = \nu'_1]$. 

In this case let $\hat{\mathbf f}$ be the isomorphism induced 
by $\mathbf f$ from $\bbQ_{\bar{\cT}}$ onto $\bbQ_{\bar{\cT}'}$.
\end{definition}

\begin{definition}
\label{da.2}
For $\bar{\cT}',\bar{\cT}'' \in {\gK}$ let $\bar{\cT}' \le_{\gK} 
\bar{\cT}''$ mean:
\mn
\begin{enumerate}
\item[$(a)$]   ${\cT}'_\ell \subseteq {\cT}''_\ell$ (as trees) 
for $\ell =1,2$
\sn
\item[$(b)$]  if $\ell \in \{1,2\}$ and $\bar \eta \in {\cT}''_\ell
\backslash {\cT}'_\ell$ and $\bar \eta \upharpoonleft k \in
{\cT}'_\ell$ \then \, $k \le 1$
\sn
\item[$(c)$]  $\bbQ_{\bar{\cT}'} \lessdot \bbQ_{\bar{\cT}''}$ 
(absolute  % 2021-01-18 13:12 ly 
under c.c.c. forcing); note that by (a) + (b) we have: 
$x \in \bbQ_{\bar{\cT}'} \Rightarrow x \in
\bbQ_{\bar{\cT}''}$ and $\bbQ_{\bar{\cT}'} 
\models x \le y \Rightarrow \bbQ_{\bar{\cT}''} \models x \le y$.
\end{enumerate}
\end{definition}

\begin{remark}  
The definition is tailored such that the union of an
increasing chain will give a forcing notion which is the union.
\end{remark}

\begin{cld}
\label{da.3}
0) $\le_{\gK}$ is a partial order of ${\gK}$. 

\noindent
1) Assume $\langle \bar{\cT}[i]:
i < \delta \rangle$ is $\le_{\gK}$-increasing and 
$\bar{\cT}$ is defined by $\bar{\cT} = \bigcup\limits_{i} \bar{\cT}[i]$ 
that is ${\cT}_m = \bigcup\limits_{i < \delta} {\cT}_m[i]$ for $m=1,2$ 
\then \,
\mn
\begin{enumerate}
\item[$(a)$]  $i < \delta \Rightarrow \bar{\cT}[i] \le_{\gK} \bar{\cT}$
\sn
\item[$(b)$]  $\bbQ_{\bar{\cT}} = \bigcup\limits_{i < \delta} 
\bbQ_{\bar{\cT}[i]}$.
\end{enumerate}
\mn
2) Assume $\bar{\cT}',\bar{\cT} \in {\gK}$.  \Then \,
there is $\bar{\cT}'' \in {\gK}$ such that $\bar{\cT}' 
\le_{\gK} \bar{\cT}''$ and $\bbQ_{\bar{\cT}''}$
is isomorphic to $\bbQ_{\bar{\cT}'} * \name{\bbQ}_{\bar{\cT}}$ 
and this is absolute by c.c.c. forcing.  Moreover, there is 
such an isomorphism extending the identity map from 
$\bbQ_{\bar{\cT}'}$ into $\bbQ_{\bar{\cT}''}$.

\noindent
3) There is $\bar{\cT} \in {\gK}^{\aleph_0}_\omega$ such that 
$\bbQ_{\bar{\cT}}$ is the trivial forcing. 

\noindent
4) There is $\bar{\cT} \in {\gK}^{\aleph_0}_\omega$ 
such that $\bbQ_{\bar{\cT}}$ is the dominating real forcing.
\end{cld}

\begin{PROOF}{\ref{da.3}}
See \cite{Sh:630}.
\end{PROOF}

\begin{claim}
\label{da.4}
1) Assume $\bar{\cT}[\gamma] \in {\gK}_{\pr_1,\gamma}$ for 
$\gamma < \gamma(*)$.  \Then \, for each $\alpha
\le \gamma(*)$ there is $\bar{\cT} \langle \alpha \rangle$
such that $\bbQ_{\bar{\cT} \langle \alpha \rangle}$ is $\bbP_\alpha$ where 
$\langle \bbP_\gamma,\name{\bbQ}_\beta:
\gamma \le \gamma(*),\beta < \gamma(*) \rangle$ is 
an $\FS$-iteration and $\name{\bbQ}_\beta = 
(\bbQ_{\bar{\cT}[\beta]})^{{\mathbf V}[\bbP_\beta]}$ 
and $\bar{\cT} \langle \alpha \rangle \in {\gK}_{\pr_1,< \alpha}$ 
and $\bar{\cT} \langle \alpha_1 \rangle \le_{\gK}
\bar{\cT} \langle \alpha_2 \rangle$ for $\alpha_1 \le \alpha_2 \le
\gamma(*),\bar{\cT}[\gamma] \le_{\gK} \bar{\cT} \langle \alpha
\rangle$ for $\gamma < \alpha \le \gamma(*)$.  We write $\bar{\cT}
\langle \alpha \rangle = \sum\limits_{\gamma < \alpha} \bar{\cT}[\gamma]$.

\noindent
2) In part (1), for each $\gamma < \gamma(*)$ there is $\bar{\cT}' \in
{\gK}_{\pr_1,\gamma}$ such that $\bar{\cT}',\bar{\cT}$ 
are isomorphic over $\bar{\cT}[\gamma]$ hence $\bbQ_{\bar{\cT}'},
\bbQ_{\bar{\cT}}$ are isomorphic over $\bbQ_{\bar{\cT}[\gamma]}$. 

\noindent
3) If in addition ${\cT}[\gamma] \le_{\gK} {\cT}'[\gamma] \in 
{\gK}_{\pr_1,\gamma}$ for $\gamma < \gamma(*)$ and $\langle 
\bbP_\gamma,\name{\bbQ}'_\beta:\gamma \le \gamma(*),\beta < \gamma(*) 
\rangle$ is an FS iteration as above with $\bbP'_{\gamma(*)} = 
\bbQ_{\bar{\cT}'}$, \then \, we  
can % 2021-01-18 13:18 
find such $\bar{\cT}'$ with $\bar{\cT} 
\le_{\gK} \bar{\cT}  % 2021-01-18 13:19 '
$.
\end{claim}

\begin{PROOF}{\ref{da.4}}
Straightforward.
\end{PROOF}

\begin{claim}
\label{da.6}
Assume
\mn
\begin{enumerate}
\item[$(a)$]  $\kappa$ is a measurable cardinal
\sn
\item[$(b)$]  $\kappa < \mu = \cf(\mu) < \lambda =
\cf(\lambda) = \lambda^\kappa$ and $(\forall \alpha < \mu)
(|\alpha|^{\aleph_0} < \mu)$ for simplicity.
\end{enumerate}
\mn
\Then \, for some c.c.c. forcing notion $\bbP$ of 
cardinality $\lambda$, in $\mathbf V^{\bbP}$ we have: 
$2^{\aleph_0} = \lambda,{\gd} = {\gb} = \mu$ and ${\ga} = \lambda$.
\end{claim}

\begin{PROOF}{\ref{da.6}}
We choose by induction on $\zeta \le \lambda$ the following
objects satisfying the following conditions:
\mn
\begin{enumerate}
\item[$(a)$]  a sequence $\langle \bar{\cT}[\gamma,\zeta]:\gamma < \mu
\rangle$
\sn
\item[$(b)$]  $\bar{\cT}[\gamma,\zeta] \in {\gK}^\lambda_{\pr_1,\gamma}$
\sn
\item[$(c)$]  $\xi < \zeta \Rightarrow \bar{\cT}[\gamma,\xi] 
\le_{\gK} \bar{\cT}[\gamma,\zeta]$
\sn
\item[$(d)$]  if $\zeta$ limit then $\bar{\cT}[\gamma,\zeta] =
\bigcup\limits_{\xi < \zeta} \bar{\cT}[\gamma,\xi]$
\sn
\item[$(e)$]  if $\gamma < \mu,\zeta = 1$ \then \,
$\bbQ_{\bar{\cT}}[\gamma,\zeta]$ is the $\bbQ_{\dom}$, 
dominating real forcing = Hechler forcing
\sn
\item[$(f)$]  if $\gamma < \mu,\zeta = \xi +1 >1$ and $\xi$ is even, 
\then \, $\bar{\cT}[\gamma,\zeta]$ is isomorphic to 
$\bar{\cT}\langle \gamma +1,\xi \rangle$
over $\bar{\cT}[\gamma,\xi]$ say by $\mathbf j_{\gamma,\xi}$ where
$\bar{\cT} \langle \gamma +1,\xi \rangle =: 
\sum\limits_{\beta \le \gamma} \bar{\cT}[\beta,\xi]$ and let 
$\hat{\mathbf j}_{\gamma,\xi}$ be the isomorphism induced from
$\bbQ_{\bar{\cT} \langle \gamma +1,\xi \rangle}$ onto
$\bbQ_{\bar{\cT}}[\gamma,\zeta]$ over $\bbQ_{\bar{\cT}[\gamma,\xi]}$
\sn
\item[$(g)$]  if $\gamma < \mu,\zeta = \xi +1,\xi$ odd, \then \,
$\bar{\cT}[\gamma,\zeta]$ is almost isomorphic to 
$(\bar{\cT}[\gamma,\xi])^\kappa/D$ over $\bar{\cT}_{[\gamma,\xi]}$ 
which means that we say $\mathbf j_{\gamma,\xi}$ is an 
isomorphism from $(\bar{\cT}[\gamma,\xi])^\kappa/D$
onto $\bar{\cT}[\gamma,\zeta]$ such that by 
$\mathbf j_{\gamma,\xi}$, % 2020-12-24 13:47  \, 
$ \langle x:\varepsilon < \kappa \rangle/D$ is mapped onto $x$.
\end{enumerate}
\mn
There is no problem to carry the definition.  Let $\bbP_\zeta = 
\bbQ_{{\bar{\cT}} \langle \mu,\zeta \rangle}$ where $\bar{\cT}
\langle \mu,\zeta \rangle =: \sum\limits_{\gamma < \mu} \bar{\cT}
[\gamma,\zeta]$ for $\zeta \le \lambda,
\bbP = \bbP_\lambda$ and $\bbP_{\gamma,\zeta} = 
\bbQ_{{\bar{\cT}} \langle \gamma,\zeta \rangle}$.  

Now
\mn
\begin{enumerate}
\item[$\boxtimes_1$]  $|\bbP| \le \lambda$.
\end{enumerate}
\mn
[Why?  As we prove by induction on $\zeta \le \lambda$ that: each
$\bar{\cT}[\gamma,\zeta]$ and
$\sum\limits_{\gamma \le \mu} \bar{\cT}[\gamma,\lambda]$ has
cardinality $\le \lambda$.  Hence for $\gamma < \mu$ we have: 
the forcing notion
$\bbQ_{\bar{\cT}[\gamma,\lambda]}$ in the universe 
$\mathbf V^{{\bbQ}_{\bar{\cT} \langle \gamma,\lambda \rangle}}$ has
cardinality $\le \lambda^{\aleph_0} = \lambda$.]
\mn
\begin{enumerate}
\item[$\boxtimes_2$]  in 
$\mathbf V^{\bbP}$ we have ${\gb} = {\gd} =
\mu $ 
   % 2021-01-18 13:20 \lambda$.
\end{enumerate}
\mn
[Why?  Let $\name\eta_\gamma$ be the 
$\bbQ_{\bar{\cT}[\gamma,1]}$-name of the dominating 
real (see clause (e)).  
As $\bar{\cT}[\gamma,1] \le_{\gK} \bar{\cT}[\gamma,\lambda]$, 
clearly $\name \eta_\gamma$ is also a
$\bbQ_{\bar{\cT}[\gamma,\lambda]}$-name of a dominating real,
but this is preserved by (forcing by) $\bbP_\gamma$ hence 
$\Vdash_{\bbP_{\gamma+1}} ``\name \eta_\gamma$ dominates % 2020-12-24 13:23 s
$({}^\omega \omega)^{\mathbf V[{\bbP}_{\gamma,\lambda}]}"$.   But $\langle
\bbP_{\gamma,\lambda}:\gamma < \mu \rangle$ is $\lessdot$-increasing with
union $\bbP$ and $\cf(\mu) = \mu > \aleph_0$ so
$\Vdash_{\bbP} ``\langle \name \eta_\gamma:\gamma < \mu \rangle$ is
$<^*$-increasing and dominating", so the conclusion follows.]

We shall 
prove below that ${\ga} \ge \lambda$, together this finishes the
proof (note that it implies $2^{\aleph_0} \ge \lambda$ hence as
$\lambda = \lambda^{\aleph_0}$ by
$\boxtimes_1$ we get $2^{\aleph_0} = \lambda$)
\mn
\begin{enumerate}
\item[$\boxtimes_3$]  $\Vdash_{\bbP} ``{\ga} \ge \lambda"$.
\end{enumerate}
\mn
So assume $p \Vdash ``\name{\cA} = \{\name A_i:i < \theta\}$ is a MAD family,
i.e. ($\theta \ge \aleph_0$ and)
\mn
\begin{enumerate}
\item[$(i)$]  $\name A_i \in [\omega]^{\aleph_0}$,
\sn
\item[$(ii)$]  $i \ne j \Rightarrow |\name A_i \cap
\name A_j| < \aleph_0$ and 
\sn
\item[$(iii)$]   $\name{\cA}$ is maximal under $(i) + (ii)$".
\end{enumerate}
\mn
Without loss of generality $\Vdash_{\bbP} ``\name A_i \in
[\omega]^{\aleph_0}"$. 

As always ${\ga} \ge {\gb}$, by $\boxtimes_2$ we know that 
$\theta \ge \mu$, and toward contradiction assume $\theta < \lambda$.  
For each $i < \theta$ and $m < \omega$ there is a maximal anti-chain 
$\langle p_{i,m,n}:n < \omega \rangle$ of $\bbP$ 
and a sequence $\langle \mathbf t_{i,m,n}:n < \omega \rangle$ 
of truth values such that $p_{i,m,n} 
\Vdash_{\bbP} ``n \in \name A_i$ iff $\mathbf t_{i,m,n}$ is truth".  
We can find a countable $w_i \subseteq \mu$
such that: $[\bigcup\limits_{m,n < \omega} \Dom(p_{i,m,n}) \subseteq
w_i],p_{i,m,n} \in \bbQ_{\Sigma\{\bar{\cT}[\gamma,\lambda]:\gamma
\in w_i\}}$, moreover,  $\gamma \in \Dom(p_{i,m,n}) 
\Rightarrow p_{i,m,n}(\gamma)$ is a $\bbQ_{\sum\{\bar{\cT}
[\beta,\lambda]:\beta \in \gamma \cap w_i\}}$-name.  

Note that $\bbQ_{\sum\{\bar{\cT}[\beta,\lambda]:\beta \in \gamma \cap w_i,
i < \theta\}} \lessdot \bbQ_{\sum\{\bar{\cT}_\beta:\beta < 
\gamma\}}$, see \cite{Sh:630}.

Clearly for some even $\zeta < \lambda$, we have $\{p_{i,m,n}:i < \theta,m <
\omega$ and $n < \omega \} \subseteq \bbQ_{\sum\{\bar{\cT}[\beta,\zeta]:
\beta < \mu\}}$.  Now for some stationary 
$S \subseteq \{\delta < \mu:\cf(\delta) = \kappa\}$ and $w^*$ we have:
$\delta \in S \Rightarrow w_\delta \cap \delta = w^*$ and $\alpha <
\delta \in S \Rightarrow w_\alpha \subseteq \delta$.  Let $\langle
\delta_\varepsilon:\varepsilon < \kappa \rangle$ be an increasing sequence of
members of $S$, and $\delta^* = \bigcup\limits_{\varepsilon < \kappa} 
\delta_\varepsilon$.  The definition of $\langle \bar{\cT}
[\gamma,\zeta +1]:\gamma < \mu \rangle,\langle \bar{\cT}
[\gamma,\zeta + 2]:\gamma < \mu \rangle$ was made to get a name of 
an infinite $\name A \subseteq \omega$
almost disjoint to every $\name A_\beta$ for $\beta <
\theta$ (in fact $(\sum\limits_{\gamma < \mu} \bbQ_{\bar{\cT}[\gamma,\zeta]})
^\kappa /D$ can be $\lessdot$-embedded into $\sum\limits_{\gamma < \mu}
\bbQ_{\bar{\cT}[\gamma,\zeta + 2]})$.
\end{PROOF}

\begin{remark}  In later proofs in \S2 we give more details.
\end{remark}
\newpage

\section {On Con$({\ga} > {\gd})$ revisited with FS, with ideal
memory, non-well ordered length}
\bigskip

(Pre 2020 introduction to this section)
We first define the FSI-templates, telling us how do we iterate along a
linear order $L$; we think of having for each $t \in L$, a forcing notion
$\bbQ_t$, say adding a generic $\name \nu_t$, and $\bbQ_t$ will
really be $\cup\{\bbQ^{\mathbf V[\langle \name \nu_s:s \in A \rangle]}:
A \in I_t\}$ 
where $I_t$ 
is % 2021-01-18 13:30 
an ideal of subsets of $\{s:s <_L t\}$;
so $\bbQ_t$ in the nice case is a definition, e.g. as in
\ref{da.1}(2A).  
In our application this definition
is constant, but we treat a more general case, so 
$\name{\bbQ}_t$ may be defined using parameters from
$\mathbf V[\langle \name \nu_s:s \in K_t \rangle],K_t$ a
subset of $\{s:s <_L t\}$ so the reader may consider only the case
$t \in L \Rightarrow K_t = \emptyset$.  
In part (3) of Definition \ref{ad.1} instead distinguishing
``$\zeta$ successor, $\zeta$ limit" we can consider the two cases for each
$\zeta$.  The depth of $L$ is the ordinal on which 
our induction rests (as \otp$(L)$ is inadequate).

Now (2020)  we allow uncountable $ K_t$-s (and similarly
$ \name{ \eta }, \nu $), a non-essential change.
\bigskip

\begin{definition}
\label{ad.1}
1)  An FSI-template (= finite support iteration template) ${\gt}$ is a 
sequence $\langle I_t:t \in L \rangle = \langle I^{\gt}_t:
t \in L^{\gt} \rangle = \langle I_t[{\gt}]:t \in L[{\gt}]\rangle$  such that:
\mn
\begin{enumerate}
\item[$(a)$]  $L$ is a linear order (or partial, 
it does not really matter); 
but we may write $x \in {\gt}$
instead of $x \in L$ and $x <_{\gt} y$ instead of $x <_L y$
  % 2021-01-18 13:33 )
\sn
\item[$(b)$] $I_t$ is an ideal of % 2020-11-30 11:31 countable 
subsets of $L_t  = % 2021-01-18 13:35 
\{s:L \models s < t\}$, % 2021-02-02 11:58 
(but see \ref{ad.3}(4)(b)).
\end{enumerate}
\mn
% 2020-12-24 13:27 %
% 1A) Let $ \partial ( \mathfrak{t} )= \sup \{|A|^+:
% A \in I^ \mathfrak{t} _t $ for some $ t \in L^\mathfrak{t} \} $
% % Let $c \ell(I^{\gt}_t) = I^{c \ell,{\gt}}_t = 
% % \{A:[A]^{\le \aleph_0} \subseteq I^{\gt}_t\}$.

\noindent
2) Let ${\gt}$ be an FSI-template.
\mn
\begin{enumerate}
\item[$(c)$]  We say $\bar K = \langle K_t:t \in L^{\gt} \rangle$ is
a ${\mathfrak t}$-memory choice 
(or $ (\mathfrak{t} , \bar{ K } ) $ is an FSI-template)
\If \, 
\sn
\begin{enumerate}
\item[$(i)$]  $K_t \in I^{\gt}_t$ % 2020-11-30 11988:35 is countable 
\sn
\item[$(ii)$]  $s \in K_t \Rightarrow K_s \subseteq K_t$.
\end{enumerate}
\sn
\item[$(d)$]  We say $L \subseteq L^{\gt}$ is $\bar K$-closed if
$t \in L \Rightarrow K_t \subseteq L$
\sn
\item[$(e)$]  for $\bar K$ a ${\gt}$-memory choice and $L \subseteq
L^{\gt}$ which is $\bar K$-closed we say $\bar K' = \bar K
\restriction L$ if $\Dom(\bar K ' ) = L$ and
$K  ' _t$ is $K_t$ for $t \in L$, (it is a 
$({\gt} \restriction L)$-memory choice, see part (5)).
% 2020-11-30 19:51 
\sn 
\item[$(f)$]  We say that $ A$  is $ \bar{ K } $-countable 
(or, pedantically $ (\mathfrak{t}, \bar{ K})$-countable)
\when  \,  % 2021-02-02 11:59 
$(A = \emptyset$  or) there are $ t_n \in L^ \mathfrak{t} $ 
and $ \bar{ K } $-closed,  $ A_n \in I^ \mathfrak{t} _ {t_n}$
for $ n < \omega $  such that 
$ A= \cup \{A_n \cup \{ t_n \}: n < \omega  % 2020-12-04 15:35 
\} $. 
We define similarly $ \bar{ K } $-finite or 
 $ (\bar{ K } , < \partial )$-finite  % 2021-02-02 12:00 small
 )  % 2020-12-24 13:30 K,
 for any (infinite) $ \partial $
\sn
   % 2020-12-04 15:31 
\item[(g)] Let $ K^{\dagger}_t $ be $ K_t \cup \{ t\}  $
\sn
\item[(h)]  We let $ \partial ( \mathfrak{t} ) = \sup \{
 |A|^+  + {\aleph_0}   % 2020-12-24 13:29 
   : A \in I_t $ for some $ t \in L \} $  
 and $ \partial ( \mathfrak{t}, \bar{ K })  = \partial  ( \mathfrak{t} )$.  % 2020-12-04 15:34 
 Let $ \partial ( \bar{ K }) = \sup \{ | \bar{ K }_t |: t \in L^ \mathfrak{t}  \} $.
\end{enumerate}
\mn
3) For an FSI-template ${\gt}$ and ${\gt}$-memory choice 
$\bar K$
and $\bar K$-closed $L \subseteq L^{\gt}$ we define 
$\Dp_{\gt}(L,\bar K)$, 
the ${\gt}$-depth (or $({\gt},\bar K)$-depth)
of $L$ by defining by induction on the 
ordinal $\zeta$ when $\Dp_{\gt}(L,\bar K) \le \zeta$.
\medskip

\noindent
\underline{For $\zeta=0$}:  $\Dp_{\gt}(L,\bar K) \le \zeta$ 
when $L = \emptyset$.
\medskip

\noindent
\underline{For $\zeta$ a successor ordinal}:  
$\Dp_{\gt}(L,\bar K) \le \zeta$ iff:
\mn
\begin{enumerate}
\item[$(a)$]  there is $L^*$ such that: $L^* \subseteq L,|L^*| \le
1, (\forall t \in L)
(\forall A \in I^{\gt}_t)(A \cap L^* = \emptyset)$ hence 
$L \backslash L^*$ is $\bar K$-closed and $\Dp_{\gt}
(L \backslash L^*,\bar K) < \zeta$ and for every $t \in L^*$ 
we have:
\sn
\begin{enumerate}
\item[$\boxtimes_{t,L}$]  $L \backslash L^* \in 
I^{\gt}_t$ and \footnote{
   we can use less, it seems not needed
at the moment.  We can go deeper to names of depth $\le \varepsilon$
inductively on 
$\varepsilon < \omega_1$, as in \cite[\S3]{Sh:619}, or in a more
particular way to make the point that is used here true, and/or make
$I^{\gt}_t$ only closed under unions (but not subsets), etc. 
Note that e.g. $\Lim_{\gt}(\bar{\bbQ})$ is well defined when
      $L^{\gt}$ is well ordered.} 
it is $\bar K$-closed.
\end{enumerate}
% 2020-11-30 14:57 \end{enumerate}
\end{enumerate} 
\medskip

\noindent
\underline{For $\zeta > 0$ a limit ordinal}:  
$\Dp_{\gt}(L,\bar K) \le \zeta$ iff:
\mn
% 2020-11-30 14:57 \begin{enumerate}
\begin{enumerate} % 2020-12-24 13:46 
\item[$(b)$]   there is a 
directed partial order $M$ and a sequence $\langle
L_a:a \in M \rangle$ with union $L$ such that the sequence is
increasing, i.e., $M \models \lqq  a \le b  % 2020-12-24 13:35 lqq
\Rightarrow L_a \subseteq L_b"$, each $L_b$ is $\bar K$-closed, 
$(\forall b \in M)(\zeta > \Dp_{\gt}(L_b,\bar K)$) 
and $t \in L \wedge A \in I_t \wedge A \subseteq L \Rightarrow 
(\exists a \in M) \, A \subseteq L_a$.
% 2020-11-30 11:44 \end{enumerate}
% 2020-11-30 11:44 
\mn

\noindent 
So $\Dp_{\gt}(L,\bar K) = \zeta$ iff $\Dp_{\gt}(L,\bar K) 
\le  % 2021-01-18 13:42 \ge
\zeta \wedge 
(\forall \xi < \zeta) \Dp_{\gt}(L,\bar K) \nleq \xi$ 
% 2020-12-04 15:36 and 
\end{enumerate} % 2020-12-24 13:46 
\noindent 
3A) 
$\Dp_{\gt}(L,\bar K) = \infty$ iff ($\forall$ ordinal $\zeta$) 
[$\Dp_{\gt}(L,\bar K) \nleq \zeta]$.
% 2020-12-24 13:46 \end{enumerate} 

\noindent
4) We say $\bar K$ is a smooth ${\gt}$-memory choice or $(\gt,\bar K)$
is smooth \If \, $\Dp_{\gt}(L^{\gt},\bar K) < \infty$ 
and $\bar K$ a ${\gt}$-memory choice 
(and $ \mathfrak{t} $ is an FSI-template). % 2020-11-30 12:42 

\noindent
5) If $\bar K$ is omitted we mean it is the trivial $\bar K$, that is
$K_t = \emptyset$ for $t \in L^{\gt}$.
We say ${\mathfrak t}$ is smooth if the trivial $\bar K$ is a smooth
${\gt}$-memory choice.  For $L \subseteq L^{\gt}$ let ${\gt}
\restriction L = \langle I_t \cap {\cP}(L):t \in L \rangle$. 

\noindent
6) Let $L_1 \le_{\gt} L_2$ mean $L_1 \subseteq L_2 \subseteq
L^{\gt}$ and $t \in L_1 \wedge  A \in I^{\gt}_t \Rightarrow A 
\cap L_2 \subseteq L_1$.
% 2020-11-30 15:00 \end{enumerate} 
\end{definition}

\begin{definition}
\label{ad.2}
Let ${\gt} = \langle I_t:t \in L^{\gt} \rangle$ be a 
FSI-template and $\bar K$ a ${\gt}$-memory choice.  

\noindent
1) We say $\bar L$ is a $({\gt},\bar K)$-representation of $L$ (or
$(\gt,\bar K)-0$-representation of $L$) \If \,:  % 2020-12-24 13:55 
\mn
\begin{enumerate}
\item[$(a)$]  $L \subseteq L^{\gt}$ is $\bar K$-closed
\sn
\item[$(b)$]  $\bar L = \langle L_a:a \in M \rangle$
\sn
\item[$(c)$]  $M$ is a directed partial order
\sn
\item[$(d)$]  $\bar L$ is 
increasing, that is $a <_M b \Rightarrow L_a \subseteq L_b$
\sn
\item[$(e)$]  $L = \bigcup\limits_{a \in M} L_a$
\sn
\item[$(f)$]  each $L_a$ is $\bar K$-closed
\sn
\item[$(g)$]  if $t \in L,A \in I^{\gt}_t,A \subseteq L$ then
$(\exists a\in M)(A \subseteq L_a)$
% 2020-11-30 12:43 \sn
% 2020-11-30 12:43 \item[$(h)$]  clause (c) of \ref{ad.1}(3) holds.
\end{enumerate}
\mn
2) We say $L^*$ is a $({\gt},\bar K)-^*$representation or a
$({\gt},\bar K)-1$-representation + % 2020-12-24 13:56  of $L$)  
 of $L$ \If \,:
\mn
\begin{enumerate}
\item[$(a)$]  $L \subseteq L^{\gt}$ is $\bar K$-closed
\sn
\item[$(b)$]  $L^* \subseteq L,L^*$ a singleton
\sn
\item[$(c)$]  if $t \in L$ and $A \in I^{\gt}_t$ then $A \cap L^* =
\emptyset$ (so $(L \backslash L^*) \le_{\gt} L$, see Definition
\ref{ad.1}(6))
\sn
\item[$(d)$]  if $t \in L^*$ then $L \backslash L^* \in I^{\gt}_t$.
\end{enumerate}
\end{definition}

\begin{claim}
\label{ad.3}
Let ${\gt}$ be an $\FSI$-template and $\bar K$ a ${\gt}$-memory
choice. 

\noindent
0) The family of $\bar K$-closed sets is closed under (arbitrary) unions and
intersections.  Also if $L \subseteq L^{\gt}$ then 
$L \cup \bigcup\{K_t:t \in L\}$ is $\bar K$-closed. 

\noindent
1) If $L_2 \subseteq L^{\gt}$ is $\bar K$-closed and $L_1$ is an initial
segment of $L_2$, \then \, $L_1$ is $\bar K$-closed. 

\noindent
2) If $L_1 \subseteq L_2 \subseteq L^{\gt}$ are $\bar K$-closed \then \,
\mn
\begin{enumerate}
\item[$(\alpha)$]  $\Dp_{\gt}(L_1,\bar K) \le \Dp_{\gt}(L_2,\bar K)$, 
moreover
\sn
\item[$(\beta)$]  $(\exists t \in L_2)[L_1 \in I^{\gt}_t]$
implies that $\Dp_{\gt}(L_1,\bar K) < \Dp_{\gt}(L_2,\bar K)$ or 
both are $\infty$.
\end{enumerate}
\mn
3)  If $L_1 \subseteq L_2 \subseteq L^{\gt}$ are $\bar K$-closed
\then \, ${\gt} \restriction L_2$ is an $\FSI$-template, $L_1$ is
$({\gt} \restriction L_2)$-closed and $\Dp_{\gt \restriction L_2}
(L_1,\bar K \restriction L_2) = \Dp_{\gt}(L_1,\bar K)$.

\noindent
4)  If $(\gt,\bar K)$ is smooth and $A \in I_t,t \in L^{\gt}$ \then \,:
\mn
\begin{enumerate}
\item[$(a)$]  there is a $\bar K$-closed $B \in I_t$ such that $A
  \subseteq B$
\sn
\item[$(b)$]  if $s \in L_1 \in I_t$ and $L_2 \in I_s$ then $L_1 \cup
  L_2 \in I_t$.
\end{enumerate}
\end{claim}

\begin{PROOF}{\ref{ad.3}}
0), 1) Trivial - read the definitions. 

\noindent
2) We prove by induction on the ordinal $\zeta$ that
\mn
\begin{enumerate}
\item[$(*)_\zeta(\alpha)$]  if $\Dp_{\gt}(L_2,\bar K) = \zeta$ 
(and $L_1  \subseteq % 2020-12-04 15:39 
L_2$ are $\bar K$-closed subsets of $L^{\gt}$) 
\then \, $\Dp_{\gt}(L_1,\bar K) \le \zeta$
\sn
\item[$(\beta)$]  if in addition $(\exists t \in L_2)(L_1 \in
  I^{\gt}_t)$ \then \, $\Dp_{\gt}(L_1,\bar K) < \zeta$.
\end{enumerate}
\mn
So assume $\Dp_{\gt}(L_2,\bar K) = \zeta$, so $\Dp_{\gt}(L_2,\bar K) 
\ngeq \zeta +1$ hence one of the following cases occurs.
\bigskip

\noindent
\underline{First Case}:  $\zeta = 0$.  

Trivial; note that clause $(\beta)$ is empty.
\bigskip

\noindent
\underline{Second Case}:  $\zeta$ is a successor, hence
$L_2$ has a $({\gt},\bar K)$-$^*$representation $L^*$ such that 
$\Dp_{\gt}(L_2 \backslash L^*,\bar K) < \zeta$; see Definition \ref{ad.2}(2).

Let $L^-_2 =: L_2 \backslash L^*$; if $L_1 \subseteq L^-_2$ then by the
induction hypothesis $\Dp_{\gt}(L_1,\bar K) \le \Dp_{\gt}(L^-_2,\bar K) <
\zeta$, so assume $L_1 \nsubseteq L^-_2$ and so only clause $(\alpha)$ is
relevant.  Now letting
$L^-_1 = L_1 \backslash L^*$ we have $[L^-_1,L^-_2$ are $\bar K$-closed]
$\wedge L^-_1 \subseteq L^-_2$ and $\Dp_{\gt}(L^-_2,\bar K) < \zeta$ 
hence $\Dp_{\gt}(L^-_1,\bar K) < \zeta$ by the induction hypothesis.
Let $L^*_1 = L_1 \cap L^*$, so $L^*_1 \subseteq L_1,L_1$ is $\bar K$-closed,
$L_1 \backslash L^*_1 = (L_2 \backslash L^*_2) \cap L_1$ is $\bar K$-closed,
$\Dp_{\gt}(L_1 \backslash L^*_1,\bar K) = \Dp_{\gt}(L^-_1,\bar K) 
< \zeta$ and necessarily $L^*_1$ has exactly one element.  
Also easily: $t \in L^*_1$ implies $L^-_1 \in
I^{\gt}_t$ so $L^*_1$ is a $({\gt},\bar K)-^*$representation 
of $L_1$.  So clearly $\Dp_{\gt}(L_1,\bar K) \le \Dp_{\gt}
(L^-_1,\bar K)+1 \le \zeta$.
\bigskip

\noindent
\underline{Third Case}:  $\zeta$ is limit and
$\langle L_a:a \in M \rangle$ is a
$({\gt},\bar K)$-representation of $L_2$ such that $a \in M
\Rightarrow \Dp_{\gt}(L_a,\bar K) < \zeta$.

Let $L^2_a =: L_a$ and $L_a^1 =: L_a \cap L_1$, so 
$\langle L^1_a:a \in M \rangle$ is increasing, 
$\bigcup\limits_{a \in M} L_a^1 = L_1$ and each $L^1_a$ is 
$\bar K$-closed (as $L^2_a,L_1$ are $\bar K$-closed, see part (0)) 
and easily % 2020-12-24 % 2020-12-24 14:03 
 $t \in L_1 \wedge A \in I^{\gt}_t \wedge
A \subseteq L_1 \Rightarrow (\exists a \in M)(A \subseteq L^2_a 
\cap L_1 = L^1_a)$.   Also by the definition of $ \Dp$ 
 at limit ordinals       % 2020-12-24 14:04 induction hypothesis,
$b \in M \Rightarrow 
\Dp_{\gt}(L^2_b,\bar K) < \zeta$. 
Hence by the induction hypothesis 
$ \Dp_\mathfrak{t} (L^1_b  % 2021-01-18 13:49 t
, \bar{ K } ) < \zeta $.  % 2020-12-24 14:06 
By 
the last two sentences (and Definition \ref{ad.1}) we get 
$\Dp_{\gt}(L_1,\bar K) \le \zeta$, as required in clause
$(\alpha)$.  For clause $(\beta)$ we know that there is $t \in L_2$
such that $L_1 \in I^{\gt}_t$, hence by clause (g) of  % 2021-01-18 13:49 f++>g
Definition \ref{ad.2}(1)) for some $b \in M$ we have 
$L_1 \subseteq L_b$ and we can use the induction hypothesis on $\zeta$
for $L_1,L_b$. 

\noindent
3) Easy. 

\noindent
4) By induction on the depth $ \zeta $.  % 2020-12-04 15:40 
The case $ \zeta = 0 $ is trivial; % 2020-12-24 14:11 the  
and  the   case 
$ \zeta $ is a limit ordinal is easy. Lastly  % 2020-12-24 14:10 
for the % 2020-12-24 14:10 
% 2020-12-04 15:41 the point is the 
successor case of \ref{ad.1}(3)   % 2020-11-30 12:45 
recall $ \boxplus _{t,L}$  there. % 2020-12-04 15:42 
\end{PROOF}
\bigskip

\begin{claim}
\label{ad.4}
1) If for $\ell = 1,2$ we have $\bar L^\ell$ is a 
$({\gt},\bar K)$-representation of $L$ and
$\bar L^\ell = \langle L^\ell_a:a \in M_\ell \rangle$ and $M = 
M_1 \times M_2$ \then \, $\bar L = \langle L_a \cap L_b:(a,b) \in M 
\rangle$ is a $({\gt},\bar K)$-representation of $L$. 

\noindent
2) If $L^*_\ell$ is a $({\gt},\bar K)-^*$representation of 
$L$ for $\ell = 1,2$ \then \, $L^*_1 = L^*_2$.

\noindent 
3) 
If $ A $ is $ (\bar{ K}, % 2020-12-24 14:17 
  < \partial )$-countable % 2021-02-02 12:01 small 
  \then  \, it is $ \bar{ K } $-closed. % 2020-12-04 15:44 

\noindent 
4) If $ L \subseteq L^ \mathfrak{t} $ is $ \bar{ K}$-closed and 
$ L_1 \subseteq L $ has cardinality  $ < \partial$ % 2020-12-24 14:19  ( \mathfrak{t} )$
\then \, % 2020-12-24 14:20 for 
some $ (\bar{ K } , \partial )$-finite  % 2021-02-02 12:01 small % 2020-12-24 14:20 -closed 
set $ L_2 \subseteq L$
   % 2020-12-24 14:20 and 
includes $ L_1$.  % 2020-12-04 15:52 

\noindent 
5) If $ (\mathfrak{t} , \bar{ K } )$ is an FSI-template then so are 
$ \mathfrak{t} ', \mathfrak{t} ''$  where % 2020-12-24 14:50 
% 2020-12-04 15:54 vacirah   % 2020-12-04 17:28 
$ \mathfrak{t} ', \mathfrak{t}'' $ are FSI-templates satisfying
\begin{enumerate} 
\item[(a)] $ L^{\mathfrak{t} '} = L^{\mathfrak{t} ''}= L^{\mathfrak{t} }$
\item[(b)] for $ t \in L^\mathfrak{t} $  let
     $ I^{\mathfrak{t}'}_t=  \{ A \in I^ \mathfrak{t} _t : 
    A $ is $ \bar{ K}$-countable$ \} $   % 2020-12-24 14:52 |A| < \partial ( \mathfrak{t} \} $
\item[(c)] for $ t \in L^\mathfrak{t} $  let
    $ I^{\mathfrak{t}''}_t =  \{ A \subseteq L^\mathfrak{t} :$ 
     the set  $ \{ B \subseteq A : B \in I^\mathfrak{t} _t $ is 
     $ \bar{ K } $-closed$ \} $ is cofinal in $ [A]^{< \partial ( \mathfrak{t}) } \} $ 
\end{enumerate} 
\end{claim}

\begin{PROOF}{\ref{ad.4}}
1) Straightforward, e.g. if $ t \in L^ \mathfrak{t} , A \in I_t $  % 2020-12-04 15:45 
and $ A \subseteq L $  then for $ {\ell} = 1,2 $ we can choose
$ a_{\ell} \in M_ {\ell} $ such that $ A \subseteq L^{\ell} _{a_{\ell} }$ 
and $ t \in L_{a_{\ell} }$. Clearly $ A \cup \{t % 2020-12-24 14:54  
   % 2021-01-18 13:53 y 
\} \subseteq L^1_{a_1} \cap L^2_{a_2} $.  

\noindent
2)-5)  Easy, too.  
\end{PROOF}

\begin{discussion}
\label{ad.4a} 
This discussion % 2020-12-24 14:55 
is from the old version, so some ``we may" 
are  actually  done in the new version.

\noindent 
1) Our next aim is to define iteration for any $\bar K$-smooth FSI-template
${\gt}$; for this we define and prove the relevant things; of course,
by induction on the depth.
In the following Definition \ref{ad.5}, in clause (A)(a), we
avoid relying on \cite{Sh:630}; moreover the reader may 
consider only the case $K_t = \emptyset$, omit $\name \eta_t$ and have 
$\name{\bbQ}_{t,\bar \varphi'_t}$ be the dominating real 
forcing = Hechler forcing.  

\noindent
2)  We may more generally than here allow
$\name \eta_t$ to be e.g. a sequence of ordinals,
and members of $\name{\bbQ}_{t,\varphi,\name \eta_t}$ be 
$\subseteq {\cH}_{< \aleph_1}(\Ord)$, and even $K_t$ large 
but increasing $L$, we need more ``information" from
$\name \eta_t \restriction \Lim_{\gt} (\bar{\bbQ} \restriction L)$.  
We may 
require more by % 2020-11-30 12:47 
changing  to: $\name{\bbQ}_t$ is a definition of 
nep c.c.c. forcing (\cite{Sh:630}) or just 
``Souslin c.c.c. forcing (= snep)" or just absolute 
enough c.c.c. forcing notion.
All those cases do not make real problems 
 (but when the parameter
$\name \eta_t$ have length $\ge \kappa$ 
(or just has no bound $ < \kappa $) % 2020-11-30 12:49 
it is changed in  % 2020-11-30 12:57 
the ultra-power! i.e. $\mathbf j(\name \eta_t)$ has
length $>$ length of $\name \eta_t$). 

\noindent
3) If we restrict ourselves to $\sigma$-centered forcing notions
(which is quite reasonable) probably we can in Definition \ref{ad.1}(3)(a)
omit $\boxtimes_{t,L}$ if in Definition \ref{ad.5} below in (A)(b)
second case we add that $t \in L^* \Rightarrow p \restriction (L
\backslash L^*)$ forces a value to $\name f_t(p(t))$ where 
$\name f_t:\name{\bbQ}_t \rightarrow \omega$ witnesses
$\sigma$-centerness and 
is absolute enough (or just assume $\bbQ_t \subseteq \omega \times 
\bbQ'_t,f_t(p(t))$ is the first coordinate).  More carefully probably we 
can do this with $\sigma$-linked instead $\sigma$-centered.
\end{discussion}

\begin{dc}
\label{ad.5}
Let ${\gt}$ be an FSI-template and $\bar K = \langle K_t:t \in 
L^{\gt} \rangle$ be a smooth ${\gt}$-memory choice.

By induction on the ordinal $\zeta$ we shall define and prove:
\mn
\begin{enumerate}
\item[$(A)$]  [Def] $\quad$ for $L \subseteq L^{\gt}$ which is
$\bar K$-closed of $({\gt},\bar K)$-depth $\le \zeta$ we define 
\sn
\begin{enumerate}
\item[$(a)$]  when $\bar{\bbQ} = \langle 
\name{\bbQ}_{t,\bar\varphi_t,\name \eta_t}:t \in L \rangle$
is a $({\gt},\bar K)$-iteration of def-c.c.c. forcing notions, but we
can let $\name \eta_t$ code $\bar \varphi_t$, 
say as $ \bar{ \varphi }= \name{ \eta } (0) \in 
{\mathscr H} ( {\aleph_0} ) $; 
so we may
omit $\bar \varphi_t$; note that ``def. - c.c.c." is defined below
\sn
\item[$(b)$]  $\Lim_{\gt}(\bar{\bbQ})$ for $\bar{\bbQ}$ 
as in (A)(a),    % 2020-12-24 16:03 
pedantically we should write 
   $\Lim_{\gt, \bar{ K } }(\bar{\bbQ})$
\item[(c)] $ \name {\bar\nu}$ is the 
 sequence of generics  of $ \bar{ \mathbb{Q} } $
\item[(d)] $ \bar{ u } = \bar{ u } [\bar{ \mathbb{Q} } ]$
   is the parameters domain sequence of $ \bar{ \mathbb{Q} } $
\item[(e)] the class $ \mathbf{Q} = \mathbf{Q}_{\fsi} $  
  of  fsi-templates as well  as  some 
related classes
\item[(f)] $  \partial ( \mathbf{q} ), \partial ^-(\mathbf{q} )$   % 2020-12-01 13:54 
for $ \mathbf{q} \in \mathbf{Q} $. 
% \item[(g)] Let $ K^{\dagger}_t = K_t \cup \{ t  \} $
\item[(g)] $ \mathbf{q} $ is $ \lim \langle \mathbf{q} _ a : a \in M \rangle $ where $ M $  is a directed partial order
and $ \mathbf{q} _ a   % 2021-01-27 11:09 
\in \mathbf{Q} $ is increasing with $ a $.
\end{enumerate}
\sn
\item[$(B)$]  [Claim] $\quad$ for $L_1 \subseteq L_2 \subseteq
L^{\gt}$ which are $\bar K$-closed of $({\gt},\bar K)$-depth $\le \zeta$ 
and $({\gt},\bar K)$-iteration of def-c.c.c. forcing notions 
$\bar{\bbQ} = \langle \name{\bbQ}_{t,\name{\bar \varphi}_t, % 2020-11-30 13:11 
\name \eta_t  } :t \in L_2 \rangle$ we prove:
\sn
\begin{enumerate}
\item[$(a)$]  $\bar{\bbQ} \restriction L_1$ is a $({\gt},
\bar K \restriction L_1)$-iteration of def-c.c.c. forcing notions
\sn
\item[$(b)$]  $\Lim_{\gt}(\bar{\bbQ} \restriction 
L_1) \subseteq \Lim_{\gt}(\bar{\bbQ})$ as quasi orders
\sn
\item[$(c)$]   if $L_1 \le_{\gt} L_2$ (see Definition \ref{ad.1}(6)) and
$p \in \Lim_{\gt}(\bar{\bbQ})$, \then \, 
$p \restriction L_1 \in \Lim_{\gt}(\bar{\bbQ} \restriction L_1)$ and
$\Lim_{\gt}(\bar{\bbQ}) \models ``p \restriction L_1 \le p"$ 
\sn
\item[$(d)$]  if $L_1 \le_{\gt} L_2$ and $p \in \Lim_{\gt}(\bar{\bbQ})$ and
$\Lim_{\gt}(\bar{\bbQ} \restriction L_1) \models ``(p \restriction L_1) 
\le q"$ \then \, $q \cup (p \restriction (L_2 \backslash
L_1))$ is a lub of $\{p,q\}$ in $\Lim_{\gt} (\bar{\bbQ})$; hence 
$\Lim_{\gt}(\bar{\bbQ} \restriction L_1) \lessdot \Lim_{\gt}
(\bar{\bbQ})$, (used in the proof of clause (B)(j))

% 2020-12-03 16:27 BDOQ 2020-11-30 13:11 
\sn
\item[$(e)$]  $\Lim_{\gt}(\bar{\bbQ} \restriction L_1) \lessdot 
\Lim_{\gt}(\bar{\bbQ})$, that \footnote{here we do not assume
$L_1 \le_{\gt} L_2$,} is
\sn
\item[${{}}$]  $\qquad (i) \quad 
p \in \Lim_{\gt}(\bar{\bbQ} \restriction L_1)
\Rightarrow p \in \Lim_{\gt}(\bar{\bbQ})$
\sn
\item[${{}}$]  $\qquad (ii) \quad 
\Lim_{\gt}(\bar{\bbQ} \restriction L_1) \models p \le q 
\Rightarrow \Lim_{\gt} (\bar{\bbQ}) \models p \le q$
\sn
\item[${{}}$]  $\qquad (iii) \quad$ if ${\cI} \subseteq
\Lim_{\gt}(\bar{\bbQ} \restriction L_1)$ is predense in 
$\Lim_{\gt}(\bar{\bbQ} \restriction L_1)$, \then \,

\hskip35pt  ${\cI}$ is predense in $\Lim_ \mathfrak{t} (\bar{\bbQ})$
\sn
\item[${{}}$]  $\qquad (iv) \quad$ if $p,q \in \Lim_{\gt}(\bar{\bbQ})$ are 
incompatible in $\Lim_{\gt}(\bar{\bbQ} \rest L_1)$ then they 

\hskip35pt are incompatible in $\Lim_ \mathfrak{t} (\bar{\bbQ})$)
\sn
\item[$(f)$]   assume $L_0 \subseteq L_2$ is $\bar K$-closed,
$L = L_0 \cap L_1$; if $p \in \Lim_{\gt}(\bar{\bbQ} \restriction L_0)$ 
and $q \in \Lim_{\gt}(\bar{\bbQ} \restriction L)$
satisfies $(\forall r \in \Lim_{\gt}(\bar{\bbQ} \restriction L))
[q \le r \rightarrow p,r$ are compatible in 
$\Lim_{\gt}(\bar{\bbQ} \restriction L_0)]$ 
\then \, $(\forall r \in \Lim_{\gt}(\bar{\bbQ} \rest L_1))
[q \le r \rightarrow p,r$ are compatible in 
$\Lim_{\gt}(\bar{\bbQ} \restriction L_  2  % 2021-01-18 13:55 1
)]$  % 2020-12-24 16:10 2--> 1
\sn

\small 
[explanation: this means that if $q$ forces for
$\Vdash_{\Lim_{\gt}(\bar{\bbQ} \restriction L_0)}$ that
$p \in \Lim_{\gt}(\bar{\bbQ} \restriction L_0)/\Lim_{\gt}(\bar{\bbQ} 
\restriction L)$ \then \, $q$ forces for 
$\Vdash_{\Lim_{\gt}(\bar{\bbQ} \restriction L_1)}$ that 
$p \in \Lim_{\gt}(\bar{\bbQ})/\Lim_{\gt} (\bar{\bbQ} \restriction L_1)$.]
\normalsize
\sn
\item[$(g)$]  if $\langle L_a:a \in M_1 \rangle$ is a
$({\gt},\bar K)-$representation of $L_1$ \then \,
$\Lim_{\gt}(\bar{\bbQ} \restriction L_1) = \bigcup\limits_{a \in M_1} 
\Lim_{\gt}(\bar{\bbQ} \restriction L^1_a)$
\sn
\item[$(h)$]   if $L^*$ is a $({\gt},\bar K)$-$^*$representation
of $L_1$  and $ L^* = L \cup \{ t \} $,  % 2021-01-18 15:22 
\then \, 
$ \Lim _ \mathfrak{t} (\bar{ \mathbb{Q} }\upharpoonright L_1)= 
  \Lim _ \mathfrak{t} (\bar{ \mathbb{Q} } \upharpoonright (L_1 \setminus L^*) * \name{ \mathbb{Q}} _ {t, \name{\eta }_ t}
  ) $ % 2021-02-10 04:54 
% $\Lim_{\gt} (\bar{\bbQ} \restriction L_1)$ is 
% as defined in (A)(b) of our definition below, second case, from $L^*$
\sn
\item[$(i)$]  $(\alpha) \quad $if $p_1,p_2 \in \Lim_{\gt}(\bar{\bbQ})$ and
$t \in \Dom(p_1) \cap \Dom(p_2) \Rightarrow p_1(t) =$ 

\hskip32pt  $p_2(t)$, \then \, $q=p_1 \cup p_2$ 
(i.e. $p_1 \cup (p_2 \backslash (\Dom(p_1)))$ belongs to 

\hskip32pt  $\Lim_{\gt}(\bar {\mathbb{Q} })$ % 2020-11-30 13:14 
and is a l.u.b. of $p_1,p_2$  in it % 2020-11-30 13:15 
\sn
\item[${{}}$]  $(\beta) \quad p \in \Lim_{\gt}(\bar{\bbQ})$
\Iff \, $p$ is a function with domain a finite subset of 

\hskip32pt $L_2$ such that for every $t \in \Dom(p)$ for some 
$A \in I^{\gt}_t,A$ is 

\hskip32pt  $\bar K$-closed and $K_t \subseteq A$ and
$\Vdash_{\Lim_{\gt}(\bar{\bbQ} \restriction A)} ``p(t) \in
\bbQ_{t,\name \eta_t}"$.  

\hskip32pt  [So if $p \in \Lim_{\gt}(\bar{\bbQ})$ then for some
$ \bar{ K }$-countable (even $ \bar{ K } $-finite, 

\hskip32pt  see \ref{ad.1})(2)(f)),  % 2020-11-30 13:16 
$L \subseteq L_2$ 
% 2020-12-04 17:58 \hskip32pt  
we have $p \in \Lim_{\gt}(\bar{\bbQ} \restriction L)]$
\sn
\item[${{}}$]  $(\gamma) \quad \Lim_ \mathfrak{t} % 2020-11-30 13:17 t
(\bar{\bbQ}) \models 
p \le q$ \Iff \, $p,q \in \Lim_t(\bar{\bbQ})$ and for every 
$t \in \Dom(p)$ 

\hskip32pt   we have $t \in \Dom(q)$ and for some $\bar K$-closed
$A \in I^{\gt}_t$ we have

\hskip32pt   $q \restriction A \in \Lim_{\gt}
(\bar{\bbQ} \restriction A)$ and $q \Vdash_{\Lim_{\gt}
(\bar{\bbQ} \restriction A)} ``p(t) \le q(t)$ 

\hskip32pt   in $\bbQ_{t,\name \eta_t}$ (as interpreted in 
$\mathbf V^{\Lim_{\gt}(\bar{\bbQ} \restriction A)}$ of course)"
\sn
\item[$(j)$]  

\begin{enumerate} 
\item[$(\alpha )$] $\Lim_{\gt}(\bar{\bbQ})$ is a c.c.c. forcing notion 
  % 2020-12-24 16:18 and 
  \item[$(\beta )$] 
$\Lim_{\gt}(\bar{\bbQ}) = \cup\{\Lim_{\gt}(\bar{\bbQ} \restriction L):
L \subseteq L_2 $   % 2020-12-24 16:19 
% 2020-11-30 15:40 \in [L_2]^{\le \aleph_0}$ 
is $\bar K$-finite% 2020-11-30 13:19 closed
$\}$
\end{enumerate} 
\sn
\item[$(k)$]  
\begin{enumerate}
\item[($ \alpha $)] $\Lim_{\gt}(\bar{\bbQ})$ has cardinality
$\le |L_2|^{\aleph_0}  + \partial ^- (\mathbf{q} ) $   % 2020-12-01 14:47 \theta 
\item[($ \beta $)] for every $ \Lim _ \mathfrak{t} (\bar{ \mathbb{Q} })$-name
$ \name{ \rho }$  of a real there is a $ \bar{ K } $-countable set $ L $ 
such that $ \name{ \rho } $ is a 
$ \Lim _ \mathfrak{t} (\bar{ \mathbb{Q}  } % 2021-01-18 14:06 
\upharpoonright L )$-name
% 2020-11-30 15:10 \item[(*)] 
% (here we use the assumption that $\name \eta_t$ 
% and members of $\name{\bbQ}_{t,\name \eta_t}$ 
% are reals; see definition in (A)(a)(i)+(i)) below).

\end{enumerate} 
\end{enumerate}
\end{enumerate}
\end{dc}

Let us carry the induction. 
\medskip

\noindent
\underline{Part (A)}: [Definition]

So assume $\Dp_{\gt}(L,\bar K) \le \zeta$.  If $\Dp_{\gt}(L) < \zeta$ we 
have already defined being $({\gt},\bar K)$-iteration and 
$\Lim_{\gt}(\bar{\bbQ} \restriction L)$, so assume $\Dp_{\gt}(L) = \zeta$.
\medskip

\noindent
\underline{Clause (A)(a)} For every $ t \in L^  \mathfrak{t} $ we have:
\mn
\begin{enumerate}
\item[$(i)$]   $\name \eta_t$ is a $\Lim_{\gt}(\bar{\bbQ}
\restriction K_t)$-name of a real (i.e. from ${}^\omega 2$, used as a
parameter) and $ u_t = \omega $
\underline{or} 
(see (A)(d))  % 2020-12-04 18:05 
% 2021-01-18 14:09
a function from   % 2021-02-10 04:57 
a set of ordinals % 2021-01-18 14:09 
$ u_t$  
 ($ u_t $  an object, not a name) % 2020-12-24 16:21 
into 
$ \{ 0,1\} $ or into $ {\mathscr H} ( {\aleph_0} ) $, 
(legal as $K_t $ is 
a  % 2021-01-18 14:00 
$ \bar{ K}$-closed subset of $ L $ and 
  % 2020-11-30 13:57 \subseteq L \wedge
  $ K_t \in I_t$ and 
$t \in L$ hence by \ref{ad.3}(2), clause $(\beta)$ we have
$\Dp_{\gt}(K_t,\bar K) < \Dp_{\gt}(K_t \cup \{t\},\bar K) \le
\Dp_{\gt}(L,\bar K) \le \zeta$ so $\Lim_{\gt}(\bar{\bbQ} \rest L_t)$ is a 
well defined forcing notion by the
induction hypothesis and \ref{ad.3}(2), clause $(\beta)$)
\sn
\item[$(ii)$]  $\bar \varphi_t$ is a pair of formulas 
which from   % 2021-02-10 04:58 from which % 2020-11-30 13:58 with 
the
parameters $\name \eta_t$ define  % 2020-11-30 13:58 ing 
in $\mathbf V^{\Lim_{\gt}(\bar{\bbQ} 
\restriction K_t)}$ a forcing notion
denoted by $\bbQ_{t,\bar \varphi_t,\name \eta_t}$ whose
set of elements is $\subseteq {\cH}(\aleph_1)$
or $ \subseteq {\mathscr H} _{{\aleph_1} }(u_t)$
\sn
\item[$(iii)$]  in $\mathbf V_1 = \mathbf V^{\Lim_{\gt}(\bar{\bbQ}
\restriction K_t)}$, if $\bbP' \lessdot \bbP''$ 
are c.c.c. forcing notions\footnote{
So the definition $ \bar{ \varphi }_t $ still defines a forcing notion; 
We may restrict ourselves to forcing notions which 
occur in our proof; but does not seem to matter for now.
} 
\then \, $\bbQ = \bbQ_{t,\bar \varphi_t,
\name \eta_t}$ as interpreted in $\mathbf V_2 = (\mathbf V_1)^{\bbP'}$ 
is a c.c.c. forcing notion there, and $\bbP' * 
\name{\bbQ}_{t,\bar \varphi_t,\name \eta_t}$ is a $\lessdot$-sub-forcing
of $\bbP'' * \name{\bbQ}_{t,\bar \varphi_t,
\name{ \eta } _t } $  % 2021-01-18 \eta_t}$
where 
$\name{\bbQ}_{t,\bar \varphi_t,\name \eta_t}$ 
mean as interpreted in $(\mathbf V^{\Lim_{\gt}
(\bar{\bbQ} \restriction K_t})^{\bbP'}$ or 
in $(\mathbf V^{\Lim_{\gt}(\bar{\bbQ} \restriction K_t)})^{\bbP''}$ 
respectively (i.e. $``p \le q"$,``$p,q$  
are  % 2021-01-18 14:01 
incompatible", 
``$\langle p_n:n < \omega \rangle$ is predense" (so the sequence is 
from the smaller universe) are preserved)
\sn
\item[$(iv)$]    assume that $\Lim_{\gt}(\bar{\bbQ} \upharpoonright K_t) 
\lessdot {\bbP}_0 \lessdot {\bbP}_\ell \lessdot {\bbP}_3$
are c.c.c. forcing notions for $\ell = 1,2 $
and 
$ {\bbP}_1 \cap {\bbP}_2 
= {\bbP}_0$.
% 2020-11-30 14:09 \end{enumerate}
\mn
Let $\name{\bbQ}_\ell$ be the $\bbP_\ell$-name of 
$\name{\bbQ}_{t, \name{ \eta } _t}$ as  % 2020-11-30 14:11 
interpreted in $\mathbf V^{\bbP_\ell}$.

If $(p_\ell,\name q_\ell) \in \bbP_\ell * \name{\bbQ}_\ell$ for
$\ell=0,1,2$ and $(p_0,\name q_0) \Vdash ``(p_\ell,q_\ell) \in
(\bbP_\ell * \name{\bbQ}_\ell)/(\bbP_0 * \name{\bbQ}_0)"$ for
$\ell=1,2$ 
and $p_3 \in \bbP$ is above $p'_1,p'_2$ 
\then \, there are $(p_\ell,\name q_\ell)  % 2021-01-18 14:11 '  '
\in \bbP_\ell *
\name{\bbQ}_\ell$ above $(p_\ell,q_\ell)$ for $\ell=0,1,2$ satisfying
$(p'_0,q'_0) \Vdash ``(p'_\ell,q'_\ell) \in (\bbP_\ell *
\name{\bbQ}_\ell)/(\bbP_0 * \name{\bbQ}_0)"$ for $\ell=1,2$ such that:
\mn
\begin{enumerate}
\item[$\bullet$] 
% 2020-12-25 14:46 if $p_3 \in \bbP$ is above $p'_1,p'_2$ then 
$p_3
  \Vdash_{\bbP_3} ``\name q_1,\name q_2$ are compatible in $\name{\bbQ}_3"$.
\end{enumerate} 
% \item[(c)] $\bar{ \name{\nu }}$ is the 
%  sequence if generics  of $ \bar{ \mathbb{Q} } $
% \item[(d)] $ \bar{ u } = \bar{ u } [\bar{ \mathbb{Q} } ]$
%   is the parameters domain sequence of $ \bar{ \mathbb{Q} } $
% \item[(e)] the class $ \mathbf{Q} $ of $ FSIT-s
% % \item[(f)] $  \theta ( \mathbf{q} $ for $ \mathbf{q} \in \mathbf{Q} $. \item[(*)] 6666666
% \item[(*)] 
% \item[(*)] 
% \item[(*)] 

\end{enumerate} 
\bigskip

\noindent
\underline{Clause (A)(b)}:
\bigskip

\noindent
\underline{First Case}:  $\zeta =0$.

Trivial.
\bigskip

\noindent
\underline{Second Case}:  $\zeta$ is a successor.

So let $L^*$ be a $({\gt},\bar K)$-$^*$representation of $L$.

Define $p \in \Lim_{\gt}(\bar{\bbQ})$ 
\Iff \,
$ p$ is a finite function, $\Dom
(p) \subseteq L,p \restriction (L \backslash L^*) \in 
\Lim_{\gt}(\bar{\bbQ} \restriction (L \backslash L^*))$ and if 
$t \in L^* \cap \Dom(p)$, then $p(t)$ is a $\Lim_{\gt}(\bar{\bbQ} 
\restriction (L \backslash L^*)$-name of a member of 
$\bbQ_{t,\bar\varphi_t,\name \eta_t}$ 
and the order is $\Lim_{\gt}(\bar{\bbQ}) \models p \le q$ \Iff \,
\mn
\begin{enumerate}
\item[$(i)$]  $\Lim_{\gt}(\bar{\bbQ} \restriction (L \backslash L^*)) \models
``(p \restriction (L \backslash L^*) \le 
(q \restriction (L \backslash L^*))"$ and
\sn
\item[$(ii)$]   if $t \in L^* \cap \Dom(p)$ then
for some $ \bar{ K } $-closed $ A \in I^ \mathfrak{t} _ t $ 
we have % 2021-01-18 14:19 
$q \restriction A   % 2021-01-18 14:19 _t
\Vdash_{\Lim_{\gt}(\bar{\bbQ} \restriction (L \backslash L^*))} 
``p(t) \le q(t)"$.
\end{enumerate}
\mn
Clearly $\Lim_{\gt}(\bar{\bbQ})$ is a quasi order.
But we should prove that Lim$_{\gt}(\bar{\bbQ})$ is well defined,
which means that the definition does not depend on the 
representation.  

So we prove
\mn
\begin{enumerate}
\item[$\boxtimes_1$]   if $\Dp_{\gt}(L,\bar K) = \zeta$ and for
$\ell=1,2$ we have $L^*_\ell$ is a 
$({\gt},\bar K)$-$^*$representation of $L$ with $\Dp_{\gt}
(L \backslash L^*_\ell,\bar K) < \zeta$ and $\bbQ^\ell$ is 
$\Lim_{\gt}(\bar{\bbQ} \restriction L)$ as defined by $L^*_\ell$
above, \then \, $\bbQ^1 = \bbQ^2$.
\end{enumerate}
\mn
This is immediate by Claim \ref{ad.4}(2) and 
the induction hypothesis clause (B)(h).
\bigskip

\noindent
\underline{Third Case}:  $\zeta$ limit.

So there are a directed partial 
order $M$ and $\bar L = \langle L_a:a \in M \rangle$ a 
$({\gt},\bar K)$-representation of $L$ such that
$a \in M \Rightarrow \Dp_{\gt}(L_a,\bar K) < \zeta$.  By
the induction hypothesis, $a \le_M b \Rightarrow L_a \subseteq L_b$ and
$\Lim_{\gt}(\bar{\bbQ} \restriction L_a) \subseteq
\Lim_{\gt}(\bar{\bbQ} \restriction L_b)$.

We let $\Lim_{\gt}(\bar{\bbQ} \restriction L) = \bigcup\limits_{a \in M}
\Lim_{\gt}(\bar{\bbQ} \restriction L_a)$, so we have to prove
\mn
\begin{enumerate}
\item[$\boxtimes_2$]  the choice % 2021-01-18 14:14 is 
of $\bar L$ is immaterial.
\end{enumerate}
\mn
So we just assume that for $\ell =1,2$ we have: $M_\ell$ is a directed
partial order, $\bar L^\ell = \langle L^\ell_a:a \in M_\ell 
\rangle,L^\ell_a \subseteq L,M_\ell \models a \le b 
\Rightarrow L^\ell_a \subseteq L^\ell_b$ and
$(\forall t \in L)(\forall A \in I_t)[A \subseteq L \rightarrow 
(\exists a \in M_\ell)(A \subseteq L^\ell_a)$ and 
$\Dp_{\gt}(L^\ell_a,\bar K) < \zeta$. 

We should prove that $\bigcup\limits_{a \in M_1} \Lim_{\gt}(\bar{\bbQ}
\rest L^1_a),\bigcup\limits_{a \in M_2} \Lim_{\gt}(\bar{\bbQ}
\rest L^2_a)$ are equal, as quasi orders of course.

Now let $M =: M_1 \times M_2$ with $(a_1,a_2) \le (b_1,b_2) 
\Leftrightarrow a_1 \le_{M_1} b_1 \wedge a_2 \le_{M_2} b_2$, 
clearly it   % 2020-12-24 16:24 
is a directed partial order.  We let
$L_{(a_1,a_2)} = L^1_{a_1} \cap L^2_{a_2}$, so clearly $L_{(a_1,a_2)}
\subseteq L^{\gt},\Dp_{\gt}(L_{(a_1,a_2)},\bar K) < \zeta$ and $(a_1,a_2) \le_M
(b_1,b_2) \Rightarrow L_{(a_1,a_2)} \subseteq L_{(b_1,b_2)}$ and
$\langle L_{(a_1,a_2)}:(a_1,a_2) \in M \rangle$ is a 
$({\gt},\bar K)$-representation of $L$ by Claim \ref{ad.4}(1). So by
transitivity of equality, it is enough to prove for $\ell =1,2$ that
$\bigcup\limits_{a \in M_\ell} \Lim_{\gt}(\bar{\bbQ} 
\restriction L^\ell_a),\bigcup\limits_{(a,b) \in M} \Lim_{\gt}
(\bar{\bbQ} \restriction L_{(a,b)})$ are equal as quasi orders.  
By the symmetry in the situation without loss of generality $\ell=1$.

Now for every $a \in M_1,\bar L= \langle L_{(a,b)}:b \in M_2 \rangle$ 
satisfies: $L^1_a \subseteq L,\Dp(L^1_a,\bar K) < \zeta,
% 2020-11-30 15:44 L^1_a,
L^1_a 
= \bigcup\limits_{b \in M_2} L_{(a,b)},b_1 \le_{M_2} b_2 \Rightarrow
L_{(a,b_1)} \subseteq L_{(a,b_2)}$.  
 
%  111
%   Also we know that $(\forall t \in
% L)(\forall A \in I^{\gt}_t)(\exists b \in M_2)(A \subseteq L
% \rightarrow A \subseteq L^2_b)$ so $(\forall t \in L^1_a)(\forall A 
% \in I^{\gt}_t)(A \subseteq L^1_a \rightarrow 
% (\exists b \in M_2)(A \subseteq L_{(a,b)}))$.
% 2020-12-24 17:20 2222
% Hence by the induction hypothesis for clause (B)(g) we have

Fix $ a \in M_1 $ and notice that  % 2021-01-18 14:20 
$\Lim_{\gt}(\bar{\bbQ} \restriction L^1_a),\bigcup\limits_{b \in L_2}
\Lim_{\gt}(\bar{\bbQ} \restriction L_{(a,b)})$ are equal as quasi
orders.  
% 2020-12-24 17:19 3333 Well 
Next we have to verify that for every 
$ t \in L^1_a$  and $ A \in I^\mathfrak{t} _t$ 
for some $ b \in L_2 $
 we have $ t \in L_{a,b}$ and $ A \subseteq L_{a,b}$.
 By the assumption on $ \langle L^2_b : b \in M_2\rangle $ 
for some $ b \in M_2 $ we have 
$ t \in L_b \wedge A \subseteq L^2_b$,   
hence $ t \in L^1_a \cap L^2_b $ and $ A \subseteq 
  L^1_a \cap L^2_b = L_{a,b} $ 
so  % 2020-12-24 17:18 so clearly 
  this $ b $
is as required.
% 2020-12-24 17:19 444  
Hence by the induction hypothesis for clause (B)(g) we have
$\Lim_{\gt}(\bar{\bbQ} \restriction L^1_a),\bigcup\limits_{b \in L_2}
\Lim_{\gt}(\bar{\bbQ} \restriction L_{(a,b)})$ are equal as quasi
orders
 
As this holds for every $a \in M_1$ and $M_1$ is directed
we get $\bigcup\limits_{a \in M_1} \Lim_{\gt}
(\bar{\bbQ} \restriction L^1_a),\bigcup\limits_{a \in M_1} \, 
\bigcup\limits_{b \in M_2} \Lim_{\gt}(\bar{\bbQ} \rest 
L_{(a,b)})$ are equal as quasi orders.  But 
the second is equal to $\bigcup\limits_{(a,b) \in M} 
\Lim_{\gt}(\bar{\bbQ} \restriction L_{(a,b)})$ so we are done.
\bigskip

% % 2020-12-03 16:30 \begin{enumerate} 
\underline{Clause (A)(c)} % 2020-12-03 16:33 \item[(c)] 

${ \name{\bar \nu }}$ is the 
 sequence of generics  of $ \bar{ \mathbb{Q} } $ 
 means 
 \begin{enumerate} 
 \item[($ \alpha $)] ${ \name{\bar \nu }}= \langle  \name{ \nu }_t:
    t \in L^\mathfrak{t} \rangle $,
 \item[($ \beta $)] for each $ t \in L^\mathfrak{t} $, $ \name{ \nu }_t$  is 
    a $ \Lim(\mathbb{Q} \upharpoonright {K^{\dagger}_t} )$-name of a function
   from a set of ordinals to $ \{ 0,1\}$ or to
    $ {\mathscr H} ({\aleph_0} ) $; for simplicity with domain $ u_t $  
    recalling $ K^{\dagger}_t = K_t \cup \{ t \}$    % 2020-12-04 18:07 $ where
 \item[($ \gamma $)]
 % 2020-12-04 18:06 $ K^{\dagger}_t = K_t \cup \{ t \} $\item[(*)] 
 if $ L \subseteq L^ \mathfrak{t} $ is $ \bar{ K } $-closed
 then $ \bar{ \nu } \upharpoonright L $ is a generic for
 $ \Lim (\bar{ \mathbb{Q} } \upharpoonright L  )  $
 \end{enumerate} 
 
\underline{Clause (A)(d)}  % 2020-12-03 16:35 \item[(d)] 

$ \bar{ u } = \bar{ u } [\bar{ \mathbb{Q} } ]
= % 2021-01-18 14:22 
   \langle u_t : t \in L^\mathfrak{t}  \rangle $  % 2020-12-24 17:21 
   is the parameters domain sequence of $ \bar{ \mathbb{Q} } $
   means that % 2021-01-18 14:23 
   each $ u_t = u(t) $ is a set of ordinals, for simplicity, 
   recalling that $\Vdash_{\Lim _ \mathfrak{t} (\bar{ \mathbb{Q} } \upharpoonright K_t )}
     {`}{`}\name{ \eta }_t $ is a function from $ u_t $ to 
     $ {\mathscr H} ( {\aleph_0} ) "$    % 2021-01-18 14:24 ", )
   
\underline{Clause (A)(e)}  % 2020-12-03 17:15 \item[(e)] 

The class $ \mathbf{Q} = \mathbf{Q} _\fsi$ of  fsi-templates  % 2020-12-03 16:38 FSIT-s 
is the class of objects 
  $ \mathbf{q} $ of the form 
   $ (\mathfrak{t} , \bar{ K } , \bar{ u } , \bar{ \mathbb{Q} }  ) $ 
   which are as above  with $ \dom (\bar{ \mathbb{Q} })
   =   % 2021-01-18 14:25 
   L^ \mathfrak{t} $  
   and $ \mathbb{P} _ \mathbf{q} = 
            \Lim _ \mathfrak{t} (\bar{ \mathbb{Q} })$. 
   We may write $ \mathbf{q} $ instead $ \mathfrak{t} ,
   (\mathfrak{q} , \bar{ K } )$  or $ \bar{ \mathbb{Q} }$.
   
   For $ \partial $ regular uncountable let $ \mathbf{Q} ^{\fsi}_ \partial $
   be the class of $ \mathbf{q} \in \mathbf{Q} _{\fsi}$  satisfying
   $ \partial ( \mathbf{q} ) \le \partial $, similarly in other cases.
   
   Let $ \mathbf{Q} _{\dom}$ be the class of $ \mathbf{q} \in \mathbf{Q} $ 
   such that $ K^\mathbf{q} _t = \emptyset , 
     \mathbb{Q} ^\mathbf{q} _{t}$  is dominating real= 
   Hechler forcing  and $ \name{ \eta }_t $ the generic.  % 2021-02-10 05:00 .
   
   On $ \mathbf{Q} _{\cln }$  see \ref{ad.15}(2), $ \boxplus ^1 _ \mathbf{q} $. 
   % 2020-12-24 17:26 2020-12-03 16:41 XOB DEBT BDOQ xxxQQQ
   
\underline{Clause (A)(f)}  % 2020-12-03 16:36 \item[(f)]  

Let % 2021-02-10 04:59 
$ \partial ( \mathbf{q} )$ for $ \mathbf{q} \in \mathbf{Q}$  
be the minimal infinite cardinal $ \partial $ which is strictly 
bigger then $ |A|, |u_t| $  for 
 $ t \in L^\mathfrak{t} , A \in I^\mathfrak{t} _t$. 

 We define $ \partial ^- (\mathbf{q} )$ similarly but requiring only   % 2020-12-01 14:47 \theta 
 "bigger or equal".  % 2020-12-03 16:37 
% 2020-12-03 16:30 \end{enumerate} 

% 2020-12-03 16:42 ]

\noindent
\underline{Part (B)}:
\bigskip

\noindent
\underline{\bf First Case}:  $\zeta=0$.

Trivial.
\bigskip

\noindent
\underline{\bf Second Case}: $\zeta$ successor.

Similar to usual iterations, so easy using the definition and the
induction hypothesis except clause (f) which we prove in details. 

\noindent
\underline{ Clause (f)}:

Let $p,q,L,L_0$ be as in the assumption of clause (f).  Let $r \in
\Lim_{\gt}(\bar{\bbQ} \restriction L_1)$ be above $q$ there
% 2020-12-04 18:20 
and we should prove the $p,r$ are compatible 
in $\Lim_{\gt}(\bar{\bbQ} \restriction L_2)$. 
Let $t_*$ be the maximal member of  $L_2,L^-_\ell := 
L_\ell \setminus \{t_*\}$ hence $L^-_2 = L_2 \backslash \{t_*\} \in 
% 2020-12-04 18:20 c\ell(
I_{t_*} $   % 2020-12-04 18:21 
and $\dep(L^-_2) < \zeta,L^- :=  L \setminus \{t_*\}$.  
If $(t_* \notin L_0 \vee t_* \notin L_1)$ or just $t_* \notin \Dom(p) \cap \Dom(r)$  
then by the induction hypothesis applied to $L^-_1,L^-_2,L^-,L^-_0,
p \upharpoonright L^-_0,q \upharpoonright L^-,r \upharpoonright L^-_2$ 
we can find a common upper bound $r^*$ of $p 
\upharpoonright L^-_0,r \upharpoonright L^-_1$ in $\Lim_{\gt}
(\bar{\bbQ} \upharpoonright L^-_2)$ and $r^* \cup p \upharpoonright \{t_*\}  
\cup r \upharpoonright \{t_*\}$ is a common upper bound of $p,r$ as required.

So assume that $t_* \in \Dom(p) \cap \Dom(r) \subseteq L_0 \cap L_1$ and let 
${\bbP}_0 := \Lim_{\gt}(\bar{\bbQ} \upharpoonright L^-)$
and ${\bbP}_{\ell +1} := \Lim_{\gt}(\bar{\bbQ} \upharpoonright
L^-_\ell)$ for $\ell = 0,1,2$.

% 2021-01-27 11:10 XAJOB -2- 2020-12-25 16:41 
Now we get $p',r'$ by applying 
the   definition in % 2021-02-10 05:01 proof of   % 2020-12-24 17:32 
clause (A)(a)(iv) for $t_*$ with 
$(p \upharpoonright L^-_0,p(t_*)),(r \upharpoonright L^-_1,r(t_*)),(q
\rest L^-,q(t_*))$ here standing for $(p_1,\name q_1),(p_2,\name
q_2),(p_0,\name q_0)$ there getting $(p'_\ell,q'_\ell)$ for $\ell < 3$
as there.

By the induction hypothesis in ${\bbP}_3$ for the conditions $p'_0,p'_1,p'_2$
we can find a common upper bound $p^*$,
so % 2020-12-24 17:32 
by $(A)(a)(iv)$ conclusion we
are done.
\bigskip

\noindent
\underline{\bf Third Case}:  $\zeta$ limit.

So let $\langle L^2_a:a \in M \rangle$ be a $({\gt},
\bar K)$-representation of $L_2$ with $a \in M \Rightarrow 
\Dp_{\gt}(L^2_a,\bar K) < \zeta$ and let   
   % 2020-12-24 17:56 ^2
  $L^1_a = L_1 \cap L^2_a$.
\bigskip

\noindent
\underline{Clause (B)(a)}: 

Trivial.
\bigskip

\noindent
\underline{Clause (B)(b)}: 

Clearly $\Dp_{\gt}(L_1,\bar K) \le \zeta$ by Claim 
\ref{ad.3}(2) $(\alpha)$ hence $\Lim_{\gt}(\bar{\bbQ} \restriction
L_1)$ is  well defined by (A)(b) which we have already 
proved % 2020-11-30 15:53 
above, that is 
$\Lim_{\gt}(\bar{\bbQ}) = \Lim_{\gt}(\bar{\bbQ} \restriction L_2) = 
 \bigcup\limits_{a \in M_2}  
 \Lim_{\gt}(\bar {\mathbb{Q}} % 2020-11-30 15:54 Q 
 \restriction L^2_a)$ as 
quasi orders. 

Clearly $\langle L^1_a = L_1 \cap L^2_a:a \in M  
\rangle$ is a $({\gt},\bar K)$-representation of $L_1$ hence 
by the induction hypothesis (if $\Dp_{\gt}(L_1,\bar K) < \zeta$) or
by the uniqueness proved in (A)(b) (if $\Dp_{\gt}(L_1,\bar K) = 
\zeta$) we know that $\Lim_{\gt}(\bar{\bbQ} \restriction L_1) = 
\bigcup\limits_{a \in M} \Lim_{\gt}(\bar{\bbQ} \rest L^1_a)$ as 
quasi orders and by the induction hypothesis for (B)(b) we know 
$\Lim_{\gt}(\bar{\bbQ} \restriction L^1_a) \subseteq \Lim_{\gt}(\bar{\bbQ}
\restriction L^2_a)$ as quasi orders (for $a \in M$), and 
we can easily finish.
\bigskip

\noindent
\underline{Clause (B)(c),(d)}: 
 
Use the proof of clause (B)(b) noting that $L^1_a \le_{\gt} L^2_a$
and so we can use the induction hypothesis,  % HOSEP  BDOQ
but we elaborate.
For clause (c), let $p \in  \Lim_ \mathfrak{t} 
  (\bar{ \mathbb{Q}}  \upharpoonright L_2 ) $ 
so there  is   an element $ a \in M$
such that % 2020-12-05 07:54 
 $ p \in \Lim_ \mathfrak{t} (\bar{ \mathbb{Q}} \upharpoonright 
 L^   2   % 2021-01-18 14:28 1
 _a)$.
 Now $ p \upharpoonright L_1 = p \upharpoonright L^1_a $ and
 as $ L^1_a \le _ \mathfrak{t}  L^2_a$  by the induction hypothesis
 we have $ p \upharpoonright L^1_a  
 \le _  {\Lim _\mathfrak{t} (\bar{ \mathbb{Q}} \upharpoonright L_1)} p $   % 2020-12-24 17:59 
    as promised.
    
    For clause (d) we assume in addition that 
    $ p \upharpoonright L_1 \le q$  in
      $ \Lim _ \mathfrak{t} (\bar{ \mathbb{Q}} \rest  L_1)$. 
    Let $ r_1 = 
      q \cup (p \upharpoonright (L_2 \setminus L_1)$.
    Easily $ r_1 $ is an upper bound of $ p, q $  in 
    $ \Lim _\mathfrak{t} (\bar{ \mathbb{Q} }\upharpoonright L_2 ) $.
    Assume further that $ r_2$ is another common upper bound of $ p, q$.
    As $M$ is directed we can choose $ a \in M $  such that 
    $ p, q, r_1, r_2 \in  \Lim _\mathfrak{t} (\bar{ \mathbb{Q}} \upharpoonright L^2_a  ) $
    but $ q , p \upharpoonright L^1\in 
      \Lim _\mathfrak{t} (\bar{ \mathbb{Q}} \upharpoonright L^1 _a ) $. 
      % 2021-01-18 14:29 By the induction  
    Hence  by the induction hypothesis  $ r_1 \le  r_2$  in 
      $ \Lim _\mathfrak{t} (\bar{ \mathbb{Q}} \upharpoonright L^2_a )  $ 
      so we can finish .
% (i.e. if $p \in 
% \Lim_{\gt}(\bar{\bbQ} \restriction L_2)$, as $M$ is directed there is
% $a \in M$ such that $\Dom(p) \subseteq L^2_a$, now $a \le_M b
% \Rightarrow p \restriction L^1_b = p \restriction L^1_a$ and we can\bigskip

\noindent
\underline{Clause (B)(e)}: 

The statements (i) + (ii) hold by clause (b). 

The statement (iii) holds: let ${\cI}$ be a predense subset of
$\Lim_{\gt}(\bar{\bbQ} \restriction L_1)$, let $p \in
\Lim_{\gt}(\bar{\bbQ})$, so for some $a \in M$ we have $p \in 
\Lim_{\gt}(\bar{\bbQ} \restriction L^2_a)$.  

By the induction hypothesis applying clause (B)(e) to 
$L^1_a,L^2_a$ we have $\Lim_{\gt}(\bbQ \restriction L^1_a) 
\lessdot \Lim_{\gt}(\bbQ \restriction L^2_a)$, hence as $p \in 
\Lim_{\gt}(\bar{\bbQ} \restriction L^2  % 2020-12-04 18:22 1
_a)$ clearly there is $q \in
\Lim_{\gt}(\bar{\bbQ} \restriction L^1_a)$ such that
$p$ is compatible with $r$ in $\Lim_{\gt}(\bbQ \restriction L^2_a)$ 
whenever $\Lim_{\gt}(\bbQ \restriction L^1_a) \models ``q \le r"$.  Now by the
assumption on ``${\cI} \subseteq \Lim_{\gt}
(\bar{\bbQ}\restriction L_1)$ is predense", as $q \in
\Lim_{\gt}(\bar{\bbQ} \restriction L_1)$ 
(by clause (B)(b)) we can find $q_0 \in {\cI}$
and $q_1$ such that $\Lim_{\gt}(\bar{\bbQ} \restriction L_1) \models
q_0 \le q_1 \wedge  q \le q_1$, so for some $b \in M$ we have
$q,q_0,q_1 \in L^1_b$ and $a \le_M b$ (as $M$ is directed).  Now we
consider $p,q,L^1_a,L^2_a,L^1_b,L^2_b$ and apply clause (B)(f).
\bigskip

\noindent
\underline{Clause (B)(f)}: 

Easy to check using clause (f) for the $L^2_a$'s, which holds by the
induction hypothesis.
\bigskip

\noindent
\underline{Clause (B)(g)}: 

Let $M_2 =: M$ (and recall $M_1$ that is from clause $(B)(g))$.
For each $a_1 \in M_1$, clearly $\Dp_{\gt}(L_{a_1},  % 2021-01-18 14:50 a_1
\bar K) \le \zeta$ as
$L_{a_1}  % 2021-01-18 14:50 
\subseteq L_2$ and $\langle L_{a_1} \cap L^2_{a_2}:a_2 \in M_2
\rangle$ is a $({\gt},\bar K)$-representation of $L_{a_1}$   % 2021-01-18 14:51 
and
$\Dp_{\gt}(L_{a_1} \cap L^2_{a_2},\bar K) \le \Dp_{\gt}(L^2_{a_2},\bar
K) < \zeta$ hence by (A)(b) we know $\Lim_{\gt}(\bar{\bbQ}
\restriction L_{a_1}) = \bigcup\limits_{a_2 \in M_2} \Lim_{\gt}(\bar{\bbQ} 
\restriction (L_{a_1} \cap L^2_{a_2}))$.
The rest should be clear.
\bigskip

\noindent
\underline{Clause (B)(h)}: 

Easy.
If $ t \in L^*$  then $ L_1 \setminus L^* \in I^ \mathfrak{t} _{t}$
hence for some $ a \in M $ we have $ L_1 \setminus L ^* \subseteq L_a$,
   % 2021-01-18 15:24 and \wilog \, $ t $
and the rest should be clear; and if $ L^* $ 
is empty % 2021-01-18 14:48 
this is easier. 
% (Note that every if we choose to assume in the definition 
%  that  $I^\mathfrak{t} _t = \{ A: $ for some  $\bar{ K } $ -closed\in L_a $ 
 % 2021-01-18 15:25 $ B \subseteq L_1$ we have $ A \subseteq B \cup K_t \}  $  
 % 2021-01-18 15:25 this would have worked.
% If $t^* \in L_1 \backslash L^*$, then $[L^*]^{\le \aleph_0}
% \subseteq I^{\gt}_{t^*}$, hence ecountablevery  $A \subseteq L^*$
% is included in some $L_a$ so this holds for $L_1$ too hence every
% $\Lim_{\gt}(\bar{\bbQ} \restriction L^*_1)$-name of a member of
% $\name{\bbQ}_{t,\bar \varphi_t,\name \eta_t}$ is a 
% $\Lim_{\gt}(\bar{\bbQ} \restriction (L_1 \cap L_1))$-name for 
% some $a \in M$.  
% The rest is easier.
\bigskip

\noindent
\underline{Clause (B)(i)}:

Easy.
\bigskip

\noindent

\underline{Clause (B)(j)}:

Sub-clause $ (\beta  )$ is clear by  % 2020-12-24 18:09 
the definition of  $ \Lim_ \mathfrak{t} (\bar{ \mathbb{Q} })$ so we shall deal with sub-clause $ (\alpha  )$.

So let $p_\alpha \in \Lim_{\gt}(\bar{\bbQ})$ for $\alpha <
\omega_1$; let $w_\alpha = \Dom(p_\alpha)$ and without loss of generality
$\langle w_\alpha:\alpha < \omega_1 \rangle$ is a $\Delta$-system with heart
$w$.  

A natural way fails because if $\langle L^*_\alpha:\alpha <
\omega_1\rangle$ is increasing continuous, $L^*_\alpha \subseteq L_2$
is  $ \bar{ K } $ -closed,  % 2020-11-30 16:20 (countable), 
$\subseteq$-increasing continuous then
$\langle \Lim_{\gt}(\bar{\bbQ} \rest L^*_\alpha)  % 2020-11-30 16:20 
:\alpha < \omega_1\rangle$ is
% 2020-11-30 16:21 not 
$ \lessdot $-increasing but not necessarily  % 2020-11-30 16:34 
continuous.  % 2020-11-30 16:18 contradiction.

Let $\{t_0,\dotsc,t_{n-1}\}$ list 
$ w $ % 2020-11-30 16:34 
without repetitions such that\footnote{As $ L^ \mathfrak{t} $
is a linear order, this mean $ t_0 < t_1 < \dots $.
}
$\ell
< k < n \Rightarrow \neg(t_k \le  t_\ell)$ 
and let $L^*_\ell$ be defined by
induction on $\ell \le 2n+1 = 2(n-1) +3$ as follows:
\mn
\begin{enumerate}
\item[$\bullet$]  $L^*_0 = \emptyset$
\sn
\item[$\bullet$]  $L^*_1 = \{s:s < t_k$ for every $k < n\}$
\sn
\item[$\bullet$]  $L^*_{2 \ell +2} = L^*_{2 \ell + 1} \cup \{t_\ell\}$ 
\when \, % 2021-02-10 05:06 if
  $ {\ell} < n $  equivalently % 2020-12-05 08:17 
  $2 \ell +2 < 2n  + 1$
\sn
\item[$\bullet$]  $L^*_{2 \ell +3 } = L^*  % 2021-02-10 05:05  
_{2 \ell +2} \cup \{s:s < t_k$
  for every $k \in \{\ell +1,\dotsc,n-1\}\}$
   \when \, % 2021-02-10 05:06 if 
   $ {\ell} < n $ equivalently % 2020-12-05 08:18 
   $ 2 {\ell} + 3 \le % 2020-12-05 08:18 <
   2n + 1 $  % 2020-11-30 16:38 
\end{enumerate}
\mn
So $L^*_{2n+1} = L_2$.

Clearly
\mn
\begin{enumerate}
\item[$\oplus$]  $\langle L^*_\ell:\ell \le 2n+1\rangle$ is
  % 2020-12-24 18:14 $\le_{\gt}$
  a $ \subseteq $-increasing  
  sequence of initial segments  of $ L_2$ 
  hence is $ \le _ \mathfrak{t} $-increasing 
\end{enumerate}
\mn
We prove by induction on $\ell \le 2n+1$ that
\mn
\begin{enumerate}
\item[$(*)_\ell$]  for some $q  % 2020-11-30 16:39 
_\ell \in \Lim_{\gt}(\bar{\bbQ} \rest
  L^*_{\ell} )  $ 
  we have $ q _ {\ell}  \Vdash_{\Lim_{\gt}(\bar{\bbQ} \rest L^*_\ell)} ``p_\alpha
  \upharpoonright L^*_ {\ell} % 2020-12-05 08:20 \alpha  % 2020-11-30 16:40 
  \in \name{\mathbf G}$ for $\aleph_1$ ordinals $\alpha"$.
\end{enumerate}
\bigskip

\noindent
\underline{Case 1}:  $\ell=0$ trivial, 
(e.g. the empty $ q \in \Lim(\bar{ \mathbb{Q} }\upharpoonright L^*_0)$).
  % 2021-01-18 15:59 
\bigskip

\noindent
\underline{Case 2}:  $\ell=1$

% % 2021-02-10 05:07 [XAJOB - jimvoni mitnaged xob\; lamah 
% $ t_0 \in L^*_1 $ we-laken 5im 
% $\alpha < \beta < \omega _1  $ 5az
% $ t_0 \in w_ \alpha \cap w_\beta \cap L^*_1$]

As $\langle w_\alpha \cap L^*_1:\alpha < \omega_1\rangle$ are pairwise
disjoint, 
every $q \in \Lim_{\gt}(\bar{\bbQ} \rest L^* _1)$ is compatible
with $p_\alpha \rest L^*_1$ for all but finitely many $\alpha$'s, so
this follows.
\bigskip

\noindent
\underline{Case 3}:  $\ell = 2i+3$

Recall $q  % 2020-11-30 16:42 
_{2i +2} \in \Lim_{\gt}(\bar{\bbQ} \rest L^*_{2 i  % 2021-01-18 16:00 \ell
+2  % 2020-11-30 16:43 
})$ has
been chosen.  Now assume % 2020-11-30 16:43 if 
$q_{2i+2} \le q \in \Lim_{\gt}(\bar{\bbQ}
\rest L^*_{2i+2})$ and $\alpha < \omega_1$.  % 2020-12-05 08:21 
  \Then\, $w^* = \{\gamma <
\omega_1   % 2020-11-30 16:45 
:w_\gamma \cap \dom(q) \nsubseteq w\}$ is finite, hence
recalling $\Lim_{\gt}(\bar{\bbQ} \rest L^*_{2i+2}) \lessdot
\Lim_{\gt}(\bar{\bbQ} \rest L^*_{2i +3})$ by % 2020-12-05 08:21 
the induction hypothesis,
there is $\beta \in (\alpha,\omega_1) \backslash w^*$ such that
$p_\beta  % 2020-11-30 16:46    \alpha  
\rest L^*_{2i+2},q \rest L^*_{2i+2}$ are compatible hence
there is $q_1 \in \Lim_{\gt}(\bar{\bbQ} 
  \rest L^*_{2i+2  % 2020-11-30 16:47 i
})$ above both.

It suffices to prove that $q_1,q,p_\beta$ has a common upper bound 
(as $ \alpha $ was an arbitrary countable ordinal).  % 2020-12-05 08:22 

We define a function $r$ by:
\mn
\begin{enumerate}
\item[$\bullet$]  $\dom(r) = \dom(q_1) \cup \dom(q) \cup \dom(p_\beta)$
\sn
\item[$\bullet$]  if $s \in \dom(q_1)$ then $r(s  % 2020-11-30 16:47 r 
) = q_1(s)$
\sn
\item[$\bullet$]  if $s \in \dom(q) \backslash L^*_{2i+2}$,
  equivalently $s \in \dom(q) \backslash \dom(q_1)$ then $r(s) = q(s)$
\sn
\item[$\bullet$]  if $s \in (\dom(p_\beta) \backslash L^*_{2i+2})$,
equivalently $s \in \dom(p_\beta  % 2021-01-18 16:02 \alpha
) 
\backslash \dom(q_1)$ then $r(s)
 = p_\beta % 2021-01-18 16:02 \alpha
 (s)$.
\end{enumerate}
\mn
It is easy to verify the ``equivalently" % 2020-11-30 16:48 ce" 
and as $\dom(q_1) \cap
\dom(p_\beta  % 2021-01-18 16:03 \alpha
) \subseteq L^*_{2i+2}$, the function $r$ is a well
defined function.  Also $r \in \Lim_{\gt}(\bbQ \rest L^*_{2i+3})$ as
its domain belongs to $[L^*_{2i+3}]^{< \aleph_0}$ and each $r(s)$ is
as required.

Why is $r$ above $q_1$?  Because $r \rest L^*_{2i+2} = q_1$.

Why is $r$ above $p_\beta $?  % 2021-01-18 16:03 \alpha$? 
By \ref{ad.5}(B)(d) recalling $\oplus$
above.

Why is $r$ above $q$?  By \ref{ad.5}(B)(d).

So we are done proving this case.
\bigskip

\noindent
\underline{Case 4}: $\ell = 2i+2$

We can find $\mathbf G \subseteq \Lim_{\gt}(\bar{\bbQ} \rest
L^*_{2i+1})$ generic over $\mathbf V$ such that 
$W = \{\alpha:p_\alpha \rest L^*_{2i+1} \in \mathbf G\}$ 
is uncountable. 
 Let $ t=t_ i $,  % , % 2021-02-10 05:09 {\ell} $,  % 2021-01-18 2021-01-18 16:05 
for each $\alpha \in  % 2020-11-30 16:52  % 2021-02-10 07:37 
W$ there is a $\bar K$-closed
$A_\alpha \in I^\mathfrak{t}_t $  % 2021-01-18 16:07 {t\_i }$ 
   % 2020-11-30 16:53 \ell}$ 
such that 
$p_\alpha(t)$ is a $\Lim_{\gt}(\bar{\bbQ} \rest A_\alpha)$-name.  So
as $\Lim_{\gt}(\bar{\bbQ} \rest A_\alpha) \lessdot
\Lim_{\gt}(\bar{\bbQ} \rest L^*_{2i+1})$ clearly  % 2020-12-05 08:23  
in $\mathbf V[\mathbf
G],q'_\alpha = p_\alpha % 2020-11-30 16:55 p
(t_\ell)[\mathbf G]$ is well defined and by absoluteness
(i.e. (A)(a)) is a member 
$\bbQ^{\mathbf V[\mathbf G]}_{t, \name{ \eta }_t}$.

Also $\mathbf V[\mathbf G] \models ``
\name{\bbQ}_{t,\name{ \eta }_t}^{V[\mathbf G]}$ satisfies
the c.c.c."  hence for some $\alpha_1 \ne \alpha_2$ from
$W, q_{\alpha _1}, q_{\alpha _2 }  $ % 2020-11-30 17:02 q_\alpha,q_\beta$ 
are compatible in $\bbQ_{t, \name{ \eta }_t}^{\mathbf V[\mathbf G]}$,
but 
% 2020-12-05 08:24 BDOQ  2020-11-30 17:00 
$\bbQ_{t, \name{ \eta }_t}   % 2020-11-30 17:01 t
^{\mathbf V[\mathbf G]}$ is ``too big".

Let $A = A_{\alpha_1} \cup A_{\alpha_2}$ so $A$ is a 
   % 2020-11-30 17:03 countable
$ \bar{ K}$-closed subset of   % 2020-11-30 17:39 $\Lim_{\gt}(
$ L^*_{2i+1}$ and it belongs to
$I_t$.  So
% 2020-11-30 17:45 $\Lim_{\gt}(\bar \mathbb{Q} % 2020-11-30 17:04 Q
$\Lim_{\gt}(\bar{\mathbb{Q}} \rest A_{\alpha_\iota }) \lessdot  
\Lim_{\gt}(\bar{\bbQ} \rest A)$ 
for $ \iota = 1,2$
hence by absoluteness
$q_{\alpha_1},q_{\alpha_2}$ belong to $\Lim_{\gt}(\bar{\bbQ} \rest A)$
and as $\Lim_{\gt}(\bar{\bbQ} \rest A) \lessdot \Lim_{\gt}(\bar{\bbQ}
\rest L^*_{2i+1})$ they are compatible, so we can finish easily.

So we have carried the induction hence the proof of (B)(j).
\bigskip

\noindent
\underline{Clause (k)}:  Easy.

\begin{claim}
\label{ad.6}
1) Assume
\mn
\begin{enumerate}
\item[$(a)$]  ${\gt}$ is an $\FSI$-template, $\Dp_{\gt}(L,\bar K) 
< \infty$ i.e. $\bar K$ is a smooth ${\gt}$-memory choice
\sn
\item[$(b)$]  $\bar{\bbQ} = \langle \name{\bbQ}_{t,\name \eta_t}:t 
\in L \rangle$ is a $({\gt},\bar K$)-iteration of def-c.c.c. forcing notions, 
   so $ L \subseteq L_ \mathfrak{t} $ is $ \bar{ K } $-closed
\sn
\item[$(c)_1$]  $L_1,L_2 \subseteq L$ and $L_1 < L_2$ (that is $(\forall
t_1 \in L_1)(\forall t_2 \in L_2)(L^{\gt} \models t_1 < t_2))$ and
$t \in L_2 \Rightarrow L_1 \in I^{\gt}_t$ 
and $ L_1, L_2 $ are $ \bar{ K } $-closed
% 2020-11-30 17:04 or at least $t \in L_2 \wedge  L' \subseteq L_1 \wedge  |L'| \le \aleph_0 \Rightarrow L' \in I^{\gt}_t$
and $L = L_1 \cup L_2$.
\end{enumerate}
\mn
\Then
\mn
\begin{enumerate}
\item[$(\alpha)$]  $\Lim_{\gt}(\bar{\bbQ})$ is actually a definition of
a forcing (in fact a c.c.c. one) so meaningful in bigger universes, moreover
for extensions 
(by c.c.c. forcings)  % 2020-11-30 17:06 
$\mathbf V_1 \subseteq \mathbf V_2$ of
$\mathbf V = \mathbf V_0$ (with the same ordinals of course), 
we \footnote{of course possibly $L_1 = \emptyset$} 
get $[\Lim_{\gt}(\bar{\bbQ})]^{\mathbf V_1} 
\subseteq_{\ic} % 2020-11-30 17:07 
[\Lim_{\gt}(\bar{\bbQ})]^{{\mathbf V}_2}$
  (see \ref{z2}(3))  % 2020-12-05 08:26 
and every maximal antichain ${\cI}$ of $\Lim_{\gt}(\bar{\bbQ})$ 
from $\mathbf V_1$ is a maximal antichain of 
$\Lim_{\gt}(\bar{\bbQ})$ (in $\mathbf V_2$).  % 2021-02-10 07:38 

% 2021-02-10 05:10 [XAJOB jimvoni : mikol maqom $ L$  sgurah legabey $ \bar{ K } $]
Recall  that  $ L $ is $ \bar{ K } $-closed.
\sn
\item[$(\beta)$]  $\Lim_{\gt}(\bar{\bbQ})$ is in fact 
$\bbQ_1 * \name{\bbQ}_2$ where $\bbQ_1 = \Lim_{\gt}
(\bar{\bbQ} \restriction L_1)$ and $\bbQ_2 = 
[\Lim_{\gt}(\bar{\bbQ} \restriction L_2)]^{\mathbf V}
[\name G_{\bbQ_1}]$ (composition).
\end{enumerate}
\mn
2) Assume clauses (a), (b) of part (1) and
\mn
\begin{enumerate}
\item[$(c)_2$]  $L$ has a last element $t^*$ and let $L^- = L 
\backslash \{t^*\}$.
\end{enumerate}
\Then \, for any $\mathbf{G} ^- \subseteq \Lim_{\gt}(\bar{\bbQ} 
\restriction L^-)$ generic over $\mathbf V$, letting $\eta_{t^*} = 
\name \eta_{t^*}[\mathbf{G} ^-] \in \mathbf V[\mathbf{G} ^-]$ we have: 
the forcing notion $\Lim_{\gt}(\bar{\bbQ})/\mathbf{G} ^-$ 
is equivalent to 
$\cup \{\bbQ^{\mathbf V[\mathbf{G} ^-_A]}_{t^*,\eta_{t^*}}:
A \in I^{\gt}_{t^*}$ is $\bar K$-closed$\}$ where $\mathbf{G} 
^-_A =: \mathbf{G} ^- \cap 
\Lim_{\gt}(\bar{\bbQ} \restriction A)$ and $\eta_{t^*_1} =
\name \eta_{t^*}[\mathbf{G} ^-]$. 

\noindent
3) Assume clauses (a), (b) of part (1) and
\mn
\begin{enumerate}
\item[$(c)_3$]  $\langle L_i:i < \delta \rangle$ is an increasing
continuous sequence of initial segments of $L$ with union $L$ and 
$\delta$ is a limit ordinal.
\end{enumerate}
\mn
\Then \, $\Lim_{\gt}(\bar{\bbQ})$ is $\bigcup\limits_{i < \delta} 
\Lim_{\gt}(\bar{\bbQ} \restriction L_i)$, moreover $\langle \Lim_{\gt}
(\bar{\bbQ} \restriction L_i):i < \delta \rangle$ is 
$\lessdot$-increasing continuous.  

\noindent
4) If ${\gt}$ is not smooth then ${\gt} \restriction L$ is not
smooth for some 
countable $L \subseteq L^{\gt}$, moreover for every $ L'$ 
satisfying $ L \subseteq L' \subseteq L^\mathfrak{t} $.
% 2020-11-30 17:09 of cardinality $< \kappa$.

\noindent
5) Assume  $ \mathfrak{t} $ is smooth and 
$ \bar{ \mathbb{Q} } = \langle \mathbb{Q} _{t, \name{ \eta }_t}:
t \in L^ \mathfrak{t} \rangle $.
If $\bar{\bbQ}$ is not a $({\gt},\bar K)$-iteration of
def-c.c. forcing notions, \then \, $\bar{\bbQ} \rest
L$ is not a $(\gt \rest L,\bar K \rest L)$-iteration of 
c.c.c.-definition forcing notions for some $\bar K$-closed 
$L \subseteq L^{\gt}$ which 
is % 2021-01-18 16:10 
the union of $ \le  {\aleph_1} $ 
$ \bar{ K } $-countable sets.
% 2020-11-30 17:15 of cardinality $\le 2^{\aleph_0}$.
\end{claim}

\begin{PROOF}{\ref{ad.6}}
Straightforward (or read \cite{Sh:630}).  
\end{PROOF}

We now give sufficient conditions for: ``if we force by $\Lim_{\gt}
(\bar{\bbQ})$ from \ref{ad.5}, then some cardinal invariants are small or
equal/bigger or equal to % 2020-12-05 08:27 than 
some $\mu$".   % 2020-11-30 17:17 
The necessity of such a 
claim in our framework is obvious; we deal with two-place 
relations only as this is the case in the
popular cardinal invariants, in particular those we deal with.
\begin{claim}
\label{ad.7}
Assume ${\gt}$ is a smooth $\FSI$-template and $\bar K = 
\langle K_t:t \in L^{\gt} \rangle$ and $\bar{\bbQ} = 
\langle \name{\bbQ}_{t,\name \eta_t}:t \in L^{\gt} \rangle$ are as 
in \ref{ad.5} and $\bbP = \Lim_{\gt}(\bar{\bbQ})$. 

\noindent
1) Assume
\mn
\begin{enumerate}
\item[$(a)$]  $R$ is a Borel \footnote{
  here and below just enough
  absoluteness is enough, of course
}
two-place relation\footnote{
  Why not $ {}^{ \omega } 2$? Just as 
  this notation is more natural for $ \mathfrak{d}, \mathfrak{b}$, 
  our main concern here.
}
on ${}^\omega \omega$
(we shall use $<^*$  
for $ \mathfrak{b} $ and $ \mathfrak{d} $,  $ \subseteq ^*$ for 
$ \mathfrak{u} $   % 2020-11-30 17:21 
and for $ \mathfrak{s} $ % 2020-12-05 08:28 
we use $ \eta R_{\spl} \nu $ meaning  $\Rang( \nu ) \cap \Rang( \eta  ),
\Rang(\nu ) \setminus \Rang(\eta ) $ are both infinite; 
the intention is to use this 
   % 2021-01-18 17:47 % 2021-01-18 16:11 
for $ \mathfrak{s}$)
\sn
\item[$(b)$]  $L^* \subseteq L^{\gt}$
\sn
\item[$(c)$]  for every 
$ \bar{ K } $-countable   % 2020-11-30 17:24 $\bar K$-closed 
$A \subseteq L^{\gt}$ for some $t \in L^*$ we have $A \in I^{\gt}_t$
\sn
\item[$(d)$]  for $t \in L^*$ and $\bar K$-closed $A \in I^{\gt}_t$ 
which includes $K_t$, in $\mathbf V^{\Lim_{\gt}(\bar{\bbQ} \restriction A)}$  
we have $\Vdash_{\name{\bbQ}_{t,\name \eta_t}} ``\name \nu_t \in
{}^\omega \omega$ is an $R$-cover of the old reals, that is
$\rho \in ({}^\omega \omega)^{\mathbf V[\Lim_{\gt}(\bar{\bbQ} 
\restriction A)]} \Rightarrow \rho R 
\name \nu_t$" where 
% 2021-01-18 16:12 333
$\name \nu_t$ is the generic real of 
$\name{\bbQ}_{t,\name \eta_t}$
  or just a 
  $ \Lim_ \mathfrak{t} (\bar{ \mathbb{Q} }  \rest K^\dagger  _t) $-name.
  We may use % 2020-12-25 10:10 
   % 222
$ \name{ \nu }_t $  a $ \Lim_ \mathfrak{t} (\bar{ \mathbb{Q} } 
   \rest A_t ) $ with $ A_ t \in I^ \mathfrak{t} _t $
% 111 
% $\name \nu_t$ is
% % 2020-12-05 08:34 in $ \mathbf{V} ^{\Lim _ \mathfrak{t} (\bar{ \mathbb{Q}  \upharpoonright K_t)}$
% a name in the forcing
% $\name{\bbQ}_{t,\name \eta_t}$,  i.e., all this in
% $(\bbQ_{t,\name \eta_t}[\name G])^{\mathbf V[\name G]},\name G$ the
% generic subset of $\Lim_{\gt}(\bar{\bbQ} \restriction A)$; 
% 

%   not  % 2020-12-05 08:30 
% depending on 
% $A$
% for $ A $ large enough (so can use $ \langle A_t: t \in L^* \rangle $. 
% (Usually 

% 333
% $\name \nu_t$ is the generic real of 
% $\name{\bbQ}_{t,\name \eta_t}$, 

% 444
% and hence
% $\name{\bbQ}_{t,\name \eta_t}$ is interpreted in the universe 
% $\mathbf V^{\Lim_{\gt}(\bar{\bbQ} \restriction A)}$, 

% 555
% so 
% $\name \eta_t$ is determined by the generic; normally we 
% assume this absolutely  BDOQ ).
\end{enumerate}
\mn
\Then \, $\Vdash_{\bbP} ``(\forall \rho \in {}^\omega \omega)
(\exists t \in L^*)(\rho R \name \nu_t)$, i.e. 
$\{\name \nu_t:t \in L^*\}$ is an $R$-cover, which, if
$R = <^*$ means ${\gd} \le |L^*|$". 

\noindent
1A) 
% If we weaken assumption (d) to ``for some 
% $\name \nu'_t$, a   % 2020-11-30 17:27 
% $\Lim_{\gt}(\bar{\bbQ} \restriction K^{\dagger}_t)$-name" where 
% $K^{\dagger}_t = K_t  \cup \{ t\} $   % 2020-11-30 17:26 
% % 2020-11-30 17:26 or we use $\bar K^{\dagger} = \langle K^{\dagger}_t:t \in L^* \rangle, K^{\dagger}_t \subseteq L^{\gt}_{\le t}$ 
% we get $\Vdash_{\bbP} 
% ``\{\name \nu'  % 2020-11-30 17:28 
% _t:t \in L^*\}$ is an $R$-cover".  
If we weaken 
assumption (d) to $\Vdash_{\bbP}$ ``for every $\rho \in {}^\omega
\omega$ for some $t \in L^{\gt}$ and $\nu \in 
\mathbf V({}^\omega \omega)^{\mathbf V[\Lim_{\gt}(\bar{\bbQ} 
\restriction K^{\dagger}_t)]}$ we have $\rho R \nu"$ then $\Vdash_{\bbP} 
``(\forall \rho \in {}^\omega \omega)(\exists t \in L^*)(\exists \nu
\in \mathbf V^{\Lim_{\gt}(\bar{\bbQ} \restriction K^{\dagger}_t)})[\rho R
\nu]"$.  This implies that in $\mathbf V^{\bbP}$, if $R = <^*$ then
${\gd} \le \sum\limits_{t \in L^{\gt}} 
|\Lim _ \mathfrak{t} (\bar{ \mathbb{Q}} \upharpoonright K^{\dagger}_t)|$;   % 2020-11-30 17:30 ||K^{\dagger}_t|$;
we could use
$K^{\dagger}$-s index by other sets. 

\noindent
2) Assume
\mn
\begin{enumerate}
\item[$(a)$]  $R$ is a Borel two-place relation on ${}^\omega \omega$
(we shall use $<^*$ or  $ \subseteq ^*$ as above)  % 2020-11-30 17:31 
\sn
\item[$(b)$]  $\mu$ is a cardinality
\sn
\item[$(c)$]  if\footnote{
We can weaken Clause (c) by saying: for every set $ X $ 
of $ < \mu $ names of reals 
there is $ t \in L^\mathfrak{t} $ such that 
for each such name from $ X $ \dots. 
}
$L^* \subseteq L^{\gt},|L^*| < \mu$ \then \, for some
$t \in L^{\gt}$ and $\bar K$-closed $L^{**} \supseteq L^*$ 
we have $L^{**} \in I^{\gt}_t$ and in 
$\mathbf V^{\Lim_{\gt}(\bar{\bbQ} \rest L^{**})},
\Vdash_{\bbQ_{t,\name \eta_t}}$ ``some $\nu \in {}^\omega \omega$ is 
an $R$-cover of the old reals"; (usually $\name \nu$ is
the generic real of $\bbQ_{t,\name \eta _y)}$,  % 2021-01-18 16:15 \nu_t}$,
this we assume absolutely).
\end{enumerate}
\mn
\Then \, $\Vdash_{\bbP} ``(\forall X \in [{}^\omega \omega]^{< \mu})
(\exists \nu \in {}^\omega \omega)(\bigwedge\limits_{\rho \in X} 
\rho R \nu)"$ (so for $R = <^*$ this means ${\gb} \ge \mu$). 

\noindent
3) Assume
\mn
\begin{enumerate}
\item[$(a)$]  $R$ is a Borel two-place relation 
\footnote{
    so $R$ is defined in $\mathbf V$; if $R$ is
    from $\mathbf V^{\Lim_{\gt}(\bar{\bbQ} \restriction K)}$, we need
    partial isomorphism (see below) of $({\gt},\bar{\bbQ})$ extending $\id_K$
}
on ${}^\omega \omega$ (we use $R = \{(\rho,\nu):\rho,\nu \in {}^\omega 2$ and 
$\rho^{-1}\{1\},\nu^{-1}\{1\}$ are infinite with finite intersection$\}  $, 
noting that ${}^{ \omega } 2 \subseteq {}^{\omega} \omega $)
\sn
\item[$(b)$]  $\sigma, \kappa,\theta$ are cardinals 
and $\kappa \le \theta \le \lambda$  % 2020-12-01 14:45 \theta 
and $ \sigma ^+ \ge  \partial ( \mathfrak{t} )$  % 2020-12-05 08:38 < 
 with $ \cf(\partial ( \mathfrak{t} ))> {\aleph_0} $
\sn
\item[$(c)$]  if 
 % 2020-11-30 17:50 
$ t_{i,\exn} \in
L^{\gt}$ 
% 2020-11-30 17:53 list a \bar{ K } $-countable subset  of $ L^\mathfrak{t} $
for $i < i(*),\exn <   \sigma  $ % 2020-11-30 17:54 $\sigma % 2020-11-30 17:51 \omega$ 
and
$\kappa \le i(*) < \theta$ and each $\{t_{i,\exn}:\exn < \sigma \}$ is
$\bar K$-closed, \then \, we can find $t_\exn \in L^{\gt}$ for
$\exn < \sigma $ such that $\{t_\exn:\exn < \sigma \} \subseteq L^{\gt}$ is
$\bar K$-closed and:
\sn
% 2020-11-30 15:14 \begin  % \xi 
\item[$(*)$]  for every $i < i(*)$ for some $j < \kappa,
j \ne i$ and the mapping
$t_{i,\exn} \mapsto t_{i,\exn},t_{j,\exn} \mapsto t_\exn$ is a partial % 2020-12-05 08:36 2  % \xi 
isomorphism of $({\gt},\bar{ K }, \bar{\bbQ})$ (see Definition \ref{ad.8} below).
% 2020-11-30 15:14 \end{enumerate}
\end{enumerate}
\mn
\Then \, in $\mathbf V^{\bbP}$ we have
\mn
\begin{enumerate}
\item[$\boxtimes^R_{\theta,\kappa}$]  if $\rho_i,\nu_i \in {}^\omega
  \omega$ for $i < i(*)$ and $\kappa \le i(*) < \theta$ 
and $i \ne j \Rightarrow \nu_i R \rho_j$, \then \, we can find 
$\rho \in {}^\omega \omega$ such that $i < i(*) \Rightarrow \nu_i R \rho$.
\end{enumerate}
\end{claim}

\begin{PROOF}{\ref{ad.7}}
Straightforward, but being requested we give details: 
 
\noindent
1) Let $\name \rho$ be a $\bbP$-name of a member of
$({}^\omega \omega)^{\mathbf V^{\bbP}}$, so as $\bbP$ 
satisfies the c.c.c. (see \ref{ad.5}(B)(j)$(\alpha )$), for each $n$ 
there is a maximal anti-chain
$\{p_{n,i}:i < i_n\}$ such that $p_{n,i}$ forces a value of
$\name \rho(n)$ and, of course, $i_n$ is countable.
Let $M = \{a:a$ is a 
$ \bar{ K } $-countable  % 2020-11-30 17:57  $\bar K$-closed 
subset of $L^{\gt}\}$
partially ordered by inclusion,    % 2020-12-24 18:27 
so obviously $M$ is closed under countable unions and
$\cup\{a:a \in M\} = L^{\gt}$; and let $L_a=a$ for $a \in M$ so by
\ref{ad.5}$(B)(i)(\beta)$ we have $p \in \Lim_{\gt}
(\bar{\bbQ}) \Leftrightarrow 
  p \in % 2020-12-05 08:39 
  \cup\{\Lim_{\gt}(\bar{\bbQ} \rest L_a):a \in M\}$
but $\bbP = \Lim_{\gt}(\bar{\bbQ})$, hence for $n < \omega,i < i_n$ for some
$a_{n,i} \in M$ we have $p_{n,i} \in \Lim_{\gt}(\bar{\bbQ}
\restriction L_{a_{n,i}})$.   % 2020-12-05 08:40 
But $M$ is $\aleph_1$-directed so for some $a \in
M$ we have $\{a_{n,i}:n < \omega,i < i_n\} \subseteq \{c:c \le_M a\}$.  
Also by \ref{ad.5}(B)(e) we know $\Lim_{\gt}(\bar{\bbQ} \rest L_a) 
\lessdot \Lim_{\gt}(\bar{\bbQ}) = \bbP$, 
so $\name \rho$ is a $\Lim_{\gt}(\bar{\bbQ} \restriction L_a)$-name.  
Now by assumption (c) of what we
are proving, as $L_a \subseteq L$ is 
$ \bar{ K } $-countable, we can find $t \in L^*
\subseteq L^{\gt}$ such that $L_a \in I^{\gt}_t$.  Also we
know that $K_t \in I^{\gt}_t$ (see Definition \ref{ad.1}(2)(c)
hence $A =: K_t \cup L_a$ belongs to $I^{\gt}_t$ and is 
$\bar K$-closed; and easily also $B = A \cup \{t\}$ is $\bar K$-closed.

Clearly $A \subseteq B \subseteq L^{\gt}$ are $\bar K$-closed so as
above $\Lim_{\gt}(\bar{\bbQ} \restriction A) \lessdot
\Lim_{\gt}(\bar{\bbQ} \restriction B) \lessdot \Lim_{\gt}(\bar{\bbQ}) 
= \bbP$ and $\name \rho$ is a $\Lim_{\gt}(\bar{\bbQ} \rest A)$-name 
(hence also a $\Lim_{\gt}(\bar{\bbQ} \restriction B)$-name)  % 2020-12-05 08:41 
of a member 
of ${}^\omega \omega$.

Now by assumption (d), in $\mathbf V^{\Lim_{\gt}(\bar{\bbQ}
\restriction A)}$ we have $\Vdash_{\bbQ_{t,\name \eta_t}} 
``\name \rho R \name \nu_t"$, hence by \ref{ad.6}(2) we
know that $\Lim_{\gt}(\bar{\bbQ} \restriction B) = 
\Lim_{\gt}(\bar{\bbQ} \restriction A) * 
\name{\bbQ}_{t,\name \eta_t}$, so together $\Vdash_{\Lim_{\gt}(\bar{ \mathbb{Q} } \rest B)} 
``\name \rho R \name \nu_t"$ hence 
by the 
previous sentence and
obvious absoluteness we have $\Vdash_{\bbP} 
``\name \rho R \name \nu_t"$.  So as $\name \rho$ was any $\bbP$-name 
of a member of $({}^\omega \omega)^{\mathbf V^{\bbP}}$ we are done. 

\noindent
1A)  Same proof. 

\noindent
2) So assume $p \Vdash_{\bbP} ``\name X \subseteq {}^\omega \omega$ 
has cardinality $< \mu$".  As we can increase $p$ without loss
of generality for some $\theta < \mu$ we have $p \Vdash_{\bbP} 
``|\name X| = \theta"$ so we can find a sequence $\langle \name
\rho_\alpha:\alpha < \theta \rangle$ of $\bbP$-names of members of
$({}^\omega \omega)^{\mathbf V^{\bbP}}$ such that $p \Vdash_{\bbP}
``\name X = \{\name \rho_\alpha:\alpha < \theta\}"$.  Let $\{p_{\alpha,n,i}:i <
i_{\alpha,n}\}$ be a maximal antichain of $\bbP$, with $p_{\alpha,n,i}$
forcing a value to $\name \rho_\alpha(n)$ and $i_{\alpha,n}$ countable.

Define $M = \{a \subseteq L^{\gt}:a$ is  % 2020-12-24 18:32 a
 $ \bar{ K } $-countable% 2020-11-30 17:59 $\bar K$-closed
$\}$, 
so for each $\alpha < \theta,n < \omega,i < i_{\alpha,n}$ for 
some $a_{\alpha,n,i} \in M$ we have $p_{\alpha,n,i}
\in \Lim_t(\bar{\bbQ} \restriction L_{a_{\alpha,n,i}})$.  
So\footnote{
In the weaker version for some $ t $ for every $ \alpha $
for some $ A \in I^ \mathfrak{t} _t $  \dots. 
}
for some $\bar K$-closed $L^{**} \subseteq L^{\gt}$ and $t \in 
L^{\gt}$ we have $L^{**} \in I^{\gt}_t$ and 
$L_{a_{\alpha,n,i}} \subseteq L^{**}$  % 2020-12-25 10:24 a_{\alpha, n,i}
for $\alpha < \theta,n < \omega,i < i_{\alpha,n}$.  We now continue
as in part (1). 

\noindent
3) So assume $i(*) \in [\kappa,\theta)$ and 
$\Vdash_{\bbP}
``\name \nu_i,\name \rho_i \in {}^\omega \omega$ and 
$i \ne j 
\Rightarrow \name \nu_i R \name \rho_j"$.  So as above we can find
$ \bar{  K } $-countable   % 2020-12-25 10:23 $\bar K$-closed
$K^*_i \subseteq L^{\gt}$ such that
$\name \nu_i,\name \rho_i$ are $\Lim_{\gt}(\bar{\bbQ} 
\restriction K^*_i)$-names; \wilog \, $K^*_i \ne \emptyset$ and 
$ K^*_i$  has cardinality $ < \partial ( \mathfrak{t} )$ hence $ \le \sigma $.  % 2020-12-05 08:44 
% 2020-12-05 08:43  $|K^*_i| = \aleph_0$; this is impossible only if $L^{\gt}$ is finite and then all is trivial. 
Let $\langle 
  t_{i, \exn}  % 2021-01-18 16:41 
  : \exn< \sigma  % 2020-11-30 18:00 t_{i,n}:n < \omega  % \xi 
\rangle$ be a list of the members of $K^*_i$ possibly with % 2020-12-05 08:45 no 
repetitions.  
Let $f_{i,j}$ be the mapping from $K^*_j$ to $K^*_i$ defined by 
$f_{i,j}(t_{j,\exn}) = t_{i,\exn}$ if well defined.  % 2020-12-25 10:25 

We define two-place relations $E_1,E_2$ on $i(*)$ and on $i(*) \times
i(*)$ respectively by:
\mn
\begin{enumerate}
\item[$(a)$]   $i E_1 j$ \Iff \, $f_{i,j}$ is a 
well defined   % 2020-12-25 10:25 
partial isomorphism of 
(${\gt},\bar{ K}, \bar{\bbQ})$ such that $\hat f_{i,j}$ 
(see claim (B) of \ref{ad.8} below) maps 
$(\name \rho_j,\name \nu_j)$ to $(\name \rho_i,\name \nu_i)$
\sn
\item[$(b)$]  $(i_1,i_2) E_2 (j_1,j_2)$ \Iff \, $i_1 E_1 j_1,i_2 E_
1   % 2021-01-18 16:20 2
  j_2$ and $f_{i_1,j_1} \cup f_{i_2,j_2}$ is a partial isomorphism of
$({\gt},\bar{\bbQ})$.
\end{enumerate}
\mn
Easily
\mn
\begin{enumerate}
\item[$\otimes(i)$]  $E_1,E_2$ are equivalence relations over their
domains
\sn
\item[$(ii)$]   $f_{j,i} = f^{-1}_{i,j}$ or both are  
not % 2021-01-18 16:43 
well defined.  % 2020-12-25 10:26 
\end{enumerate}
\mn
% 2020-12-24 18:34 As 
% $|i(*)/E_1| < \cf(\kappa)$ (by clause (c) of the
% assumption) and we can replace $i(*)$ by $i(*) + \kappa$, \wilog \, 
% $i < \kappa \Rightarrow 0E_1i$.   % 2020-12-05 08:46   
Now we apply
assumption (c), and get $\langle t_\exn:\exn < 
\sigma \rangle$ 
and let $ K^* = \{t_\exn:\exn < \sigma   
\}$. 
By $(*)$ of
clause (c) and clause (A)(b) of Definition \ref{ad.8} below
for any $i,j < i(*)$ clearly $K^*_i \cup K^*_j$ 
and $K^*_i \cup K^* $
    % 2020-12-25 10:31 < i(*)
  % 2020-12-25 10:30 \{t_\exn:\exn < \sigma  % 2020-12-01 12:58 \xi 
are $\bar K$-closed (see the 
definition below).  For any $i < i(*)$ let $j_i <
\kappa$ be as in $(*)$ of clause (c) which means: $j_i \ne i$ and the
following mapping $g_i$ is a partial isomorphism of 
$({\gt}, \bar{ K } ,\bar{\bbQ}):  % 2020-12-05 08:47 
\Dom(g_i) = \{t_{i,\exn},t_{j_i,\exn}:\exn < \sigma % 2020-12-01 12:57 \xi 
\},g_i(t_{i,\exn})
= t_{i,\exn},g_i(t_{j,\exn}) = t_\exn$.  % \xi 

Let $\name \nu,\name \rho$ be $\Lim_{\gt}(\bar{\bbQ} \rest K^*)$-names
such that for some, equivalently any $i,\hat g_i$ maps 
$\name \nu_{j_i},\name \rho_{j_i}$ to $\name \nu,\name\rho$ 
respectively (this is O.K. as for any
$i_1,i_2$ we have $j_{i_1} E_1 j_{i_2}$ 
because % 2020-12-25 10:33 $j_{i_1} E_1 j_{i_2}$ hence 
$g_{i_2} \circ f_{j_{i_2},j_{i_1}}  = g_{i_1} \restriction
K^*_{j_{i_1}}$).  Now for any $i < i(*)$, as $j_i \ne i$, we know
$\Vdash_{\Lim_{\gt}(\bar{\bbQ} \restriction (K^*_i \cup
K^*_{j_i}))} ``\name \nu_i R \name \rho_{j_i}"$, so applying $g_i$ 
we have $\Vdash_{\Lim_{\gt}
    (\bar{ \mathbb{Q} }(K^*_i \cup K^*))} ``\name \nu_i R 
\name \rho"$.  So we have proved $\boxtimes^R_{\theta,\kappa}$.
\end{PROOF}
\bigskip

In \ref{ad.8} below we note that isomorphisms 
(or embeddings) of ${\gt}$'s tend to induce 
isomorphisms (or embeddings) of $\Lim_{\gt}(\bar{\bbQ})$, and
deal (in \ref{ad.9},  \ref{ad.10}) with some natural operations.  In
\ref{ad.8} we could use two ${\gt}$'s, but this can trivially be
reduced to one.
\begin{dc}
\label{ad.8}
Assume that ${\gt},\bar K$ and $\bar{\bbQ} = \langle 
\name{\bbQ}_{t,\name \eta_t}:t \in L^{\gt} \rangle$ are as in \ref{ad.5}.  
By induction on $\zeta$ we define and 
prove \footnote{if $K_t = \emptyset$ and all
$\name{\bbQ}_{t,\eta}$ have the same definition of forcing notion, 
as in our main case, we can separate the definition and claim.}
\mn
\begin{enumerate}
\item[$(A)$]  [Def] $\quad$ we say $f$ is a partial isomorphism of
$({\gt},\bar{\bbQ})$ of Depth $\le \zeta$ if: 

\hskip25pt (omitting $\zeta$ means for some ordinal 
$\zeta$; writing ${\gt}$ instead of 

\hskip25pt $({\gt}, \bar{ K } ,\bar{\bbQ})$ means we assume 
$\bbQ_{t,\name \eta_t} = \bbQ$, i.e. constant, $K_t = \emptyset$ 

\hskip25pt for every $t \in L^{\gt}$ and may say 
``${\gt}$-partial isomorphism")
\sn
\begin{enumerate}
\item[$(a)$]  $f$ is a partial one-to-one function from $L^{\gt}$ 
to $L^{\gt}$
\sn
\item[$(b)$]  $\Dom(f),\Rang(f)$ are $({\gt},\bar K)$-closed 
sets of depth $\le \zeta$
\sn
\item[$(c)$]  for $t \in \Dom(f)$ and $A \subseteq \Dom(f)$ we have 
$A \in I^{\gt}_t \Leftrightarrow f''(A) \in I^{\gt}_{f(t)}$
\sn
\item[$(d)$]  for $t \in \Dom(f)$, we have: $f$ maps $K_t$
onto $K_{f(t)}$ and $f \restriction K_t$ maps $\name \eta_t$
to $\name \eta_{f(t)}$, more exactly the isomorphism $\hat f$ which 
$f$ induces from $\Lim_{\gt}(\bar{\bbQ} \restriction K_t)$ onto $\Lim_{\gt}
(\bar{\bbQ} \restriction K_{f(t)})$ does this.
\end{enumerate}
\sn
\item[$(B)$]  [Claim] % 2020-12-24 18:38 $(a) 
$ \quad f$ induces naturally an isomorphism
which we call $\hat f$ from $\Lim_{\gt}(\bar{\bbQ} \rest \Dom(f))$ 
onto $\Lim_{\gt}(\bar{\bbQ}
\rest \Rang(f))$.
\end{enumerate}
\end{dc}

\begin{proof}  
Straightforward, recalling  we are assuming that
$ \bar{ \varphi }_ t $ is definable from 
$ \name{ \eta }_t$.
\end{proof}

\begin{definition}
\label{ad.9}
1) We say ${\gt} = {\gt}^1 + {\gt}^2$ if
\mn
\begin{enumerate}
\item[$(a)$]   $L^{\gt} = L^{\gt^1}  % 2021-01-18 16:55 
    + L^{\gt^2}$ (as linear orders) 
\sn
\item[$(b)$]  for $t \in L^{\gt^1},I^{\gt^1}_t = I^{\gt}_t$
\sn
\item[$(c)$]  for $t \in L^{\gt^2},I^{\gt^2}_t = \{A \subseteq
  L^{\gt}:A \cap L^{\gt^2} \in I^{\gt^2}_t\}$.
\end{enumerate}
\mn
So ${\gt}^1 + {\gt}^2$ is well defined if ${\gt}^1,{\gt}^2$
are disjoint, i.e. $L^{\gt^1} \cap L^{\gt^2} = \emptyset$. 

\noindent
2) We say ${\gt}^1 \le_{\wk} {\gt}^2$ iff
\mn
\begin{enumerate}
\item[$(a)$]  $L^{\gt^1} \subseteq L^{\gt^2}$ (as linear orders) and 
$t \in L^{\gt^1} \Rightarrow I^{\gt^1}_t \subseteq I^{\gt^2}_t$
\sn
\item[$(b)$]   if 
% 2020-11-30 18:11 for every countable
   \footnote{we may restrict ourselves to $\FSI$-templates ${\gt}$ of
globally countable, i.e., such that 
$A \in I^{\gt}_t$ and $t \in L^{\gt} \Rightarrow |A| \le
\aleph_0$}   %     % 
   $ s \in L^{\mathfrak{t}^1} $  then $ I^ {\mathfrak{t}^1} _ s   =
\{  % 2020-12-05 09:14 
A \in I^{\mathfrak{t} ^2}_t : A \subseteq L^{\mathfrak{t} ^1 }  \} $
% or locally countable with no loss.  We use this restriction
% as in Definition \ref{ad.11}, if $A_i \subseteq A_j \in I^{\gt}_t$ for 
% $i<j<\kappa$ then in $\gt^* = \gt^\kappa/D,
% \bigcup\{j^*_{D,t}(A_i):i < \kappa\} \in I^{gt^*}_{j_{D,\gt}}(t)$ 
% even if $\bigcup\{A_i:i < \kappa\} \notin I^{\gt}_t$.
% $A \subseteq L^{\gt^1}$ and $t \in L^{\gt^1}$ we have $A \in 
% I^{\gt^1}_t \Leftrightarrow A \in I^{\gt^2}_t$.
\end{enumerate}
\mn
3) If $\langle {\gt}^\zeta:\zeta < \xi \rangle$ is 
$\le_{\wk}$-increasing, $\xi$ a limit ordinal, we define ${\gt}^\xi 
=: \bigcup\limits_{\zeta < \xi} {\gt}^\zeta$ by

\[
L^{\gt^\xi} = \bigcup\limits_{\zeta < \xi} L^{\gt^\zeta} \qquad
\text{ (as linear orders)}
\]

\[
I^{\gt^\xi}_t = \cup \{I^{\gt^\zeta}_t:\zeta < \xi \text{ and } 
t \in L^{\gt}_\zeta\}
\]

\mn
Clearly $\zeta < \xi \Rightarrow {\gt}^\zeta \le_{\wk} {\gt}^\xi$.  
Such ${\gt}^\xi$ is called the limit of $\langle {\gt}^\zeta:\zeta < 
\xi \rangle$; now a $\le_{\wk}$-increasing
sequence $\langle {\gt}^\zeta:\zeta < \xi \rangle$ is called
continuous if for every limit ordinal $\delta < \xi$ we have 
${\gt}^\delta = \bigcup\limits_{\zeta < \delta} {\gt}^\zeta$.

\noindent
4) If $\langle {\gt}^\zeta:\zeta < \xi \rangle$ are pairwise disjoint
(that is $\zeta \ne \varepsilon \Rightarrow L^{\gt^\zeta} \cap
L^{\gt^\varepsilon} = \emptyset$) we define $\sum\limits_{\zeta < \xi}
{\gt}^\zeta$ by induction on $\xi$ naturally: for $\xi = 1$ it is 
${\gt}^0$, for $\xi$ limit it is $\bigcup\limits_{\varepsilon 
    % 2021-01-18 17:02 \zeta 
< \xi} 
(\sum\limits_{\zeta < \varepsilon} {\gt}^\zeta)$ and for 
$\xi = \varepsilon + 1$ it is $(\sum\limits_{\zeta < \varepsilon}
{\gt}^\zeta) + {\gt}^\varepsilon$, so $\xi_1 \le \xi_2 \Rightarrow
\sum\limits_{\zeta < \xi_1} {\gt}^\zeta \le_{\wk} \sum\limits_{\zeta < \xi_2}
{\gt}^\zeta$ (even an initial segment). 

\noindent
5) We can replace in 0) - 4) above ${\gt}^\zeta$ by $({\gt}^\zeta,
\bar K^\zeta)$. 
\end{definition}

\begin{claim}
\label{ad.10}
Let ${\gt}$ be an $\FSI$-template. 

\noindent
1) If $L^{\gt} = \emptyset$  or just is well ordered  % 2020-12-25 10:35 
% 2020-11-30 18:15 or just $L^{\gt}$ is finite  and $\{s:L^{\mathfrak t} \models s <t\} \in I^{\gt}_t$ 
\then \, ${\gt}$ is
smooth. 

\noindent
2) If ${\gt}^1,{\gt}^2$ are disjoint $\FSI$-templates, 
\then \, ${\gt}^1 + {\gt}^2$ is an $\FSI$-template and $\ell \in 
\{1,2\} \Rightarrow {\gt}^\ell \le_{\wk} {\gt}^1 + {\gt}^2$. 

\noindent
3) If ${\gt}^1,{\gt}^2$ are disjoint smooth $\FSI$-templates \then \,
$ \mathfrak{t} = % 2020-11-30 18:15 
{\gt}^1 + {\gt}^2$ is a smooth $\FSI$-template; moreover,
$\Dp_{\gt}(L^{\gt}) \le \Dp_{{\gt}^1}(L^{{\gt}^1}) + \Dp_{{\gt}^2}
(L^{{\gt}^2})$ and $\Dp_{\gt}(L^{{\gt}^\ell}) = 
\Dp_{{\gt}^\ell}(L^{{\gt}^\ell})$. 

\noindent
4) If $\langle {\gt}^\zeta:\zeta < \xi \rangle$ is an 
$\le_{\wk}$-increasing (\ref{ad.9}(2)) sequence of 
$\FSI$-templates and $\xi$ is a limit ordinal, \then \, ${\gt}^\xi 
=: \bigcup\limits_{\zeta < \xi} {\gt}^\zeta$ is an $\FSI$-template 
and $\zeta < \xi \Rightarrow {\gt}^\zeta \le_{\wk} {\gt}^\xi$. 

\noindent
5) If $\langle {\gt}^\zeta:\zeta < \xi \rangle$ is an increasing
continuous (see Definition \ref{ad.9}(3)) sequence of smooth
$\FSI$-templates and $\xi$ is a limit ordinal, \then \, ${\gt}^\xi =:
\bigcup\limits_{\zeta < \xi} {\gt}^\zeta$ is a smooth $\FSI$-template 
and $\zeta < \xi \Rightarrow {\gt}^\zeta \le_{\wk} 
{\gt}^\xi$ and $\Dp_{{\gt}^\xi}(L^{{\gt}^\xi}) \le
% 2020-11-30 18:17 \sum\limits_
% 2020-11-30 18:18 \sup {\zeta < \xi}
\sup\{ \Dp_{{\gt}[\zeta]}(L^{{\gt}[\zeta]})+ 1 :  \zeta < \xi \} $. 

\noindent
6) If $\langle {\gt}^\zeta:\zeta < \xi \rangle$ is a sequence of pairwise
disjoint [smooth] $\FSI$-templates, \then \, $\sum\limits_{\zeta < \xi}
\gt^\zeta$ is a [smooth] $\FSI$-template and $\langle 
\sum\limits_{\zeta < \varepsilon} \gt^\zeta:
\varepsilon \le \xi % 2020-11-30 18:19 \zeta 
\rangle$ is increasing continuous.

\noindent
7) In parts (1)-(6) we  % 2020-11-30 18:20 
can expand ${\gt}^\zeta$ by $\bar K^\zeta$ . 

\noindent
8)  Assume $\mathbf J$ is a linear order, $\mathbf t_x$ is a smooth
FSI-template 
for every $x \in \mathbf J$ and $\langle L^{\gt_x}:x \in
\mathbf J \rangle$ are pairwise disjoint (for notational simplicity) and
we define ${\gt}$ by: $L^{\gt} = \sum\limits_{x \in \mathbf J}
L^{\gt_x}$ (so $L^{\gt} \models s < t$ iff $(\exists x,y)(s
\in L^{\gt_x} \wedge t \in L^{{\gt}_y} \wedge x <_  \mathbf{J} 
    % 2021-01-18 16:59 J
y) \vee
(\exists x \in \mathbf J)(L^{{\gt}_x} \models s < t))$
 and $I^{\gt}_t = \{A \subseteq L^{\gt}:(\forall s \in A)
(s <_{L^{\gt}} t)$ and letting 
$ x \in \mathbf{J} $ be such that % 2020-12-05 08:53 
$t \in {\gt}^x$ we have 
$A \cap L^{{\gt}_x} \in I^{{\gt}_x}_t$ and
$\{y:y <_{\mathbf J} x,A \cap L^{{\gt}_y} \ne \emptyset\}$ is
finite$\}$.  Then ${\gt}$ is a smooth 
FSI-template % 2020-11-30 18:21 
(we can expand by $\bar K$'s) (use in
\S3).       
\end{claim}

\begin{PROOF}{\ref{ad.10}}
Easy, e.g.   % 2020-12-25 10:39 in 
part (3) is proved by induction on
$\Dp_{\gt}(L^{\gt})$ and part (6) by induction on $\xi$ and in
part (7) let $M$ be $[\mathbf J]^{< \aleph_0}$ ordered by inclusion and
$L_{\{x(1),\dotsc,x(n)\}} = \cup \{L^{{\mathfrak t}_{x(\ell)}}:\ell
=1,\dotsc,n\}$ for any $x(1),\dotsc,x(n) \in \mathbf J$.
\end{PROOF}
\bigskip

\centerline {$* \qquad * \qquad *$}
\bigskip

\begin{discussion}
\label{ad.10a}
1)
To prove our desired result $\CON({\ga} > {\gd})$ we need to construct
an $\FSI$-template ${\gt}$ of the right form.  Now we do it using a measurable
cardinal.  The point is that if we are given $\left< \langle t_{i,n}:n <
\omega \rangle:i < i(*) \right>,L^{\gt},i(*) \ge \kappa$ and $D$ is a
normal ultrafilter on $\kappa$, then in ${\gt}^\kappa/D$ the
$\omega$-sequence $\left< \langle t_{i,n}:i < \kappa \rangle/D:
n < \omega \right>$ is as required in \ref{ad.7}(3)(c), considering 
${\gt}^\kappa/D$ an extension of ${\gt}$.  

\noindent
2) We  shall deal with $ \mathfrak{s} $ only in \ref{ad.15}(2).

\noindent 
3) 
Note that our main old conclusion 
(i.e. \ref{ad.15}(1)) has two proofs. The first is shorter
and depends on \S1  and \ref{ad.13}, \ref{ad.14}.  The second is longer but does 
not.

\end{discussion}

\begin{definition}
\label{ad.11}
1)
For a  $ \mathbf{q} \in \mathbf{Q} $  and $ \partial = \partial (\mathbf{q}) < \kappa $
% 2020-11-30 18:28 template ${\gt}$ and an $(2^{\aleph_0})^+$-complete 
and an $ \partial ^+ $-complete 
ultrafilter $D$ on $\kappa$ 
(hence $ (2^ \partial )^+$-complete),
we define ${\gt}^* =: 
{\gt}^\kappa/D,\mathbf j_{D,{\gt}}$ and $\mathbf j_{D,{\gt}}({\gt})$ as follows:
\mn
\begin{enumerate}
\item[$(a)$]   we define ${\gt}^*$ by:
\[
L^{{\gt}^*} = (L^{\gt})^\kappa/D \text{ as a linear order}
\]
and if $t^* = \langle t_i:i < \kappa \rangle/D$ where 
$t_i \in L^{\gt}$ then we let $I^{\gt^*}_{t^*} = \{A:
\text{we can find } A_i \in I^{\gt}_{t_i}$ for $i < \kappa$ 
such that $A \subseteq \prod\limits_{i < \kappa} A_i/D\}$
\sn
\item[$(b)$]  We then let $\mathbf j_{D,{\gt}}$ be the canonical embedding 
of ${\gt}$ into ${\gt}^\kappa/D$ that is $\mathbf j_{D,{\gt}}(t) = \langle t:i <
\kappa \rangle/D$ for every $t \in L^{\gt}$ and 
\sn
\item[$(c)$]  let ${\gt}' = \mathbf j_{D,{\gt}}({\gt})$ be defined by
$L^{{\gt}'} = L^{{\gt}^*} \restriction \{\mathbf j_{D,{\gt}}(s):s
\in L^{\gt}\},I^{{\gt}'}_{\mathbf j_{D,{\gt}}(s)} = 
\{\{\mathbf j_{D,{\gt}}(t):t \in A\}:A \in I^{\gt}_s\}$. 
\newline
% [We can deal with $\bar K$, if $D$ if $(\bigcup\limits_{t \in L^{\gt}} 
% |K_t|^+)$-complete which holds as here as each $K_t$ 
% is countable and can deal also with $\bar{\bbQ}$ if we have 
% $< \com(D)$ kinds of $\bar \varphi_t$ (letting
% $\name \eta_i$ vary) which too holds here).]
\end{enumerate} 

\noindent 
2) Similarly for $ \mathbf{q} \in \mathbf{Q} $ instead $ \mathfrak{t} $.
% 2020-12-05 08:54 \end{enumerate}
\end{definition}

\begin{remark} 
We may allow $ \partial (\mathbf{q} ) \ge \kappa $ but presently not worth the trouble.
\end{remark} 

\begin{claim}
\label{ad.12}
In Definition \ref{ad.11}: 

\noindent
1) ${\gt}^\kappa/D$ is also an $\FSI$- template and 
$\mathbf j_{D,{\gt}}({\gt}) \le_{\wk} {\gt}^\kappa/D$
and $\mathbf j_{D,\mathbf t}$ is an isomorphism from ${\gt}$ onto
$\mathbf j_{D,{\gt}}({\gt})$. 

\noindent
2) If ${\gt}$ is a smooth $\FSI$-template \then \, ${\gt}^\kappa/D$ is 
a smooth $\FSI$-template. 

\noindent
3) Moreover, $\Dp_{{\gt}^\kappa/D}(L^{{\gt}^\kappa/D})
\le (\Dp_{\gt}(L^{\gt}))^\kappa/D$. 

\noindent 
4) Similarly  we define  % 2020-12-25 10:42  for
$ \mathbf{q} ^\kappa /D$  for $ \mathbf{q} \in \mathbf{Q} $; 
so $ u_t $ is increased if $ u^ \mathbf{q} _t $.
 is of cardinality $  \ge \kappa $ and similarly $ K^ \mathbf{q} _t$
\end{claim}
\begin{PROOF}{\ref{ad.12}}
Straightforward.
\end{PROOF}

\noindent
Now \ref{ad.13}, \ref{ad.14} below are used only in the short proof of
\ref{ad.15} depending on \S1, so you may ignore them.
\begin{definition}
\label{ad.13}
% 2020-12-06 07:51    1)
Fix $\aleph_0 < \kappa < \mu = \cf(\mu) < \lambda = \cf(\lambda) 
= \lambda^\kappa$ and $D$ a $\kappa$-complete (or just 
$(2^{\aleph_0})^+$-complete) uniform
ultrafilter on $\kappa$.  We define by induction on $\zeta \le \lambda$,
a  % 2021-01-18 17:06 
smooth $\FSI$-template ${\gt}_{\gamma,\zeta}$ for $\gamma < \mu$ such that:
\mn
\begin{enumerate}
\item[$(a)$]  ${\gt}_{\gamma,\zeta}$ is a smooth  $\FSI$-template  % 2020-12-25 10:45 
\sn
\item[$(b)$]  if $\gamma_1 \ne \gamma_2$ then ${\gt}_{\gamma_1,\zeta},
{\gt}_{\gamma_2,\zeta}$ are disjoint, i.e. 
$L^{{\gt}_{\gamma_1,\zeta}} \cap L^{{\gt}_{\gamma_2,\zeta}} = \emptyset$
\sn
\item[$(c)$]  for $\xi < \zeta$ we have ${\gt}_{\gamma,\xi} 
\le_{\wk} {\gt}_{\gamma,\zeta}$
\sn
\item[$(d)$]  if $\zeta$ is limit then ${\gt}_{\gamma,\zeta} =
\bigcup\limits_{\xi < \zeta} {\gt}_{\gamma,\xi}$, see \ref{ad.9}(3), 
\ref{ad.10}(6).
\sn
\item[$(e)$]  if $\zeta = \xi +1$ and $\xi$ is even, \then \, there is an
isomorphism $\mathbf j_{\gamma,\zeta}$ from $\sum\limits_{\beta \le \gamma} \,
{\gt}_{\beta,\xi}$ onto ${\gt}_{\gamma,\zeta}$ which is the identity
over ${\gt}_{\gamma,\xi}$
\sn
\item[$(f)$]  if $\zeta = \xi +1$ and $\xi$ is odd, \then \, there is an
isomorphism $\mathbf j_{\gamma,\zeta}$ from $({\gt}_{\gamma,\xi})^\kappa/D$
onto ${\gt}_{\gamma,\zeta}$ which extends the inverse of
$\mathbf j_{D,{\gt}_{\gamma,\xi}}$.
\end{enumerate}

% \noindent 
% 2) Assume in addition tha $ {\aleph_0} < \partial < \kappa $  
% male5 ! BDOQ 
\end{definition}

\begin{observation}
\label{ad.14}
The definition is \ref{ad.13} is legitimate.
\end{observation}

\begin{PROOF}{\ref{ad.14}}
By the previous claims.
\end{PROOF}
  
\begin{conclusion}
\label{ad.15}
Assume: $\kappa$ is measurable,\footnote{
Instead $ \lambda = \lambda ^ \kappa $ it suffice to demand 
$ \lambda = \lambda ^{\aleph_0} = \lambda ^ \kappa /D$.
This holds for any strong limit cardinal $ > \kappa $ 
of cofinality $ \not= {\aleph_0}, \not= \kappa $.
} 
$\kappa < \mu
= \cf(\mu) < \lambda = \cf(\lambda) = \lambda^\kappa$.
% 2020-12-03 16:47 \Then \, 
1)
 For some c.c.c. forcing notion $\bbP$ of cardinality 
$\lambda$, in $\mathbf V^{\bbP}$ we have ${\ga} = \lambda,{\gb} 
= {\gd} = \mu$  
hence $ \mathfrak{s} \le \mu$.  % 2021-01-18 17:14 

\noindent 
2) 
If in addition $ \partial = \cf( \partial )   % 2020-12-03 16:48 
< \kappa $ then for some $ \mathbb{P} $ as above
 in addition we have 
   $ \mathfrak{s} \ge \partial $ (hence  % 2021-01-18 17:14 
 $ \partial  \le \mathfrak{s} \le \mu $)
     % 2020-12-06 07:53 $ \theta \le  \mathfrak{s}  \le \partial  $ % 2020-12-03 16:49 = \mathfrak{s} $ BDOQ   % \theta 

\end{conclusion}

\begin{PROOF}{\ref{ad.15}}
\smallskip

\noindent
1)  \underline{Short Proof}: (depending on \S1).  

Let ${\gt}_{\gamma,\zeta}$ (for $\gamma < \mu,\zeta \le \lambda)$ 
be as in \ref{ad.13}.  Let ${\gt} = \sum\limits_{\gamma < \mu} 
{\gt}_{\gamma,\lambda}$ and let $\bar K = \langle K_t:t \in 
L^{\gt} \rangle,K_t = \emptyset$ and let $\bar{\bbQ} = \langle
\name{\bbQ}_t:t \in L^{\gt} \rangle$ with $\name{\bbQ}_t$ being 
constantly the dominating real forcing (= Hechler forcing).  

Lastly, let $\bbP = \Lim_{\gt}(\bar{\bbQ})$.

\noindent
The rest is as in the end of \S1.  But if we like to use \ref{ad.5},
etc. we need
\mn
\begin{enumerate}
\item[$\boxplus$]  $\bbQ_{\dom}$ is as required in \ref{ad.5}(A)(a)(i)-(iv), 
   % 2020-12-25 10:52 hapnayah 
 i.e. def - c.c.c.
\end{enumerate}
\mn
We elaborate concerning why $\bbQ_{\dom}$ 
satisfying sub-clause (iv) of the full definition of % 2020-12-25 10:48 
\ref{ad.5}(A)(a).

Given $p_\ell$ assume
\mn
\begin{enumerate}
\item[$(a)$]  $\bbP_0 \lessdot \bbP_\ell \lessdot \bbP_3$ (for
  $\ell=1,2$) be c.c.c. forcing
\sn
% 2020-12-25 14:19 
\item[$(b)$]
% 2020-12-25 14:30 \item[(*)] 
$\name {\mathbb{Q}} _\ell$
the $\bbP_\ell$-name of 
$\bbQ_\dom$ 
with the generic $ \name{ \nu}_ {\ell}$, 
   (in a sense they are the same name)  % 2020-12-25 13:10 
\sn
\item[$(c)$]  $(p_\ell,\name q_\ell) \in \bbP_\ell * \name{\bbQ}_\ell$
  for $\ell = 0,1,2$
\sn
\item[$(d)$]  $(p_0,\name q_0) \Vdash ``(p_\ell,\name{q}_\ell) \in
  (\bbP_\ell * \name{\mathbb{Q} }_\ell )/(\bbP_0 * \bbQ_0)"$ for $\ell=1,2$
\snn  
\item[(e)]  $ p_3 \in \mathbb{P} _3 $ is a common upper
bound of $ p_1, p_2$.
\end{enumerate}
\mn
Of course
\mn
\begin{enumerate}
\item[$(*)_1$]  we can replace $(p_\ell,\name q_\ell)$ for $\ell < 3$
  by $(p'_\ell,q'_\ell)$ above $(p_\ell,q_\ell)$ for $\ell=0,1,2$ and
  (c),(d) still holds
\sn
\item[$(*)_2$]  \wilog \, there is $\nu_1 \in {}^{\omega >}\omega$
  such that $p_1 \Vdash ``q_1 $ % 2021-01-18 17:19 (t)$
  has trunk $\nu_1$".
\end{enumerate}
\mn
[Why?  Let $\mathbf G_1 * \mathbf G^1 \subseteq \bbP_1 * \name{\bbQ}_1$ be
generic over $\mathbf V$ such that $(p_1,\name q_1),(p_0,\name q_0) \in
\mathbf G_1$,
and $ p_3 \in \mathbb{P} _3 $ is a common upper bound of 
$ p_1,p_2 $.

We can find $\nu_1$ and $p'_1 \in \mathbf G_1$ above $p_1$
such that $p'_1 \Vdash_{\bbP_1} ``\tr(\name q_1) = \nu_1"$.  Let $q'_1
= q_1$ and choose $(p'_0,\name q'_0) \in \mathbf G_1 * \mathbf G^1$ above
$(p_0,\name q_0)$ such that $(p'_0,\name q'_0) \Vdash ``(p'_1,\name
q'_1) \in (\bbP_1 * \name{\bbQ}_1)/(\bbP_0 * \name{\bbQ}_0)"$.

Let
$(p'_2,\name q'_2) = (p_2,\name q_2)$, so clearly we are done.]
\mn
\begin{enumerate}
\item[$(*)_3$]  \wilog \, for some $\nu_2,p_2 \Vdash ``\tr(q'_2) =
  \nu_2 "$ (and $(*)_1$ still holds); we shall not repeat such
  statements.
\end{enumerate}
\mn
[Similarly as in the proof of $(*)_2$
because in the proof there    % 2020-12-25 13:18 
$ (p_2,\name q_2)$ was not
changed and we can interchange $\bbP_1,\bbP_2$.]
\mn
\begin{enumerate}
\item[$(*)_4$]  \wilog \, for some $\nu_0$ of length $\ge \ell
  g(\nu_1),\ell g(\nu_2)$ we have $p_0 \Vdash ``\tr(\name q_0) =
  \nu_0"$.
\end{enumerate}
\mn
[Why?  As we can just increase $(p_0,\name q_0)$, not change
$p_1,\name q_1,p_2,\name q_2$.]  % 2020-12-25 13:19 
\mn
\begin{enumerate} 
\item[$(*)_5$] \wilog \, $ \tr(\name{  q}_{\ell} )= \nu _0$ 
\end{enumerate} 
[Why? By the properties of $ \mathbb{Q} _{\ell} $.]

Now $ \name{ q_1 }, \name{ q_2}$ are two $ \mathbb{P} _3$-names
of members of $ \name{ \mathbb{Q}} _3$ with the same trunk
hence 
$ \Vdash _{\mathbb{P} _3} \lqq \name{ q_1 }, \name{ q_2}
$ are compatible$"$ so we are done.
% 2020-12-25 13:19 after (*)_4- the old proof  44444
% \begin{enumerate}

% \item[$(*)_5$]  from now on we don't change
% \sn
% \item[$(*)_6$]  \wilog \, $p_1$ forces a value $\rho_1$ for $\name\nu
%   \rest \ell g(\nu_3),\name\nu$ the generic of $\name{\bbQ}_1$.
% \end{enumerate}
% \mn
% [Why?  As in $(*)_2$ but we increase $p_1,p_0$ only $p'_1$ forces a
% value to $\name p_{\name q_0} \rest \ell g(\nu_3)$ says $\ell_1$ and
% is $\ge p_1$, necessarily $\nu_3 < \rho_1$.  Then we increase $p_0$ to
% $p'_0 \in \bbP_0$ not changing $\name q_0$ such that $p'_0
% \Vdash_{\bbP_0} ``p'_1 \in \bbP_1/\bbP_0"$.  So why $(p'_0,\name q_0)
% \Vdash ``(p'_1,\name q_1) \in (\bbP_1 * \name{\bbQ}_1)/(\bbP_0 *
% \name{\bbQ}_0)"$ by the properties of $\bbQ_0$.]
% \mn
% \begin{enumerate}
% \item[$(*)_7$]  \wilog \, $p_2$ forces a value $\rho_2$ to $\name\nu
%   \rest \ell g(\nu_3)$.
% \end{enumerate}
% \mn
% [Why?  Similarly to $(*)_6$ and $(*)_5$.]

% Now we can finish: we can find $p_3 \in \bbP_3$ above $p_0,p_1,p_2$
% and clearly $p_3 \Vdash_{\bbP_3} ``\{\name q_0,\name q_1,\name q_2\}$
% have a common upper bound.
% % 2020-12-03 17:37 \end{PROOF}
\bigskip

\noindent
\underline{Alternative presentation of the proof of \ref{ad.15}(1), self 
contained not depending on \ref{ad.13}, \ref{ad.14}}: 

We define an 
$\FSI$-template ${\gt}^\zeta = \mathfrak{t}[ \zeta ]$ from 
$ \mathbf{Q} _{\dom}$
for $\zeta \le \lambda$ by induction on $\zeta$.
\bigskip

\noindent
\underline{Case 1}:  For $\zeta = 0$. 

Let ${\gt}^\zeta$ be defined as follows:

\[
L^{\gt[\zeta]} = \mu  % 2020-12-03 16:56 
\]

\[
I^{\gt[\zeta]}_\alpha = \{A:A \subseteq \alpha\}  \text{ for $ \alpha < \mu$} % 2020-12-03 16:57 
\]

\bigskip

\noindent
\underline{Case 2}:  For $\zeta = \xi +1$. 

We choose ${\gt}^\zeta$ such that there is an isomorphism $\mathbf j_\zeta$
from $L^{{\gt}[\zeta]}$ onto $(L^{{\gt}[\xi]})^\kappa/D$, satisfying
$\mathbf j_\zeta \restriction L^{{\gt}[\xi]}$ is the canonical
embedding $\mathbf j_{D,{\gt}[\xi]}$, 
that  % 2020-12-03 17:02 and 
if $x \in L^{{\gt}[\zeta]}$ 
then  % 2020-12-03 17:02 ,
$\mathbf j_\zeta(x) = \langle x_\varepsilon:
\varepsilon < \kappa \rangle/D \in (L^{{\gt}[\xi]})^\kappa/D$ 
and:  % 2020-12-03 17:00 then:
$A \in I^{{\gt}[\zeta]}_x$ 
\Iff \,  
for some $\bar A = \langle
A_\varepsilon:\varepsilon < \kappa \rangle$ we have 
$A_\varepsilon \in I^{{\gt}[\xi]}_{x_\varepsilon}$ and 
$ \{  y : % 2020-12-03 17:03 \{\mathbf j_\zeta(y):
y \in A\} \subseteq \{\langle y_\varepsilon:
\varepsilon < \kappa \rangle/D:\{\varepsilon < \kappa:y_\varepsilon \in
A_\varepsilon\} \in D\}$.
\bigskip

\noindent
\underline{Case 3}:  $\zeta$ limit\footnote{
we may do one of    % 2020-12-25 13:34 
the following changes
(but not both): (a) in subcase 3B use $I^{{\gt}[\zeta]}_x =\{A$:
for some $\xi < \zeta,x \in
L^{{\gt}^\xi}$ and $A \in 
I^{\gt^\xi}_x\}$ and/or (b) in sub-case 3A behave as in sub-case 3B.
}.

We choose ${\gt}^\zeta$ as follows:

\[
L^{{\gt}[\zeta]} = \bigcup\limits_{\xi < \zeta} L^{{\gt}[\xi]}
\text{ as linear orders}.
\]

\mn
$I^{{\gt}[\zeta]}_x$ is
\bigskip

\noindent
\underline{Subcase 3A}:  If $x \in L^{{\gt}[0]}$ then
$\{A:A \subseteq \{s:L^{{\gt}[\zeta]} \models ``s < x"\}\}$.
\bigskip

\noindent
\underline{Subcase 3B}:  If $x \notin L^{{\gt}[0]}$ but $x \in
L^{\gt[\zeta]}$ then $I^{\gt[\zeta]}_x$ is 
\footnote{ 
so members of $L^{{\gt}[0]}$ have  % 2020-12-25 13:30 s  % 2020-12-06 07:57 
the ``veteranity privilege", 
i.e. ``founding father right"; members $t$ of $L^{{\gt}^0}$ 
have the maximal $I^{{\gt}[\zeta]}_t$.
}
(we rely on $ L^{\mathfrak{t} ^0}$  is well ordered):

\begin{equation*}
\begin{array}{clcr}
\{A: &\text{ for some } \xi < \zeta \text{ we have } x \in L^{{\gt}^\xi}
\text{ and if } y = \Min\{y \in L^{{\gt}^0}:L^{{\gt}[\zeta]} 
\models ``x < y"\} \\
  &\text{ which is } \in L^{{\gt}[0]} 
\text{ (and is always well defined see clause (b) of } 
 \oplus % 2020-12-03 17:27 (*) 
\text{ below) then} \\
  &A \backslash \{t \in L^{{\gt}[\zeta]}:L^{{\gt}[\zeta]} 
\models ``t < z" \text{ for some } z \text{ such that }
L^{{\gt}^0} \models ``z < y"\} \text{ belongs to } I^{{\gt}^\xi}_x \\
  &\text{(hence is } \subseteq L^{{\gt}^\xi})\}.
\end{array}
\end{equation*}

\mn
We now prove by induction on $\zeta \le \lambda$ that:
\mn
\begin{enumerate}
\item[$\oplus $]  $(a) \quad {\gt}^\zeta$ is an FSI-template
\sn
\item[${{}}$]  $(b) \quad L^{{\gt}[0]}$ is 
a cofinal % 2020-12-25 13:28 an unbounded 
subset of
$L^{{\gt}[\zeta]}$
\sn
\item[${{}}$]  $(c) \quad {\gt}^\zeta$ is smooth
\sn
\item[${{}}$]  $(d) \quad {\gt}^\xi \le_{\wk} {\gt}^\zeta$ for $\xi < \zeta$
\sn
\item[${{}}$]  $(e) \quad$ if $x \in L^{{\gt}[\zeta]}$
then $\{z:\text{for some } y \in L^{{\gt}^0}$ we have 
$L^{{\gt}[\zeta]} \models ``z \le y < x"\} \in I^{{\gt}[\zeta]}_x$
\sn
\item[${{}}$]  $(f) \quad L^{{\gt}[\zeta]}$ 
has cardinality $\le (\mu + |\zeta|)^\kappa$
% \sn
% \item[${{}}$]  $(g) \quad \langle {\gt}^\varepsilon:\varepsilon \le \zeta
% \rangle$ is $\le_{\wk}$-increasing continuous
\sn
\item[${{}}$]  $(g) \quad$    % 2020-12-06 07:59 h
we have ${\gt}^\zeta = \sum\limits_{\gamma < \mu}
{\gs}^{\gamma,\zeta}$ where ${\gs}^{\gamma,\zeta} = {\gt}^\zeta \rest 
\{x \in L^{{\gt}[\zeta]}:L^{{\gt}_\zeta} \models x < \gamma$ and 
$\beta \le x$ if $\beta  < % 2020-12-03 17:29 \in 
\gamma\}$.
\snn 
\item[${{}}$]  $(h) \quad$ the sequence $ \langle \mathfrak{s} ^{\gamma, \zeta }:
   \zeta \le \lambda \rangle $ is $\le_{\wk}$-increasing continuous.
\end{enumerate}
\mn
[Why?  Easy, e.g. why clauses (a)+(c) hold?  For $\zeta=0$ by
\ref{ad.10}(1). 
% 2020-12-25 13:38 +(6).
For $\zeta = \xi+1$ by \ref{ad.12}(2) 
noting that for $ t \in L^{\mathfrak{t} [0]}$ 
the desired value of $ I^ \mathfrak{t} _t$ holds. % 2020-12-25 13:37 
For
$\zeta$ limit, for any $t \in L^{{\gt}^0}$ clearly 
${\gs}^{\gamma,\zeta}$ is the union of the increasing continuous sequence
$\langle {\gs}^{\gamma,\varepsilon}:\varepsilon < \zeta \rangle$
hence is a smooth FSI-template by clause (h) and
\ref{ad.10}(5).  Now also ${\gt}^\zeta$ is a smooth FSI-template by
\ref{ad.10}(6).  
So $ \oplus $ holds indeed.]   % 2020-12-03 17:30 

Of course, we let $\bar K^\zeta = \langle K^\zeta_t:t \in
L^{\gt^\zeta} \rangle,K^\zeta_t = \emptyset$ and $\bbQ_t$ is the dominating
real forcing.

Lastly, let for $\zeta \le \lambda,\bbP_\zeta = \Lim_{\gt}
(\bar{\bbQ} \restriction L^{\gt^\zeta})$.
% 2020-12-03 17:31 Now
\mn
\begin{enumerate} 
\item[$\odot $] Now 
\begin{enumerate}
\item[$(\alpha)$]  $\bbP_\lambda$ is a c.c.c. forcing notion of
cardinality $\le \lambda^{\aleph_0}$ hence $\mathbf V^{\bbP_\lambda}
\models 2^{\aleph_0} \le \lambda$ by \ref{ad.4}(B)(j) as $\lambda =
\lambda^\kappa$ 
\sn
\item[$(\beta)$]  in $\mathbf V^{\bbP_\lambda}$ we have ${\gd} \le
\mu$, by \ref{ad.7}(1) applied with $R = <^*$ and $L^* = L^{\gt[0]}$ 
using $(*)(b)+(e)$
\sn
\item[$(\gamma)$]  in $\mathbf V^{\bbP_\lambda}$ we have ${\gb} \ge
\mu$ by \ref{ad.7}(2) applied with $R = <^*$
\sn
\item[$(\delta)$]  ${\gb} = {\gd} = \mu$ and ${\ga} \ge
\mu$ by $(\beta) + (\gamma)$ as it is well known that ${\gb} \le
{\gd}$ and ${\gb} \le {\ga}$.
\end{enumerate}
\end{enumerate} 
\mn
[Why? e.g. why clause $ ( \beta )$ holds? 
Applying \ref{ad.7}(1), we let $ R= < ^*, L^ *= 
L^{\mathfrak{t} [0]}$ and we have to verify clauses
(a)-(d) there. They are easy, e.g. for clause (c) there,
if $ A \subseteq L^\mathfrak{t} $ is $ \bar{ K } $-countable
then there is $ t \in L^* $ as promised
because $ L^{\mathfrak{t} [0]}$ is cofinal 
and is of order type $ \mu $ which is a 
regular uncountable cardinal.]

But in order to sort out the value of ${\ga}$ we intend to use
\ref{ad.7}(3) with $\theta$ there chosen as $\lambda$ here. 

But why the demand (c) from \ref{ad.7}(3) holds?  
Recall that every $ A \in L^{\mathfrak{t} [\zeta ]}$  
  is $ \bar{ K } $-closed. 
So assume $i(*)
\in [\kappa,\lambda)$ and $t_{i,n} \in L^{{\gt}^\lambda}$ for $i <
i(*),n < \omega$ be given.  As $\lambda$ is regular $> i(*)$,
necessarily for some $\xi < \lambda$ we have $\{t_{i,n}:i < i(*),n <
\omega\} \subseteq L^{{\gt}^\xi}$.  Now let $t_n \in L^{\mathfrak{t} ^{\xi
+1}}$ be such that $\mathbf j_{\xi +1}(t_n) = \langle t_{i,n}:i < \kappa
\rangle/D$; easily $\langle t_n:n < \omega \rangle$ is as required
(note that the number of isomorphism types of $\omega$-sequences
$\langle t_n:n < \omega \rangle$ in ${\gt}$ is trivially
\footnote{
  in fact, it is $\le 2^{\aleph_0}$ by the construction, but
  irrelevant here
}
$\le \beth_2$).  

So
\mn
\begin{enumerate}
\item[$(\varepsilon)$]  in $\mathbf V^{\bbP_\lambda}$ we have ${\ga}
\ge \kappa \Rightarrow {\ga} \ge \lambda$ by \ref{ad.7}(3), see there.

% 2021-02-10 05:13 (Xajob  jimvoni lamah $ \mathfrak{a} \ge \kappa $?)

\end{enumerate}
\mn
We are assuming $ \kappa \le \mu $ and by 
  $ \odot ( \gamma ) $ we have $ \mu \le \mathfrak{b} $ 
  and always $ \mathfrak{b} \le \mathfrak{a} $  so
  together $ \kappa \le \mathfrak{a} $. 
  Recallin $ ( \varepsilon ) $   % 2021-02-10 05:16 Together
  we are done.

\bigskip 
\noindent 
2) % 2020-12-03 17:38  DEBT
We indicate how to adapt the second proof  of part (1). 
For $ \partial $ a regular uncountable cardinal 
we consider only $ \mathbf{q} \in \mathbf{Q}_\partial  ^{\cln}$ which mean: 

\begin{enumerate} 
\item[$\boxplus ^1_ \mathbf{q}  $] Let $ \mathbf{q} \in \mathbf{Q} _{\cln}$  
mean
\begin{enumerate}
\item[(a)] $ \mathbf{q} \in \mathbf{Q}_\fsi $
\item[(b)] $ \partial (\mathbf{q} ) \le \partial $ 
    % 2020-12-25 13:52 {\aleph_1} $
\item[(c)] for every $ t \in L^ \mathbf{q} $ one of the following 
  occurs
  \begin{enumerate} 
  \item[($\alpha $)] $ K^ \mathbf{q} _t = \emptyset $ and $ \mathbb{Q} ^ \mathbf{q} _ t$ 
     is dominating real forcing= Hechler forcing
   \item[($\beta $)]  $K^\mathbf{q} _ t $  has cardinality $ < \partial $ 
     and $ I^\mathbf{q} _t = {\mathscr P} ( K^ \mathbf{q} _t )  $  and 
     $ \mathbb{Q} ^ \mathbf{q} _ {t, \name{ \eta }}$
     is  an explicitly linked $ (< \partial )$-forcing notion
      with universe $ \gamma ^\mathbf{q} _t < \partial $; 
    see below
 \end{enumerate} 
\end{enumerate} 
\end{enumerate} 
  
  Where 
  
  \begin{enumerate} 
  \item[$ \boxplus _2$]  We say that % 2020-12-06 08:05 
  the forcing notion $ \mathbb{Q} $ 
  is  an explicitly linked $ (< \partial )$-forcing notion with universe $ \gamma $ 
  when:  % 2020-12-06 08:06 
  \begin{enumerate}  
  \item[(a)] the set of elements of $ \mathbb{Q} $ is the ordinal $ \gamma $
  \item[(b)] for each $ n < \omega $ the set
    $ \{\omega \alpha + n : \omega \alpha + n < \gamma  \} $
     is a set of pairwise compatible elements of $ \mathbb{Q} $
  \end{enumerate} 
  \end{enumerate} 
  
  Next 
  \begin{enumerate} 
  \item[$\boxplus _3$]  the relevant claims  
  \ref{ad.5}-\ref{ad.10}  apply % 2020-12-25 13:54   hold 
for all
  $ \mathbf{q} \in \mathbf{Q} ^ {\cln}_ \partial $ with minor changes; mainly recalling $ \boxplus ^1_\mathbf{q} (d)(\beta )$.
  % 2020-12-25 13:53 that is   BDOQ XOB  to FILL?  xajob 2020-12-25 13:56  bevecem 5uly raq ccc 
  \end{enumerate} 
  
%     \begin{enumerate} 
%   \item[$\boxplus _4$ ] 
 We choose  $ \mathfrak{t} ^ \zeta, \langle\mathfrak{s} ^{\alpha, \zeta }:
     \alpha < \mu \rangle $     % 2020-12-06 08:08 
     % 2020-12-06 08:07 $ \mathbf{q} ^ \zeta $ 
 by induction on $ \zeta \le \lambda $  % 2020-12-06 08:08 
 as we have defined $ \mathfrak{t} ^ \zeta $ in the second proof 
 of part (1), but the second case splits to two, that is:
%  \begin{enumerate} 
% 2020-12-03 18:18  \item[(*)]
 
 \underline{Case 1}   $ \zeta = 0 $ 
 
 As above
 % 2020-12-03 18:18 \item[(*)] 
 
 \underline{Case 2} $ \zeta = \xi + 1 $  and $ \xi $  is even.
 
 As in the successor case above
 % 2020-12-03 18:19 \item[(*)] 
 
 \underline{Case 3} $ \zeta $ is a limit ordinal
 
 As above
 
 % 2020-12-03 18:19 \item[(*)] 
 \underline{Case 4}  $ \zeta = \xi + 1 $  and $ \xi $ is odd
 
%  For each $ \alpha < \mu $ let 
%  $ \mathbf{S}^1 _{\xi, \alpha }  = \{ (\xi, \varepsilon): \varepsilon < 
%      (\mu + |\xi |)^ \kappa $ 
%       2020-12-06 08:10 L, \gamma, \name{ \mathbf{Q}} ):
%  L $ is a $ K^{\mathfrak{t} [\xi ]}$-closed  subset of $ L^{s[\xi, \alpha ]}$  
%  of cardinality $ < \partial , \gamma < \partial $  and  
% $\name{  \mathbb{Q}} $ 
%  is a canonical $ \Lim _ {\mathbf{q} [\zeta ]}(\mathbb{Q} ^{\mathbf{q} [\partial ]}
%   \upharpoonright L)$-name of a forcing notion as in $ \boxplus _2$  with universe $ \gamma 
%  \} $.
 
 Now let us define $ \mathbf{q} _ \zeta $. We let
 \begin{enumerate} 
 \item[$ \odot$] 
 \begin{enumerate} 
 \item[(a)] 
 $ L^{\mathbf{q} [\zeta ]}= L^{\mathbf{q}  [\xi ]}\cup \{ (\mathbf{q} ^ \xi ,
 \alpha, \varepsilon ): \varepsilon 
   < (\mu + |\xi |)^ \kappa \} $ 
 and the order is defined by (in addition to the old order)
 % 2020-12-25 13:58 thar 
 \begin{enumerate} 
 \item[$(\alpha )$] 
 $ t = ( \mathbf{q} ^ \xi , \alpha, \varepsilon) $ 
 is below $ \alpha + 1 \in \mu = L^{\mathfrak{t} [0]}$
 above $ \alpha $
 \item[$(\beta )$] moreover $ t $ is above any 
 $ s \in  L^{\mathfrak{t} [ \xi] }  $ which is 
 below $\alpha + 1 $
 \item[$ (\gamma )$]  we let 
 $(\mathbf{q} ^ \xi , \alpha _1), \varepsilon _1) < 
    (\mathbf{q} ^ \xi , \alpha _2), \varepsilon _2)$ 
    iff $(\alpha _ 1 <  \alpha _2 )\vee (\alpha _1 = \alpha _2 
    \wedge \varepsilon _1 < \varepsilon _2 ) $
 \end{enumerate} 
%  % 2020-12-06 08:16 L, \name{ \mathbb{Q} } i): 
%   (L, \gamma, \name{ \mathbf{Q}} ) \in \mathbf{S} _{\xi , \alpha }, 
%  \alpha $ successor $ \} $
%  etc.   
%  % 2020-12-03 18:19 \end{enumerate} 
%  BDOQ XOB FILL 
% 2020-12-03 18:33   \end{enumerate} 
   \item[(b)] $ \mathfrak{t} ^\zeta , \langle \mathfrak{s} ^{\alpha, \zeta }
     : \alpha < \mu \rangle $  are as in $\oplus $ above
     \item[(c)] We define 
  $ \mathbf{q} ^ \zeta $ by induction on $ \zeta \le \lambda $. The new point is when 
  $ \zeta = \xi + 1 , \xi $ odd. 
  \item[(d)] In this case for $  \alpha < \mu $ 
      %   2020-12-06 08:10 L, 
      let $ \langle ( \xi, \alpha,\gamma, L_ {\xi, \alpha, \varepsilon}, \name{ \mathbb{Q} }_{\xi, \alpha, \varepsilon })
      : \varepsilon < ( \mu + |\xi |)^\kappa \rangle $ list the quintuples
      $ ( \xi, \alpha,\gamma, L, \name{ \mathbb{Q} })$
      % % 2020-12-06 08:27 \gamma, \name{ \mathbf{Q}} ):
      such that 
 $ L $ is a $ \bar{ K}^{\mathfrak{t} [\xi ]}$-closed  subset of $ L^{\mathfrak{s} [\xi, \alpha ]}$  
 of cardinality $ < \partial , \gamma < \partial $  and  
$\name{  \mathbb{Q}} $ 
 is a canonical $ \Lim _ {\mathbf{q} [\zeta ]}(\mathbb{Q} ^{\mathbf{q} [\partial ]}
  \upharpoonright L)$-name of a forcing notion as in $ \boxplus _2$  with universe
$ \gamma  
 \} $.
 \item[(e)]  lastly, if $ t = ( \mathbf{q} ^ \xi , \alpha, \varepsilon )$ 
 we let $ K_t = L_{\xi, \alpha, \varepsilon}$ and 
 $ \name{ \mathbb{Q}} ^ {\mathbf{q} [\xi ]}_t =  
 \name{ \mathbb{Q} }_{\xi, \alpha, \varepsilon }$ 
  \end{enumerate}
    \end{enumerate} 
    
    The rest should be clear.
  \end{PROOF}

% \end{proof} 
\newpage

\section {Eliminating the measurable} 
\bigskip

Without a measurable cardinal our problem is to verify condition (c) in
\ref{ad.7}(3).  Toward this it is helpful to show that for some 
$\aleph_1$-complete filter $D$ on $\kappa$, for any $i(*) \in 
[\kappa,\lambda)$ and $t_{i,\exn} \in L^{\gt}$, for $i < i(*),
\exn <  \sigma $, % 2020-12-01 14:10 \omega$, 
we have: for some $B \in D^+$ for every $j < i(*)$ 
some $A \in D$ satisfies: 
$ \lqq $for any $i_0,i_1 \in A \cap B$, the mapping
$t_{j,\exn} \mapsto t_{j,\exn}$; $t_{i_0,\exn} \mapsto t_{i_1,\exn}$ is a partial
isomorphism of ${\gt}$".  So $D$ behaves as an $\aleph_1$-complete
ultrafilter for our purpose. 

 [If you know enough model theory, this is the problem of finding
convergent sequences, see 
  % 2021-02-10 05:17 \cite[Ch.I,\S2, II]{Sh:300};  (and later 
  \cite[Ch.I,\S2, II]{Sh:c}, 
\cite[\S2]{Sh:300a}, \cite{Sh:300b}). 
The later  % 2021-02-10 05:19 This
had generalize what we know on 
% 2020-12-27 11:32 for 
stable first order $T$ with $\kappa = \kappa_r (T)$ 
(see \cite[Ch.II]{Sh:c} $ \kappa $ is regular and $ \le |T|^+$)
any indiscernible
 sequence (equivalently set) 
$\langle \bar a_\alpha:\alpha < \alpha^* \rangle$ of cardinality 
$\ge \kappa$, is convergent; why? as for any 
$\bar{\mathbf{b}} \in {}^{\kappa >} {\gC}$, for all but $< \kappa$ ordinals 
$\alpha < \alpha^*,\bar{\mathbf b} \char 94 \bar a_\alpha$ has a 
fixed type so average is definable.  
% 2021-02-10 05:23 In
% \cite[Ch.I,\S2, II]{Sh:300}; 
 % (=\cite[\S2]{Sh:300a}, \cite{Sh:300b}),  % 2020-12-01 14:17 \cite[Ch.II]{Sh:300}, 
 % 2021-02-10 05:24 In  we  deal with it in general,  
  The present is closed to \cite{Sh:72}, \cite{Sh:482}.
(The general case is % 2021-02-10 05:23 so 
harder to prove existence
which we do there under the relevant assumptions).]

 \bigskip

\begin{claim}
\label{abd.1}
Assume $2^{\aleph_0} < \mu = \cf(\mu) < \lambda = \cf(\lambda) 
= \lambda^{\aleph_0}$.  \Then \, for some $\bbP$ we have
\mn
\begin{enumerate}
\item[$(a)$]  $\bbP$ is a c.c.c. forcing notion of cardinality $\lambda$
\sn
\item[$(b)$]  in $\mathbf V^{\bbP}$ we have ${\gb} = {\gd} = \mu$ and
${\ga} = 2^{\aleph_0} = \lambda$.
\end{enumerate}
\end{claim}
 
 \begin{remark} 
\label{abd.1a}
About combining \ref{abd.1} with the end of \S2,
that is adding   $ \partial = (2^ \sigma )^+< \mu   $ 
and getting also $ \sigma \le \mathfrak{s} \le \mu $ 
(and even $ \mathfrak{s} < \mu $) see \cite{Sh:F2009},
\cite{Sh:F2029}
and more)
and \cite{Sh:F2023}.
 \end{remark} 
 
\begin{PROOF}{\ref{abd.1}} 
We rely on \ref{ad.5} + \ref{ad.7}.  Let $L^+_0$ be a
linear order isomorphic to $\lambda$, let $L^-_0$ be a linear order
anti-isomorphic to $\lambda$ (and $L^-_0 \cap L^+_0 = \emptyset$) and let
$L_0 = L^-_0 + L^+_0$.

Let $\mathbf J$ be the following linear order:
\mn
\begin{enumerate}
\item[$(a)$]  its set of elements is ${}^{\omega >}(L_0)$
\sn
\item[$(b)$]  the order is: $\eta <_{\mathbf J} \nu$ \Iff \, for
some $n < \omega$ we have $\eta \restriction n = \nu \restriction n$ and 
$\ell g(\eta) = n \wedge  \nu(n) \in L^+_0$ or $\ell g(\nu) = n \wedge 
\eta(n) \in L^-_0$ or we have
$\ell g(\eta) > n \wedge  \ell g(\nu) > n \wedge  L_0 \models \eta(n) < \nu(n)$.
\end{enumerate}
\mn
[See more on such orders Laver \cite{Lv71}  % 2020-12-01 13:27 , \cite{Sh:220,AP}
and
\cite[\S2]{Sh:E62}, \cite[\S5]{Sh:734}  % 2020-12-01 14:19 \cite[Ch.XIII,\S2]{Sh:e},
but we are self contained.] 

Note that
\mn
\begin{enumerate}
\item[$\boxdot_1$]  every interval of $\mathbf J$ as well as $\mathbf J$ 
itself has cardinality $\lambda$
\sn
\item[$\boxdot^+_1$]  if $\aleph_0 < \theta = \text{ cf}(\theta) <
\lambda$ or $\theta = 1$ or $\theta = 0$ and 
$\langle t_i:i < \theta \rangle$ is a strictly decreasing
sequence in $\mathbf J$ then $\mathbf J \restriction \{y \in 
\mathbf J:(\forall i < \theta)(y <_{\mathbf J} t_i)\}$ has cofinality
$\lambda$ if it is non-empty
\sn
\item[$\boxdot^-_1$]  the inverse of $\mathbf J$ satisfies
$\boxdot^+_1$, moreover is isomorphic to $\mathbf J$
\sn
\item[$\boxdot_2$]  if $\theta = \cf(\theta) > \aleph_0$ and
$s_\alpha,t_\alpha \in \mathbf J$ for $\alpha < \theta$ then we can find
a function $f:\theta \rightarrow \theta$ which is regressive and a
club $E$ of $\theta$ such that: if $\alpha_\ell < \beta_\ell$ are from
$E$ for $\ell=1,2$ and $f(\alpha_1) = f(\beta_1) = f(\alpha_2) =
f(\beta_2)$ then: $t_{\alpha_1} <_{\mathbf J} s_{\beta_1} \Leftrightarrow
t_{\alpha_2} <_{\mathbf J} s_{\beta_2}$ and $t_{\alpha_1} = s_{\beta_1}
\Leftrightarrow t_{\alpha_2} = s_{\beta_2}$ 
(we can add $t_{\alpha_1} <_{\mathbf J} t_{\beta_1} \Leftrightarrow
t_{\alpha_2} <_{\mathbf J} t_{\beta_2}$, etc., but this can be deduced
using the above several times).
\end{enumerate}
\mn
We now define by induction on 
$\zeta < \mu$ an $\FSI$-templates ${\gt}_\zeta 
=  \mathfrak{t} [\zeta ]  $ % 2020-12-01 14:33  
 such that
\mn
\begin{enumerate}
\item[$\odot   % 2020-12-01 14:32  (*)
^1_\zeta$] 
the set of members of $L^{\gt_\zeta}$
is a set of finite sequences starting with $\zeta$ hence disjoint to 
$ L^{\mathfrak{t} [\varepsilon ]} $ 
        % 2020-12-27 11:57 \zeta ]}  $ % 2020-12-01 14:34 $ {\gt}_\varepsilon$ 
for $\varepsilon < \zeta$; for $x \in 
L^{\gt[\zeta]}$ let $\xi(x) = \zeta$.  % 2020-12-01 14:34 ^ \zeta 
\end{enumerate}
% 2020-12-01 14:35 \end{PROOF}
\bigskip

\noindent
\underline{Defining ${\gt}_\zeta$}: 
\smallskip

\noindent
\underline{Case 1}:  $\zeta = 0$ or 
$\zeta$ successor or $\cf(\zeta) = \aleph_0$.

\begin{enumerate}
\item[$ \odot_2 $]   % 2020-12-04 08:39  \item[$(*)_2$]
Let $L^{{\gt}[\zeta]} = \{\langle \zeta \rangle\}$ and 
$I^{{\gt}[\zeta]}_{\langle \zeta \rangle   } = \{\emptyset\}$. % 2020-12-01 14:54 <>
\end{enumerate} 

\bigskip

\noindent
\underline{Case 2}:  $\cf(\zeta) > {\aleph_0} $ % 2020-12-27 12:17 \sigma $ . % 2020-12-01 14:54 \aleph_0$.

First

\begin{enumerate}
\item[$ \odot_3 $]   % 2020-12-04 08:40 \item[$(*)_3$]
Let $h_\zeta:\mathbf J \rightarrow \zeta$ be a function such that:
$\varepsilon < \zeta \Rightarrow h^{-1}_\zeta\{\varepsilon\}$ is a
dense subset of $\mathbf J$, 
  specifically $\nu = \eta   
    % 2020-12-27 12:18 \nu =% 2020-12-01 13:13 {}   % 2020-12-01 13:11 ^
  \char 94 \langle s \rangle  \in \mathbf
J \wedge (\otp(L^+_0 \rest \{ t : t % 2021-02-10 05:27 s 
< _{L^+_0}s \} % 2020-12-27 12:19  (L^+_0)_s
) = i < \zeta 
% 2020-12-27 12:20 
\vee \otp$(the  inverse of $(L^-_0)_s) = i < \zeta) 
\Rightarrow h_\zeta(\eta \, % 2020-12-01 13:14 {}^  
\char 94 \langle s \rangle )=i$
and otherwise $h_\zeta(\nu ) =0 $, % 2020-12-27 12:21 \eta)=0$.  
 Let $ h( \zeta, \eta ) = h_ \zeta (\eta )$ for $ \eta \in \mathbf{J} $.
\end{enumerate} 

Second 

\begin{enumerate}
\item[$ \odot_4 $]    % 2020-12-04 08:40 \item[$(*)_4$]
The set of elements of ${\gt}_\zeta$, that is of $ L^{\mathfrak{t} [\zeta ]}$ % 2020-12-01 15:19 
is

\[
\{\langle \zeta \rangle\} \cup \{\langle \zeta \rangle \char 94 
\langle \eta \rangle \char 94 x:\eta \in \mathbf J \text{ and } x \in
\bigcup\limits_{\varepsilon \le h_\zeta(\eta)} L^{{\gt}_\varepsilon}\}.
\]
\end{enumerate}

\mn
Third 

\begin{enumerate}
\item[$ \odot_5 $]    % 2020-12-04 08:40 \item[$(*)_5$] 
The order $<_{{\gt}_\zeta}$ defined by
$\langle \zeta \rangle \text{ is maximal}$  and:  % 2020-12-27 12:52  

% 2020-12-01 14:59 \begin{equation*}
% 2020-12-01 14:59 \begin{array}{clcr}
$ \langle \zeta \rangle {\char 94} \langle \eta_1 \rangle 
{\char 94} x_1 <_{{\gt}[\zeta]}
\langle \zeta \rangle {\char 94} \langle \eta_2 \rangle 
{\char 94} x_2$  
\Iff \,  at least one of the following holds: 
\begin{enumerate} 
\item[(a)] $\eta_1
<_{\mathbf J} \eta_2  $    % 2020-12-01 15:01 \vee (
\item[(b)] $ \eta_1 = \eta_2  \wedge  \xi(x_1) < \xi(x_2) $  
% 2020-12-01 15:01 \vee 
\item[(c)]  $ (\eta_1 = \eta_2 \\
  \wedge  \xi(x_1) = \xi(x_2) \wedge  x_1 <_{{\mathfrak t}_{\xi(x_1)}} x_2).$ 
% 2020-12-01 14:59 \end{array}
% 2020-12-01 14:59 \end{equation*}
\end{enumerate} 
\end{enumerate} 
\mn

Lastly, 

\begin{enumerate}
\item[$(*)_1 $]    % 2020-12-04 08:41  \item[$ \odot _1$ ] 
for $y \in L^{{\gt}[\zeta]}$ we define the ideal $I = 
I^{{\gt}[\zeta]}_y$:
% 2020-12-01 15:27 \end{enumerate} 
\mn

\begin{enumerate}
\item[$(\alpha)$]  if $y = \langle \zeta \rangle$ then 
$I = \bigr\{ Y:Y \subseteq L^{{\gt}[\zeta]} \backslash 
\{ \langle \zeta \rangle\}\}$
\sn
\item[$(\beta)$]  if $y = \langle \zeta \rangle \char 94 \langle
\nu \rangle \char 94 x$, \then \, $I$ is the family of countable
sets $Y$ satisfying the following conditions: 
\sn
\begin{enumerate}
\item[$(i)$]  $Y \subseteq L^{{\gt}[\zeta]}$ 
\sn
\item[$(ii)$]  $(\forall z \in Y)(z <_{{\gt}[\zeta]} y)$ 
\sn
% \item[$(iii)$]  for every $A \subseteq L^{{\gt}_{h_\zeta(v)}}$ we have 
% $A \in I^{{\gt}_{h_\zeta(v)+1}} \Leftrightarrow \{ \langle \zeta \rangle 
% ?????    % 2020-12-01 13:24 <>
%   % 2020-12-01 13:14 {}^\ 
% \char 94 \langle \nu \rangle  % 2020-12-01 13:25 <> % 2020-12-01 13:24 {}^\ 
% \char 94 y:y \in A\} \in I_y$
% \sn
\item[$(iii)$]   % 2020-12-01 iv
the set $\{\eta \in \mathbf J:(\exists x)(\langle 
\zeta \rangle \char 94 \langle \eta \rangle \char 94 x \in Y)\}$ is finite.
\item[$(iv)$]  if $ \nu < _ \mathbf{J} \eta $  and 
  $ z \in L^{\mathfrak{t} [h(\nu )]} $ % 2020-12-27 12:57 ( \zeta, \nu )}$ 
the $ \langle \zeta \rangle {\char 94} \langle \eta \rangle {\char 94} z
  \notin Y $
\item[$(v)$]  if $ \eta \le _ \mathbf{J} \nu $ then the set 
$ \{ z \in  L^{\mathfrak{t}} :  \langle \zeta \rangle {\char 94} 
  \langle \eta \rangle {\char 94} z  \in Y % 2020-12-03 17:20 
   \} $  belongs to 
   $ I^{\mathfrak{t} [h(\zeta , \eta ) ]}_x $
\end{enumerate}
\end{enumerate} 
\end{enumerate} 
% 2020-12-01 15:29 \end{enumerate} 
\mn
Why is ${\gt}_\zeta$ really an  % 2021-02-10 05:28 
$\FSI$-template?  We prove, of course, by
induction on $\zeta$ that:
\mn
\begin{enumerate}
\item[$(*)^2_\zeta$]  $(i) \quad L^{{\gt}_\zeta}$ is a linear order
\sn
\item[${{}}$]  $(ii) \quad I^{{\gt}_\zeta}_t$ is an ideal of subsets of
$\{s \in I^{{\gt}_\zeta}_t:s < t\}$
\sn
\item[${{}}$]  $(iii) \quad {\gt}_\zeta$ is an FSI-template,
\sn
\item[${{}}$]  $(iv) \quad {\gt}_\zeta$ is disjoint to 
${\gt}_\varepsilon$ for $\varepsilon < \zeta$
\sn
% \item[${{}}$]  $(v) \quad \cP^{\gt_\zeta}$ is a family
% of subsets of $L^{{\gt}_\zeta}$ closed under intersections such that

% \hskip25pt  each $L^{{\gt}_\zeta}_{<t}$ belongs to 
% $c \ell({\cP}^{{\gt}_\zeta})$.
\end{enumerate}
\mn
[Why?  By \ref{ad.10}(8) and looking at the definitions.] 

Next we prove by induction on $\zeta$, that ${\gt}_\zeta$ is a 
smooth   % 2021-02-10 05:28  sst 
$\FSI$-template. Arriving at % 2021-02-10 05:29 
$\zeta$
\mn
\begin{enumerate}
\item[$(*)^3_\zeta$]  for 
$\eta \in \mathbf J$ and $\varepsilon \le h_\zeta(\eta)+1$, we
have ${\gt}_\zeta \restriction \{\langle \zeta \rangle \char 94
\langle \eta \rangle \char 94 \rho:\rho \in 
\bigcup\limits_{\xi < \varepsilon} {\gt}_\xi\}$ is a 
smooth$ \, \FSI$-template. 
\end{enumerate}
\mn
[Why?  We prove this 
by induction on $\varepsilon$; for $\varepsilon =0$ by 
\ref{ad.10}(1), for $\varepsilon$ successor by \ref{ad.10}(3) 
for $\varepsilon$ limit by \ref{ad.10}(5) and \ref{ad.10}(6).]
\mn
\begin{enumerate}
\item[$(*)^4_\zeta$]  for $Z \subseteq \mathbf J$ we have 
${\gt}_\zeta \restriction (\bigcup\limits_{\eta \in Z} \{\langle \zeta
\rangle \char 94 \langle \eta \rangle \char 94 \rho:\rho \in 
\bigcup\limits_{\xi < h_\zeta(\eta)} \gt_\xi\})$ is a 
smooth$ \, \FSI$-template. 
\end{enumerate}
\mn
[Why?  By induction on $|Z|$, for $|Z| = 0,|Z| = n+1$ by \ref{ad.10}(3),
for $|Z| \ge \aleph_0$ by \ref{ad.10}(5).]
\mn
\begin{enumerate}
\item[$(*)^5_\zeta$]   ${\gt}_\zeta \restriction (L^{{\gt}_\zeta}
\backslash \{ \langle \zeta \rangle\})$ is a smooth$ \, \FSI$-template.
\end{enumerate}
\mn
[Why?  By $(*)^4_\zeta$ for $Z = \mathbf J$.]
\mn
\begin{enumerate}
\item[$(*)^6_\zeta$]  ${\gt}_\zeta$ is a smooth % 2020-12-27 13:00 smooth$^+$ 
FSI-template.
\end{enumerate}
\mn
[Why?  By \ref{ad.10}(3).]
\mn
\begin{enumerate}
\item[$(*)^7_\zeta$]   if $K \subseteq L^{\gt_\zeta}$ is
countable and $t \in L^{\gt_\zeta}$ then the ideal $I^{\gt_\zeta}_t \cap
{\cP}(K)$ is generated by a countable family of subsets of $K$.
\end{enumerate}
\mn
[Why?  Check by induction on $\zeta$.]

\medskip 

\mn
Now for $\zeta \le \mu$ let
\mn
\begin{enumerate}
\item[$(*)^8_\zeta $] % 2020-12-01 15:31 777777  \boxdot_2$]   
${\gs}_\zeta =: \sum\limits_{\varepsilon < \zeta}
{\gt}_\varepsilon$, i.e.
\sn
\begin{enumerate}
\item[$(i)$]  the set of elements of ${\gs}_\zeta$ is
$\bigcup\limits_{\varepsilon < \zeta} L^{{\gt}_\varepsilon}$
\sn
\item[$(ii)$]  for $x,y \in {\gs}_\zeta$ we have 
$x <_{{\gs}_\zeta} y$ iff $\xi(x) < \xi(y) \vee (\xi(x) = \xi(y) \wedge 
x <_{{\gt}_{\xi(x)}} y)$
\sn
\item[$(iii)$]  $I^{{\gs}_\zeta}_y = \{Y \subseteq ^{{\gs}[\zeta]}:
(\forall z \in Y)(z <_{{\gs}_\zeta} y)$ and 
$\{z \in Y:  % 2020-12-27 13:04 {\gs}_\zeta:
\xi(z) = \xi(y) % 2020-12-27 13:05 $   % 2020-12-27 13:04 and $\hskip80pt  z \in Y
\} \in I^{{\gt}[{\xi(z)}]}_y\}$
% \sn
% \item[$(iv)$]   $P^{{\gs}_\zeta}$ is as in \ref{ad.9}(3).
\end{enumerate}
\end{enumerate}   % 2020-12-04 
\mn  % 2020-12-04 08:45 \sn
\begin{enumerate} 
\item[$(*)^9  % 2020-12-01 15:33 8
_\zeta$]   ${\gs}_\zeta$ is a smooth$ \,
  \FSI$-template. 
\end{enumerate}  
\mn
[Why?  Just easier than the proof above.]
\mn
\begin{enumerate}
\item[$(*)^{10} % 2020-12-01 15:33 }9
_\zeta$]   if $K \subseteq L^{{\gs}_\zeta}$ is
countable and $t \in L^{{\gs}_\zeta}$, \then \, the ideal
$I^{{\gs}_\zeta}_t \cap {\cP}(K)$ of subsets of $K$ is generated
by a countable family of subsets of $K$.
\end{enumerate}
\mn
[Why?  By $(*)^7_\zeta$ and the definition of ${\gs}_\zeta$ and of
the ${\gt}_\varepsilon$-s.]
% 2020-12-04 08:47 \end{enumerate}
\bigskip

\noindent
Let \footnote{
   but if you like to avoid using $(*)^7_\zeta,(*)^{10} % 2020-12-01 15:34 9
    _\zeta$
  and ${\cW}$ below just use $\extheta  % 2020-12-27 13:23 
  = \beth^+_2$.  In fact even
  without $(*)^7_\zeta + (*)^  { 10 } % 2020-12-01 15:34 9
  _\zeta$ above, countable ${\cW}$ suffice
  but then we have to weaken the notion of isomorphisms, and no point.
}
  % 2020-12-01 15:55 
$ \sigma = {\aleph_0}, \extheta   % 2020-12-27 13:08 \theta 
= (2^{\sigma })^+$, we shall prove below by induction on 
$\zeta$ that ${\gs}_\zeta,{\gt}_\zeta$ are $(\lambda,\theta,\sigma  )$-good 
(see definition below and Sub-claim \ref{abd.2a}); 
% 2020-12-27 13:09 for part (2) we need $ ( \lambda , \theta, \sigma )$-good                                                       
% ?? qarah majehu?
  then 
we can finish the proof as in \ref{ad.15} 
    % 2020-12-27 13:09 using ${\gs}_\mu$ 
(and $ (*)^7_\zeta $  and $ (*)^{10}_\zeta $) 
% 2020-12-27 13:10 BDOQ LAMA NAPAL
% 2020-12-01 14:51  where    % 2020-12-01 15:10 $ 3523 $
\end{PROOF} 

\begin{definition}
\label{abd.2}
1) Assume 
\footnote{
 we ignore here $\bar K$ and $\{(t,\bar
 \varphi_t,\name \eta_t):t \in L^{\gt}\}$
 using the default values  % 2020-12-04 08:51 
     }  % 2021-02-10 05:30 sop hevarat julayim
$ \lambda \ge \extheta \ge \tau   % 2020-12-04 08:56 , \extheta 
> \sigma $,  
    % 2020-12-27 13:10 [[ BDOQ xob $ \theta  \vee \partial $ ]]
$\extheta$ is regular
uncountable and $(\forall \alpha < \extheta)[|\alpha|^{\sigma } < \extheta]$  
  % 2020-12-01 15:09 \end{proof} 
  and $ \mathbf{W} \subseteq {\mathscr P} ({\mathscr P} (\sigma ))$.
We say that a smooth $\FSI$-template ${\gt}$ is 
$(\lambda,\extheta,\tau, \sigma, \mathbf{W}  )$-good \If \, :
\mn
\begin{enumerate}
\item[$\oplus$]   assume that $t_{\alpha,\exn} \in L^{\gt}$ for 
$\alpha < \extheta , % 2020-12-27 13:20 \theta
 \exn < \sigma % 2020-12-27 13:21 \omega
,\{t_{\alpha,\exn}:\exn < \sigma   % 2020-12-04 09:01 \omega 
\}$ is $\bar K$-closed, 
\then \, we can find 
$ {\mathscr W } \in \mathbf{W} $  and % 2020-12-04 09:03 
a club $C$ of $ \partial$ % 2020-12-27 13:20 \theta$ 
and a pressing down 
function $h$ on $C$ such that:
\sn
\item[$\oplus'$]  if $S \subseteq C$ is stationary in  
    % 2020-12-27 13:18 $\theta,(\forall 
$ \partial $,  % 2021-02-10 05:34 $  , 
$ ( \forall  
\delta \in S)[\cf(\delta) > \sigma % 2020-12-04 09:01 \aleph_0
\wedge (\tau = \extheta \rightarrow \cf( \delta ) = \tau )]$ 
% 2020-12-27 13:18 BDOQ 
and $h \rest S$ is constant
\then \,:
\sn
\begin{enumerate}
\item[$\boxtimes^1_S$]   for every $\alpha < \beta$ in $S$ 
and $ w \in {\mathscr W } $, 
the truth value
of the following statements does not depend on $(\alpha,\beta)$: 
(but may depend on $\exn,\exm$ and $w \in {\cW}$) 
\sn
\item[${{}}$]  $(i) \quad t_{\alpha,\exn} = t_{\beta,\exm}$
\sn
\item[${{}}$]  $(ii) \quad t_{\alpha,\exn} <_{L^{\gt}} t_{\beta,\exm}$
\sn
\item[${{}}$]  $(iii) \qquad \{t_{\alpha,\exl }:\exl \in w\} \in 
I^{\gt}_{t_{\alpha,\exm}}$ 
\sn
\item[${{}}$]  $(iv) \quad \{t_{\beta,\exl}:\exl \in w\} \in 
I^{\gt}_{t_{\alpha,\exn}}$
\sn
\item[${{}}$]  $(v) \quad \{t_{\alpha,\exl}:\exl \in w\} \in 
I^{\gt}_{t_{\beta,\exn}}$
\sn
\item[$\boxtimes^2_S$]  let $\delta^* \le \extheta$ be such
that % 2020-12-27 13:17 
$\cf(\delta^*) = \tau$ and $\sup(S \cap \delta^*) = \delta^*$;
if $\extheta \le \beta^* < \lambda$ and
$s_{\beta,\exn} \in L^{\gt}$ for 
$\beta \in [ \partial ,  \beta^*)   % 2020-12-27 13:16  < \lambda
,\exn < \omega$ 
\then \, we can find $t_\exn \in L^{\gt}$ for $\exn < 
\omega$ such that for every $\beta < \beta^*$, for every large enough 
$\alpha \in S \cap \delta^*$ for some ${\gt}$-partial $\otimes$ 
isomorphism $f$ we have $f(t_\exn) = t_{\alpha,\exn},f(s_{\beta,\exn}) = s_{\beta,\exn}$.
\end{enumerate}
\end{enumerate}
\mn
2) We say ${\gt}$ is strongly $(\lambda,\theta,\tau, \sigma )$-good \If \,
above we have $\mathbf{W} = {\mathscr P} (  % 2020-12-27 13:15 
  {\mathscr P} ( \sigma )) $.

\mn 
3)
We may omit $ \mathbf{W} $ in part (1) when  
% 2020-12-04 09:15 allow 
$\mathbf{W}  =  \{{\mathscr W } : {\mathscr W } $ 
is % 2021-02-10 05:35 
an ideal
of the Boolean algebra ${\mathscr P} ( \sigma ) $ 
generated by $ \le \sigma $ sets$\} $  % 2020-12-04 09:11 {\cP}(   \omega)$ 
% 2020-12-04 09:17 (but  if $\theta > \beth_2(\sigma )$it is natural to  use $\mathbf{W} =\{ {\mathscr P} ( \sigma )  \} $ % 2020-12-04 09:14  {\mathscr P} ({\mathscr P} ( \sigma )).  % 2020-12-04 09:13 this is the same).  

\mn 
4)  % 2020-12-27 13:14 In  cases
Above we may omit $\tau$ if $\tau = \theta$.
\end{definition}

\begin{observation}
\label{abd.02.0}
In Def. \ref{abd.2},
instead ``$h$ regressive" it is enough to demand:
for some sequence $\langle X_\alpha:\alpha < \theta \rangle$ of sets, 
increasing continuous, $|X_\alpha| < \theta$ and for every (or club
of) $\delta < \theta$, if $\cf(\delta) > \aleph_0$ then $h(\delta) \in 
{\cH}_{< \aleph_1}(X_\delta)$.
\end{observation}

\begin{claim}
\label{abd.2a}
\noindent 
1)
In the proof of \ref{abd.1};
\mn
   % 2020-12-27 13:29 1)
\begin{enumerate}
\item[$(i)$]  ${\gt}_\zeta$ is strongly $(\lambda,\extheta, \sigma )$-good
\sn
\item[$(ii)$]  ${\gs}_\zeta$ is strongly $(\lambda,\extheta,\sigma % 2020-12-27 13:29    /\aleph_1
  )$-good
    % 2020-12-27 13:28 BDOQ
\sn
\item[$(iii)$]   if $\cf(\zeta) \ne \theta$ then 
${\gs}_\zeta$ is also strongly $(\lambda,\theta)$-good.
\end{enumerate}

\noindent 
2) Assume $ \lambda = \cf( \lambda) > \mu = \cf( \mu), \mathbf{J} ,
\bar{ \mathfrak{t} % 2021-02-10 05:37 t
} _ \varepsilon (\varepsilon <  \mu ), \mathfrak{s} _ \zeta 
(\zeta \le \mu )$  
are as in the proof of \ref{abd.1}. 
If $ \partial = (2^ \sigma )^+  < \mu $ 
then clauses (i), (ii), (iii) above hold.  % 2020-12-27 13:34 
% 2020-12-27 13:31 BDOQ  male5

\end{claim}
  
\begin{PROOF}{\ref{abd.2a}}
\noindent 
1)
Recall that $\partial % 2021-02-10 05:37 \theta 
= (2^\sigma    % 2020-12-04 09:21 {\aleph_0}
)^+$ (see before Definition \ref{abd.2}).

First note that

\begin{enumerate} 
\item[$(*)_1$]  for every $ \zeta < \mu $ there % 2021-02-10 05:38 \lambda $   there
    is a sequence $ \bar{  \varrho }_ \zeta $ such that:
\begin{enumerate}
  \item[(a)] $\bar{ \varrho }_\zeta = \langle  \varrho _{\zeta, s}: s \in L^{\mathfrak{t} [\zeta ]}\rangle $ 
  \item[(b)] $ \varrho _{\zeta, s }\in {}^{ \omega > } \lambda $ 
  \item[(c)]  the truth value of $ L^{\mathfrak{t} [\zeta ]}
     \models \lqq s < t " $ depends only on 
\begin{enumerate}
\item[$(\alpha )$] $ \lg ( \varrho _{\zeta, s}) $ 
\item[$(\beta )$] $\lg ( \varrho _{\zeta, t }) $ 
\item[$(\gamma )$] the truth values of 
$\varrho _{\zeta, s}(k) < \varrho _{\zeta, t}( {\ell} )$,
$\varrho _{\zeta, s}(k) = \varrho _{\zeta, t}({\ell} )$,
$\varrho _{\zeta, s}(k) > \varrho _{\zeta, t}({\ell} )$
for the relevant $k, {\ell} $
\end{enumerate} 
  \end{enumerate} 
\end{enumerate} 

[Why? Read the definition of $ \mathfrak{t} _ \zeta $.]

Second note that 

\begin{enumerate} 
\item[$(*)_2$]  there is a sequence  $ \bar{ \varrho }=
\langle \varrho _ s : s \in L^ {\mathfrak{s}[\mu ] }\rangle $
satisfying the parallel of $ (*)_1$.
\end{enumerate}

\noindent
 hence
 
 \begin{enumerate} 
 \item[$(*)_3$]  if $ \bar{ s } = \langle s_ \exn=  s( \exn): 
 \exn \le  \sigma \rangle \in  {}^{ \sigma + 1 }  % 2021-02-10 05:39 + 1 
 (L^ {\mathfrak{s} [\mu ]}) $ 
   then the truth value of $ \{s_ \exn: \exn < \sigma  \} \in I^{\mathfrak{s}[\mu ] } $
   depends only on   % 2021-02-10 07:44  s
   \begin{enumerate} 
 \item[(a)] $ \lg ( \varrho _ {s (\exn)})$ for $ \exn \le \sigma $
 \item[(b)] the truth values of 
 $\varrho _{\zeta, s( \varepsilon )}(k) < \varrho _{\zeta, s( \zeta )}  % 2021-02-10 05:41 
 ({\ell} )$,
$\varrho _{\zeta, s(\varepsilon)}(k) = \varrho _{\zeta, s( \zeta )}({\ell} )$,
$\varrho _{\zeta, s(\varepsilon ) }(k) > \varrho _{\zeta, s  % 2021-02-10 07:45  \sigma % 2021-02-10 05:42 s*
( \zeta )}({\ell})$ 
for $ \varepsilon, \zeta \le \sigma $ and relevant $ k, {\ell} $
   \end{enumerate} 
 \end{enumerate} 

\noindent 
{Why? Again look at the choice of $ \mathfrak{s} _ \mu $ }

\noindent 
Now, given $ \bar{ t } _ \alpha = \langle t_{\alpha ,\exn } = % 2021-02-10 07:45 =
    t[\alpha, \exn ]: \exn < \sigma % 2021-02-10 05:44 \theta 
    \rangle  \in
   {}^{ \sigma }(L^{\mathfrak{s} [\mu ]}) $ 
for $ \alpha  < \extheta  $
 \underline{define} 
 
 \begin{enumerate} 
 \item[$(*)_4$] 
 $ {\mathscr U } _ \alpha =
  \cup \{\Rang(\varrho _{t[\beta , \exn]}): \beta < \alpha , 
    \exn < \sigma \}  \cup % 2021-02-10 07:46 \exists 
    \{ \infty \}  $ 
\end{enumerate} 

Next % 2021-02-10 05:45 we define

\begin{enumerate}
\item[$(*)_5$]  we define the  function $ h % 2021-02-10 07:46 
      _ {\ell} , {\ell} = 0,1$ 
    % 2021-02-10 05:46 ,2$ from 
  with domain $ \extheta  \setminus \{ \emptyset \} $, so for 
  $ \alpha \in (0, \extheta  )$ we let:
  \begin{enumerate} 
  \item[(a)] $ h_0(\alpha )$ is equal to the set 
 as  $  \{ ( \exn, k , \exm, {\ell}  ):
  \varrho _{t[\alpha , \exn ]}(k) < 
  \varrho _{t[\alpha, \exm ]}({\ell}   % 2021-02-10 07:47 
  ) $ % 2021-02-10 05:47 \exm) $ 
  and both are well defined$\} 
    \} $  
\item[(b)] $h_1( \alpha ) $ is the minimal 
  non-zero member  $ \beta $ of 
$ {\mathscr U } _ \alpha $ such that (if there is no one then it is
zero):

for every  $ \exn < \sigma, k < \lg (t_{\alpha, \exn })  % 2021-02-10 05:52 
(k)$,
the following are equal:
\begin{enumerate} 
\item[$ (\alpha )$] 
the minimal member of $ {\mathscr U } _ \alpha $ which is 
$ > \varrho _{t[\alpha , \exn]}(k) $ 
% 2021-02-10 05:49 is equal  to 
\item[$ (\beta )$] 
the minimal member of $ {\mathscr U } _\beta  $ which is 
$ > \varrho _{t[\beta  , \exn]}(k)$,
\end{enumerate} 
\item[(*b)']    % 2021-02-10 07:49 and 
similarly for  $ \ge $ (and so for equal)
  \end{enumerate} 
\end{enumerate} 

Clearly 

\begin{enumerate}
\item[$(*)_6$] $ h_0$  has range of cardinality $ < \extheta $ 
and  $ h_1 $ is regressive 
\end{enumerate} 

\noindent 
Lastly

\begin{enumerate} 
\item[$(*)_7$]  if $ S \subseteq \extheta $ is stationary 
and  $ \delta \in S \Rightarrow \cf( \delta) > \sigma $
and $ h_0, h_1$  restricted  to  $ S $ are constant
then $ S$  is as required.
 \end{enumerate} 
        % 2021-02-10 07:50 \end{PROOF} 

\underline{OLD/ pre 2020  PROOF }

We prove this by induction on $\zeta$.
\bigskip

\noindent
\underline{For ${\gs}_\zeta$}:

If $\zeta=0$ it is empty.  Otherwise given $t_{\alpha,\exn} \in 
{\gs}_\zeta = \sum\limits_{\varepsilon < \zeta} {\gt}_\varepsilon$ 
for $\alpha < \theta,\exn < \omega$ let $h^*_0(\alpha)$ be the sequence 
consisting of:
\mn
\begin{enumerate}
\item[$(i)$]  $\xi_{\alpha,\exn} =: \Min\{\xi:\xi \in 
\{\xi(t_{\beta,\exm}):\beta < \delta,\exm < \omega\} 
\cup \{\infty\}$ and $\xi \ge \xi(t_{\alpha,\exn})\}$ for $\exn < \omega$ and 
\sn
\item[$(ii)$]  $u_\alpha = \{(\exn,\exm,\exl):\xi(t_{\alpha,\exn}) = 
\xi_{\alpha,\exm} \wedge  \exl=1$ or $\xi(t_{\alpha,\exn}) \le \xi(t_{\alpha,\exm})
\wedge  \exl =2\}$ and
\sn
\item[$(iii)$]   $\mathbf w_\alpha = \{(n,w):\exn < \omega,
w \subseteq \omega \text{ and }
\{t_{\alpha,\exm}:\exm \in w\} \in I^{\gt}_{t_{\alpha,\exn}}\}$ that is
$h^*_0(\alpha ) % 2021-02-10 05:53 \delta)
= \langle u_\alpha,\langle \xi_{\alpha,\exn}:\exn < \omega
\rangle,\mathbf w_\alpha \rangle$.
\end{enumerate}
\mn
If $S_y = \{\delta:\cf(\delta) \ge \aleph_1,h^*_0(\delta) = y\}$ is stationary
we define $h^*_1 \rest S_y$ such that it codes $h^*_0(\delta)$ and if
$\exn   % 2021-02-10 05:54 
(*) < \omega$ and the sequence $\langle \xi(t_{\alpha,\exn(*)}):\alpha
\in S_y \rangle$ is constant call it $\xi_{y,\exn(*)}$ let $u_{y,\exn(*)} = \{\exn:  % 2021-02-10 05:55 
\xi_{\alpha,\exn} = \xi_{y,\exn(*)}\}$, \then \, $h^*_1 \restriction S_y$ codes a
function witnessing the $(\lambda,\theta)$-goodness of 
${\gt}_{\xi_{y,\exn(*)}}$ for $\langle t_{\alpha,\exn}:\exn \in u_{y,\exn(*)},
\alpha \in S_y \rangle$.
\bigskip

\noindent
Fix $S$ as in $\oplus'$.
It is easy to check that this shows 
$\boxtimes^1_S$ even if $\cf(\zeta) = \theta$.
But assume $\cf(\zeta) \ne \theta \wedge  \delta^* = \theta$ or $\delta^*
< \theta,\cf(\delta^*) = \aleph_1$ (or just $\aleph_0 < 
\cf(\delta^*) < \theta$), $\delta^* = \sup(S \cap \delta^*)$; 
we shall prove also the statement from $\boxtimes^2_S$.  Let
$w_1 = \{\exn:  % 2021-02-10 05:57 
\text{the sequence }
\langle \xi(t_{\beta,\exn}):\beta \in S \rangle$ is strictly
increasing$\}$, $w_0 = \{\exn:  % 2021-02-10 05:56 \exn
\langle \xi(t_{\beta,\exn}):\beta \in S \rangle$
is constant$\}$, let $\xi(S,\exn) = \xi_{S,\exn} = \cup \{\xi(t_{\beta,\exn}):
\beta \in S\}$ as $\cf(\zeta) \ne \theta$ it is $< \zeta$ also when $\exn  % 2021-02-10 05:56 
\in w_1$.

Given $\langle \bar s_\beta:\beta < \beta^* \rangle,\beta^* < \lambda$
and $\bar s_\beta = \bar s = 
\langle s_{\beta,\exn}:\exn < \omega \rangle$ we have to find
$\langle t_\exn:\exn < \omega \rangle$ as required in $\boxtimes^2_S$.  If
$\exn   % 2021-02-10 05:58 
\in w_0$ let $w'_{0,\exn} = \{\exm   % 2021-02-10 05:58 
\in w_0:\xi(t_{\alpha,\exn}) =
\xi(t_{\alpha,\exm})$ for $\alpha \in S\}$ and to choose $\langle t_\exm:\exm
\in w'_{0,\exn} \rangle$ we use the induction hypothesis on 
$\gt_{\xi(S,\exn)}$.  If $\exn  % 2021-02-10 05:59 
\in w_1$ then we can find 
$t^*_\exn  % 2021-02-10 05:59 
\in {\gt}_{\xi_{S,\exn}}$ such that
$\{t:t \in {\gt}_{\xi_{S,\exn}},t \le_{{\gt}_{\xi(S,\exn)}} t^*\}$ is
disjoint to $\{t_{\beta,\exm}:\beta < \delta^*,\exm < \omega\} \cup 
\{s_{\beta,\exm}: \beta < \beta^*$ and $m < \omega\}$ this is possible
because the lower cofinality of
$L^{{\gt}_{\xi(S,\exn)}}$ is the same as that of $L_0$ and 
is $\lambda = \cf(\lambda) > \theta + |\beta^*|$.
Then we choose $\eta^* \in \mathbf J$ such that $(\forall x)
(\langle \zeta \rangle \char 94 \langle \eta^* \rangle 
\char 94  x \in {\gt}_{\xi(S,\exn)} 
\Rightarrow \langle \zeta \rangle \char 94 \langle \eta^* \rangle
\char 94 \langle x \rangle <_{{\gt}_{\xi(S,\exn)}} t^*)$ and we
choose together $\langle t_{\exn'}:  % 2021-02-10 06:01 
\exn' \in w_1,\xi_{S,\exn'} = \xi_{S,\exn} \rangle$
such that $t_\exn % 2021-02-10 06:01 
\in\{\langle \zeta \rangle \char 94 \langle \eta  % 2021-02-10 06:01 \rangle 
\rangle \char 94 \langle x \rangle \in {\gs}_\zeta:\eta <_{\mathbf J} 
\eta^*\}$ taking care of ${\cW}$, (inside $\{\exn % 2021-02-10 06:00 
\in w_1:
\xi(t_{\alpha,\exn}) = \xi_{S,\exm}\}$ and automatically for others,
i.e. considering $t_{\exn_1},t_{\exn_2}$ such that $\xi_{S,\exn_1} \ne
\xi_{S,\exn_2})$, this is immediate.
  % 2020-12-04 09:37 \end{PROOF}

\bigskip

\noindent
\underline{For ${\gt}_\zeta$}:

Similar (using $\boxdot_1 + \boxdot_2$).  
% 2020-12-27 15:48 
\end{PROOF}   % 2020-12-04 09:37 
\bigskip

\centerline {$* \qquad * \qquad *$}
\bigskip

We may like to have ``$2^{\aleph_0} = \lambda$ is singular", ${\ga} =
\lambda,{\gb} = {\gd} = \mu$.  Toward this we would like to have
a linear order $\mathbf J$ such that if $\bar x = \langle x_\alpha:\alpha <
\theta \rangle$ is monotonic, say decreasing then for any $\sigma < \lambda$
for some limit $\delta < \theta$ of uncountable cofinality the 
linear order $\{y \in \mathbf J:\alpha < \delta \Rightarrow y 
<_{\mathbf J} x_\alpha\}$ has cofinality $> \sigma$.  Moreover, 
$\delta$ can be chosen to suit $\omega$
such sequences $\bar x$ simultaneously.  So every set of $\omega$-tuples from
$\mathbf J$ of cardinality $\ge \theta$ but $< \lambda$ can be ``inflated".
\bigskip

\begin{lemma}
\label{abd.3}  
Assume
\mn
\begin{enumerate}
\item[$(a)$]  $(2^\sigma   % 2020-12-04 09:37 {\aleph_0}
)^+ < \mu = \cf(\mu) < \lambda = 
\lambda^\sigma  % 2020-12-04 10:10 {\aleph_0}
,\lambda$ singular
\sn
\item[$(b)$]  $(\forall \alpha < \mu)[|\alpha|^{\aleph_0} < \mu]$
\sn
\item[$(c)$]  $\mu \ge \aleph_{\cf(\lambda)}$ or at least
\sn
\item[$(c)^-$]  there is $f:\lambda \rightarrow \cf(\lambda)$
such that if $\langle \alpha_\varepsilon:\varepsilon < 
\mu \rangle  \in {}^{ \mu } \lambda $ 
is
(strictly) increasing continuous, $\alpha_\varepsilon < \lambda$ and $\gamma
< \cf(\lambda)$ \then \ for some $\varepsilon <  \mu $ % 2020-12-04 10:12  \tau$ 
we have
$f(\alpha_\varepsilon) \ge \gamma$.
\end{enumerate}
\mn
\Then \, for some c.c.c. forcing notion of cardinality $\lambda$ we have
$\Vdash_{\bbP} ``2^{\aleph_0} = \lambda,{\gb} = \gd = \mu   % 2020-12-04 10:12 \kappa
,\ga = \lambda"$.
\end{lemma}

\begin{PROOF}{\ref{abd.3}}
Note that $(c) \Rightarrow (c)^-$, just let $\alpha < \lambda
\wedge  \cf(\alpha) = \aleph_\varepsilon \wedge  \varepsilon < \cf(\lambda) 
\Rightarrow f(\alpha) = \varepsilon$, clearly there is
such a function and it satisfies clause $(c)^-$.  So we can assume
$(c)^-$. Let $\sigma = \cf(\lambda)$ and $\langle 
\lambda_\varepsilon:\varepsilon < \sigma \rangle$ be a strictly
increasing sequence of regular cardinals $> \mu + \sigma$ 
with limit $\lambda$.  Let
$L_0,L^+_0,L^-_0$ be as in the proof of \ref{abd.1},
$L_{0,\varepsilon}$ be the unique interval of $L_0$ of order type (the
inverse of $\lambda_\varepsilon$) $+ \lambda_\varepsilon$, so $\langle
L_{0,\varepsilon}:\varepsilon < \sigma \rangle$ be
ia $ \subseteq $-increasing
with union $L_0,L_{0,\varepsilon}$ an interval of $L_{0,\xi}$ for
$\varepsilon < \xi < \sigma$.  We
define $g:L_0 \rightarrow \cf(\lambda)$ as follows:
if $x \in L^+_0$ then $g(x) = f(\otp(\{y \in L^+_0:y <_L x\},<) ) $
and if $x \in L^-_0$ and the order type of $(\{y \in L^+_0:x <_L
y\},<_L)$ is the inverse of $\gamma$ then $g(x) = f(\gamma)$ and let

\[
\mathbf J^* = \{\eta \in {}^{\omega >}(L_0):\eta(0) \in L_{0,0}
\text{ and } \eta(n+1) \in L_{0,g(\eta(n))} \text{ for } n < \omega\}
\]

\mn
ordered as in the proof of \ref{abd.3}.

We define ${\gs}_\zeta,{\gt}_\zeta$ as there.  We then prove that
${\gs}_\zeta,{\gt}_\zeta$ are $(\tau,\theta)$-good and 
$(\lambda,\tau)$-good as there and this suffices repeating the proof of
\ref{abd.1}. 
\end{PROOF}

\begin{discussion}
\label{abd.3a}
We may like to separate ${\gb}$ and ${\gd}$.  So
below we adapt the proof of \ref{abd.1} to do this (can do it also for
\ref{abd.3}).

A way to do this is to look at 
the forcing in \ref{abd.1} as the limit of the FS iteration
$\langle \bbP^*_i,\name{\bbQ}^*_j:i \le \mu,j < \mu \rangle$, so
the memory of $\bbQ^*_j$ is $\{i:i < j\}$ where
$\name{\bbQ}^*_j$ is $\Lim_{\gt}[\langle \bbQ_t:t \in
L^{{\gt}_j} \rangle]$.  Below we will use the limit of FS iteration
$\langle \bbP^*_i,\name{\bbQ}^*_j:j < \mu \times \mu_1 \rangle,
\bbQ^*_\zeta$ has memory $w_\zeta \subseteq \zeta$ where e.g. for $\zeta
= \mu \,\alpha +i$ where $i < \mu,w_\zeta 
= \{\kappa \beta +j:\beta \le \alpha,j \le i,(\beta,j) \ne
(\alpha,i)\}$.  Let $\bbP^* = \bbP^*_{\mu \times \mu_1}$ be 
$\cup\{\bbP_i:i < \mu \times \mu_1\}$.

Of course, $\bbQ_\zeta$ will 
be defined as $\Lim_{{\gt}_\zeta}(\bar{\bbQ})$,
the ${\gt}_\zeta$ defined as above and ${\gb} = \mu,{\gd} =
\mu_1$.  Should be easy.  If $\langle \name A_\varepsilon:
\varepsilon < \varepsilon^{\bar x} \rangle$ exemplifies ${\ga}$ in
$\mathbf V^{{\bbP}^*}$, so $\varepsilon^* \ge \mu$ then for some
$(\alpha^*,\beta^*) \in \mu \times \mu_1$ for $\kappa (= \theta)$ of the
names they involve $\{\name{\bbQ}_{\mu \alpha + \beta}:
\alpha \le \alpha^*,\beta \le \beta^*\}$ only.

Using indiscernibility on the pairs $(\alpha,\beta)$ 
to making them increase we can finish.
\end{discussion}

\begin{lemma}
\label{abd.4}
1) In Lemma \ref{abd.1}, if $\mu = \cf(\mu) \le \cf(\mu_1),
\mu_1 < \lambda$, \then \, we can change in the conclusion ${\gb} 
= {\gd} = \mu$ to ${\gb} = \mu,{\gd} = \mu_1$. 

\noindent
2) Similarly for \ref{abd.3}.
\end{lemma}

\begin{PROOF}{\ref{abd.4}}
First assume $\mu_1$ regular. 
\smallskip

\noindent
\underline{First Proof}:   Let $\mu_0 = \mu$.  In the proof 
of \ref{abd.1} for $\ell \in \{0,1\}$ using
$\mu = \mu_\ell$ gives ${\gs}^\ell_{\mu_\ell}$ and \wilog \,
${\gs}^0_{\mu_0},{\gs}^1_{\mu_1}$ are disjoint.  Let 
${\gs}$ be ${\gs}_0 +' {\gs}_1$ meaning $L^{\gs} = L^{{\gs}^0_{\mu_0}} 
+ L^{{\gs}^1_{\mu_1}}$, and for $t \in L^{{\gs}^\ell_{\mu_\ell}}$ 
we let $I^{\gs}_t =: I^{{\gs}^\ell_{\mu_\ell}}_t$ (this is 
not ${\gs}_0 + {\gs}_1$ of \ref{ad.10}).  
Now the appropriate goodness can be proved
so we can prove ${\ga} = \lambda$.  Easily we get ${\gd} \ge
\mu_1 $ % 2020-12-27 13:41 \ell$ 
and ${\gb} \le \mu_0$.  This is enough to get inequality but
to get exact values we turn to the second proof.

Instead of starting with $\langle \bbQ_i:i < \mu
\rangle$ with full memory we start with $\langle \name{\bbQ}_\zeta:
\zeta < \mu \times \mu_1 \rangle,\name{\bbQ}_\zeta$ 
with the following ``memory" if $\zeta = \mu \alpha +i,
i < \kappa,w_\zeta = \{\mu \beta +j:\beta \le \alpha,
j \le i,(\beta,j) \ne (\alpha,i)\}$.
To deal with the case $\mu_1$ is singular we should use a
$\mu$-directed index set (instead $\mu_0 \times \mu_1$) as the product
of ordered sets.   
\end{PROOF}
\newpage

\section {On related cardinal invariants} 

\underline{Explanation of \S4}:  

On Th.  \ref{au.1} you may wonder: ${\gu}$
has nothing to do with order or quite directed family, so how can we
preserve small ${\gu}$?  True, using the ``directed character" of
${\gb}$ and ${\gd}$ has been the idea, i.e. in the end 
we have $\bbP = \langle \bbP_i:i < \mu \rangle$ 
is $\lessdot$-increasing, $\bbP = \cup\{\bbP_i:i < \mu\}$ 
and $\name \eta_i$ a $\bbP_{i+1}$-name of a real dominating 
$\mathbf V^{{\bbP}_i}$.  But really what we need
   % 2020-12-27 13:42 for 
for a triple   % 2020-12-04 10:15 
$(\bbP,\name{\bar \eta} ,  \mathbb{P} ')$
as $ (\mathbb{P} _i, \name{ \eta }_i, \mathbb{P} _{i + 1 })$ 
above, is that
taking ultrapower by the $\kappa$-complete ultrafilter $D$, preserve
the property of $\name{\bar \eta}$, in our present
case $\name{\bar \eta}$ has to witness ${\gu} = \mu$.  For being a 
dominating real this is very natural (\L os theorem).
But here we shall use $\langle \name D_i:i < \mu
\rangle,\name D_i$ a $\bbP_i$-name of an ultrafilter on $\omega$ and
demand $\Rang(\name \eta_i)$   % 2020-12-04 10:16 
to be mod finite 
included in every member of $\name D_i$ 
and moreover $\name \eta_i$ is generic over $\mathbf V^{{\bbP}_i}$ for 
a forcing related to $\name D_i$.  When we like to preserve 
something in inductive construction on $\alpha < \lambda$ of $\langle
\bbP^\alpha_i:i < \mu \rangle$, % 2020-12-27 13:47 ; 
it is reasonable to have
strong induction hypothesis more than needed just for the final
conclusion.  We need here a condition on $(\bbP^\alpha_{i+1},
\name \eta^\alpha_i,\bbP^\alpha_i,\name D^\alpha_i)$ preserved 
by the ultrapower (as the relevant forcing is
c.c.c. nicely enough defined this work).

Secondly, we need in limit $\alpha$: if $\cf(\alpha) > \aleph_0$ 
straightforward if  % 2020-12-04 10:17 
not, being generic for the $\bbQ_i$ has nice enough properties so
that we can complete $\bigcup\limits_{\beta < \alpha} \name D^\beta_i$ 
to a suitable ultrafilter. 

This explains to some extent the scope of possible applications, of
course, in each case the exact inductive assumption on
$(\bbP^\alpha_{i+1},\name \eta^\alpha_i,\bbP^\alpha_i,\name Y^\alpha_i)$ with
$\name Y^\alpha_i$ a relevant witness, varies.

On continuing \S2, \S3 so eliminating the measurable here
see  \cite{Sh:F2009}, \cite{Sh:F2023}.

\begin{theorem}
\label{au.1}
Assume
\mn
\begin{enumerate}
\item[$(a)$]  $\kappa$ is a measurable cardinal
\sn
\item[$(b)$]  $\kappa < \mu = \cf(\mu) < \lambda =
\cf(\lambda) = \lambda^\kappa$.
\end{enumerate}
\mn
\Then \, for some c.c.c. forcing notion $\bbP$ of cardinality 
$\lambda$, in $\mathbf V^{\bbP}$ we have: $2^{\aleph_0} = 
\lambda,{\gu} = {\gd} = {\gb} = \mu$ and ${\ga} = \lambda$.
\end{theorem}

\begin{remark}  
Recall ${\gu} = \Min\{|{\cP}|:{\cP} \subseteq [\omega]^{\aleph_0}$ 
generates a non-principal ultrafilter on $\omega\}$.
\end{remark}

\begin{PROOF}{\ref{au.1}}
The proof is broken to definitions and claims.
\end{PROOF}

\begin{definition}
\label{au.1a}
For a filter $D$ on $\omega$ (to which all co-finite subsets of 
$\omega$ belong) let $\mathbb Q(D)$ be: 

\begin{equation*}
\begin{array}{clcr}
\{T:&T \subseteq {}^{\omega >} \omega \text{ is closed under 
initial segments, and for some} \\
  &\tr(T) \in {}^{\omega >} \omega, \text{ the trunk of } T, 
\text{ we have}: \\
  &(i) \quad \ell \le \ell g(\tr(T)) \Rightarrow T \cap {}^\ell
\omega = \{\tr(T) \restriction \ell\} \\
  &(ii) \quad \tr(T) \trianglelefteq \eta \in {}^{\omega >} \omega
\Rightarrow \{n:\eta \char 94 \langle n \rangle \in T\} \in % 2021-02-10 06:04 \name
D\}
\end{array}
\end{equation*}

\mn 
ordered by inverse inclusion.
\end{definition}

\begin{definition}
\label{au.1b}
1) Assume $S \subseteq \{i < \mu:\cf(i) \ne \kappa\}$ is unbounded in
$\mu$ (the default value is $\{i < \mu:\cf(i) \ne \kappa\}$).

Let ${\gK}_{\lambda,S}$ be the family of ${\gt}$ consisting of 
$\bar{\bbQ} = \bar{\bbQ}^{\gt} = \langle \bbP_i,\name{\bbQ}_i:i < \mu \rangle
= \langle \bbP^{\gt}_i,\name{\bbQ}^{\gt}_i:i < \mu \rangle$
and $\bar D = \bar D^{\gt} = \langle \name D_i:i < \mu$ and 
$\cf(i) \ne \kappa \rangle = \langle \name D^{\gt}_i:i \in S 
\rangle$ and $\bar\tau^{\gt} = \langle
\name \tau^{\gt_i}:i < \mu \rangle$ such that:  
\mn
\begin{enumerate}
\item[$(a)$]   $\bar{\bbQ}$ is a FS-iteration of 
c.c.c. forcing notions (and $\bbP_{\gt} = \bbP^{\gt}_\mu = 
\Lim(\bar{\bbQ}^{\gt}) = \bigcup\limits_{i < \mu} \bbP^{\gt}_i)$  % 2020-12-06 08:51 
\sn
\item[$(b)$]  if $i \in S$, then $\bbQ_i = \bbQ (\name D_i)$,
see Definition \ref{au.1a} above
\sn 
\item[$(c)$]  $\name D_i$ is a $\bbP_i$-name of a non-principal 
ultrafilter on $\omega$ when $i \in S$
\sn
\item[$(d)$]  $|\bbP_i| \le \lambda$
\sn
\item[$(e)$]  for $i \in S$ let
$\name \eta_i$ be the $\bbP_{i+1}$-name of the $\name{\bbQ}_i$-generic real

\[
\name \eta_i = \cup\{\tr(p(i)):p \in \name G_{\bbP_{i+1}}\}.
\]

and we demand: for $i <j< \mu$ of cofinality $\ne \kappa$ we have

\[
\Vdash_{\bbP_j} ``\Rang(\name \eta_i) \in \name D_j"
\]
\sn
\item[$(f)$]   $\name \tau_i$ is a $\bbP_i$-name of a
function from $\name{\bbQ}_i$ to $\{h:h$ 
 is  % 2021-02-10 06:05 
a function from a 
finite set of ordinals to ${\cH}(\omega)\}$, such that: 
\newline
$\Vdash_{\bbP_i}$ ``$p,q \in \name{\bbQ}_i$ 
are compatible % 2021-02-10 06:04 (
in $\name{\bbQ}_i$) \Iff \, the
functions $\name \tau_i(p),\name \tau_i(q)$
are compatible,  % 2021-02-10 07:51 
i.e. $\tau_i(p) \restriction (\Dom(\name \tau_i(p)) 
\cap \Dom(\name \tau_i(q)) = \tau_i(q) \restriction
(\Dom(\name \tau_i(p)) \cap \Dom(\name \tau_i(q))$ and \then \, 
they have a common upper bound $r$ such that 
$\name \tau_i(r) = \name \tau_i(p) \cup \name \tau_i(q)"$
\sn
\item[$(g)$]  if $i \in S \cap \Dom(p),p \in
\bbP_j$ and $i < j \le \mu$ then 
$\name \tau_i(p(i))$ is $\{\langle 0,\tr(p) \rangle\}$; i.e. this is
forced to hold 
\sn
\item[$(h)$]  we stipulate $\bbP_i = \{p:p$ is a function with domain $a$
finite subset of $i$ such that for each 
$j \in \Dom(p),\emptyset_{\bbP_j}$ forces that $p(j) \in 
\name{\bbQ}_j$ and it forces a value to $\name \tau_j(p(j))\}$
\sn
\item[$(i)$]  $\Vdash_{\bbP_i} ``\name{\bbQ}_i \subseteq
{\cH}_{< \aleph_1}(\gamma)$ for some ordinal $\gamma"$.
\end{enumerate}
\mn
2) Let $\gamma({\gt})$ be the minimal ordinal $\gamma$ such that $i <
\mu \Rightarrow \Vdash_{\bbP_i}$ ``if $x \in \name{\bbQ}_i$
then $\dom(\name \tau_i(x)) \subseteq \gamma$". 

\noindent
3) We let $\tau^{\gt}_i$ be the function with domain $\bbP_i$ such
that $\tau^{\gt}_i(p)$ is a function with domain 
$\{\gamma({\gt})j + \beta:j \in \Dom(p)$ and
$p \restriction j \Vdash_{\bbP_j} ``\beta \in \Dom(\name
\tau_j(p(j))"\}$ and let
$\tau^{\gt}_i(\gamma(\gt)j + \beta)$ be the value 
which $p \restriction j$ forces on $\name \tau^{\gt}_j(\beta)$.
\end{definition}

\begin{convention} 
We fix $ \lambda, \mu , S $ as in $ \ref{au.1}$, \ref{au.1b}; 
so we may write $ \mathfrak{K} $ instead
 $ \mathfrak{K} _{\lambda, S}.$
\end{convention} 

\noindent
Obviously 
\begin{subc}
\label{au.A}
${\gK} \ne \emptyset$.
\end{subc}

\begin{PROOF}{\ref{au.A}}
Should be clear.
\end{PROOF}

\noindent
Recall 
\begin{subc}
\label{au.B}
If in a universe $\mathbf V,D$ is a nonprincipal ultrafilter on $\omega$
then
\mn
\begin{enumerate}
\item[$(a)$]  $\Vdash_{\bbQ(D)} ``\{\tr(p)(\ell):\ell < \ell g(\tr(p))$ and 
$p \in \name G_{\bbQ(D)}\}$  is an infinite subset of $\omega$,
almost included in every member of $D$"
\sn
\item[$(b)$]  $\bbQ(D)$ is a c.c.c. forcing notion, even $\sigma$-centered
\sn
\item[$(c)$]  $\name \eta_i = \cup\{\tr(p):p \in \name G_{\bbQ(D)}\} \in
{}^\omega \omega$ is forced to dominate $({}^\omega \omega)^{\mathbf V}$
\sn
\item[$(d)$]  $\{p \in \bbQ(D): % 2021-02-10 06:06 [D]
\tr(p) = \eta\}$ is a directed
subsets of $\bbQ(D)$  % 2021-02-10 06:06 [D]
for every $ \eta \in {}^{ \omega > } \omega $. % 2021-02-10 06:07 
\end{enumerate}
\mn
[Note that this, in particular clause (c), does not depend on
additional properties of $D$; but as we naturally add many Cohen reals
(by the nature of the support) we may add more and then can demand
e.g. $\name D_i$ ($\cf(i) \ne \kappa$) is a Ramsey ultrafilter.]
\end{subc}

\begin{definition}
\label{au.Ba}
1)  We define $\le_{\gK}$ by: ${\gt} \le_{\gK} {\gs}$ if 
(${\gt},{\gs} \in {\gK}$ and) $i \le \mu \Rightarrow 
\bbP^{\mathfrak t}_i \lessdot \bbP^{\gs}_i$ and $i < \mu$ and $\cf(i) 
\ne \kappa \Rightarrow \Vdash_{\bbP^{\mathfrak s}_i} 
``\name D^{\gt}_i \subseteq \name D^{\gs}_i"$ and $i < \mu \Rightarrow
\Vdash_{\bbP^{\gs}_i} ``\name \tau^{\gt}_i \subseteq \name \tau^{\gs}_i"$.

\noindent
2) We say ${\gt}$ is a canonical $\le_{\gK}$-u.b. of 
$\langle {\gt}_\alpha:\alpha < \delta \rangle$ \If \,:
\mn
\begin{enumerate}
\item[$(i)$]  ${\gt},{\gt}_\alpha \in {\gK}$
\sn
\item[$(ii)$]  $\alpha \le \beta < \delta \Rightarrow {\gt}_\alpha
\le_{\gK} {\gt}_\beta \le_{\gK} {\gt}$
\sn
\item[$(iii)$]  if  $ i \in \mu \setminus S $  % 2020-12-04 10:18 $     $i < \mu$ and $\cf(i) = \kappa$ 
then 
$\Vdash_{\bbP^{\gt}_i} ``\name{\bbQ}^{\gt}_i = 
\bigcup\limits_{\alpha < \delta} \name{\bbQ}^{\gt_\alpha}_i"$.
\end{enumerate}
\mn
Note that if $\cf(\delta) > \aleph_0$ then 
we can add  % 2020-12-04 10:18 
$\Vdash_{\bbP^{\gt}_i} ``\name{\bbQ}^{\gt}_i = \bigcup\limits_{\alpha < \delta} 
\name{\bbQ}^{\gt_\alpha}_i"$ for every $i < \mu$,
so ${\gt}$ is totally determined. 

\noindent
3) We say $\langle {\gt}_\alpha:\alpha < \alpha^* \rangle$ is
$\le_{\gK}$-increasing continuous if: $\alpha < \beta < \alpha^*
\Rightarrow {\gt}_\alpha \le_{\gK} {\gt}_\beta$ and for limit $\delta <
\alpha^*,{\gt}_\delta$ is a canonical $\le_{\gK}$-u.b. of 
$\langle {\gt}_\alpha:\alpha < \delta \rangle$.  Note that we have
not said ``the canonical $\le_{\gK}$-u.b." as for $\delta < \alpha^*,
\cf(\delta) = \aleph_0$ we have some freedom in completing
$\cup\{\name D^{\gt_\alpha}_i:\alpha < \delta\}$
to an ultrafilter (on $\omega$ in $\mathbf V^{\bbP^{\gt}_i}$, 
when 
$ i \in \mu \setminus S  $).  % 2020-12-04 10:19 i < \mu,\cf(i) \ne \kappa$).
\end{definition}

\begin{subc}
\label{au.C}
If $\bbP_1 \lessdot \bbP_2$ and $\name D_\ell$ 
is a $\bbP_\ell$-name of a nonprincipal ultrafilter on $\omega$ 
for $\ell =1,2$ and $\Vdash_{\bbP_2} ``\name D_1 \subseteq \name
D_2"$, \then \, $\bbP_1 * \bbQ(\name D_1) \lessdot \bbP_2 * 
\bbQ(\name D_2)$.
\end{subc}

\begin{PROOF}{\ref{au.C}}
Why?  First, we can first force with $\bbP_1$,  so \wilog \, 
$\bbP_1$ is trivial and $D_1 \in \mathbf V$ is a nonprincipal 
ultrafilter on $\omega$.  Now clearly $p \in \bbQ(D_1) \Rightarrow p \in 
\bbQ(\name D_2)$ and $\bbQ(D_1) \models p \le q \Rightarrow 
\bbQ(\name D_2) \models p \le q$ and if $p,q \in \bbQ(D_1)$ are 
incompatible in $\bbQ(D_1)$ then they are incompatible in 
$\bbQ(\name D_2)$.  

Lastly, in $\mathbf V$, let ${\cI} = \{p_\exn:\exn < \omega\} \subseteq 
\bbQ(D_1)$ be predense in $\bbQ(D_1)$, we shall prove that 
${\cI}$ is predense in $\bbQ(\name D_2)$ in $\mathbf V^{\bbP_2}$.  

For this it suffices to note
\mn
\begin{enumerate}
\item[$\boxtimes$]  if $D_1$ is a nonprincipal ultrafilter on $\omega,
{\cI} \subseteq \bbQ(D_1)$ and $\eta \in {}^{\omega >} \omega$, \then \, the
following conditions are equivalent:
\sn
\begin{enumerate}
\item[$(a)_\eta$]  there is no $p \in \bbQ(D_1)$ incompatible with 
every $q \in {\cI}$ which satisfies tr$(p)=\eta$
\sn
\item[$(b)_\eta$]  there is a set $T$ such that:
\sn
\item[${{}}$]  $(i) \quad \nu \in T \Rightarrow \eta \trianglelefteq \nu \in p$
\sn
\item[${{}}$]  $(ii) \quad \eta \trianglelefteq \nu 
\trianglelefteq \rho \in T \Rightarrow \nu \in T$
\sn
\item[${{}}$]  $(iii) \quad$ if $\nu \in T$ then either 
$\{n:\nu \char 94 \langle n \rangle \in T\} \in D_1$ or 

$\hskip30pt (\forall n)(\nu \char 94 \langle n \rangle \notin
T) \wedge  (\exists q \in {\cI})(\nu = \tr(q))$
\sn
\item[${{}}$]  $(iv) \qquad$ there 
is a strictly decreasing function $h:T \rightarrow \omega_1$
\sn
\item[${{}}$]  $(v) \quad \eta \in p$.
\end{enumerate}
\end{enumerate}
\end{PROOF}
\bigskip

\begin{proof}
\smallskip

\noindent
\underline{Proof of $\boxtimes$}:

Straightforward.

So as in $\mathbf V,{\cI} \subseteq \bbQ(D_1)$ is predense, for every
$\eta \in {}^{\omega >} \omega$ we have $(a)_\eta$ for $D_1$ hence by
$\boxtimes$ we have also $\eta \in {}^{\omega >} \omega \Rightarrow
(b)_\eta$, but clearly if $T_\eta$ witness $(b)_\eta$ in $\mathbf V$
for $D_1$, it witnesses $(b)_\eta$ in $\mathbf V^{\bbP_2}$ for $D_2$ hence
applying $\boxtimes$ again we get: $\eta \in {}^{\omega >} \omega
\Rightarrow (a)_\eta$ in $\mathbf V^{\bbP_2}$ for $D_2$, hence ${\cI}$
is predense in $\bbQ(D_2)$ in $\mathbf V^{\bbP_2}$.  So we have proved
Subclaim \ref{au.C}.  
\end{proof}

\begin{subc}
\label{au.D}
If $\bar{\gt} = \langle {\gt}_\alpha:\alpha < \delta \rangle$ is
$\le_{\gK}$-increasing continuous and $\delta < \lambda^+$ is a limit
ordinal, \then \, it has a canonical $\le_{\gK}$-u.b. 
\end{subc}

\begin{PROOF}{\ref{au.D}}
Why?  By induction on $i \le  % 2021-02-10 06:08 < 
\mu$, we define $\bbP^{\gt}_i$ and if 
$i<\mu$ we then have $\name{\bbQ}^{\gt_i},\name \tau_i$
and $\name D_i$ (if $\cf(i) \ne \kappa$) such that the relevant 
demands (for ${\gt} \in {\gK}$ and for being 
canonical $\le_{\gK}$-u.b. of $\bar{\gt}$) hold.

Defining $\bbP^{\gt}_i$ is obvious: for $i=0$ trivially, if $i = 
j+1$ it is $\bbP^{\gt}_j * \name{\bbQ}^{\gt}_j$ 
and if $i$ is limit it is $\cup\{\bbP^{\gt}_j:j < i\}$.  

If $\bbP^{\gt}_i$ has been defined and $\cf(i) = \kappa$ we let
$\name{\bbQ}^{\gt}_i = \bigcup\limits_{\alpha < \delta}
\bbQ^{\gt_\alpha}_i$ and $\name \tau^{\gt_i} =
\bigcup\limits_{\alpha < \delta} \name \tau^{\gt_\alpha}_i$,
easy to check that they are as required. 
   If $\bbP^\gt_i$ 
has been defined and $\cf(i) \ne \kappa$, then 
$\bigcup\limits_{\alpha < \delta} D^{\gt_ \alpha }_i$ is a filter
on $\omega$ containing the co-bounded subsets, and we complete it to an
ultrafilter, 
 call it $ D^{\mathfrak{t} _i}$.

Note that there is such $\name D^{\gt_i}$ because:
\mn
\begin{enumerate}
\item[$(a)$]  for $\alpha < \delta,\bbP^{\gt_\alpha} \lessdot
\bbP^{\gt}_i$ hence $\Vdash_{\bbP^{\gt}_i} ``\name D^{\gt_\alpha}_i$ 
is a filter on $\omega$ to which all co-finite subsets of $\omega$ 
belong 
and it increases with $\alpha"$.
\end{enumerate}
\mn
Note that there will be no need for new values of the
$\name \tau_i$'s nor any freedom in defining them.
As we have proved the relevant demands on $\bbP^{\gt}_j,
\name{\bbQ}^{\gt}_j$ for $j<i$ clearly $\bbP^{\gt}_i$ is c.c.c. by 
using $\langle \name \tau_j:j < i \rangle$ and clearly $\langle
\bbP^{\gt}_\zeta,\name{\bbQ}^{\gt}_\xi:\zeta \le i,\xi < i \rangle$ is 
an FS iteration.  Now we shall prove that $\alpha < \delta \Rightarrow
\bbP^{\gt_\alpha}_i \lessdot \bbP^{\gt}_i$. 

So let ${\cI}$ be a predense subset of $\bbP^{\gt_\alpha}_i$ and
$p \in \bbP^{\gt}_i$ and we should prove that $p$ is compatible with some
$q \in {\cI}$ in $\bbP^{\gt}_i$; we divide the proof to cases.
\bigskip

\noindent
\underline{Case 1}:  $i$ is a limit ordinal.

So $p \in \bbP^{\gt}_j$ for some $j < i$,  let ${\cI}' = \{q \restriction
j:q \in {\cI}\}$, so clearly ${\cI}'$ is a predense subset of
$\bbP^{\gt_\alpha}_j$ (as ${\gt}_\alpha \in {\gK}$).  
By the induction hypothesis, in $\mathbb{P}  % 2021-02-10 06:10 
^{\gt}_j$ the condition $p$ is 
compatible with some $q' \in {\cI}'$; so let $r' \in 
\bbP^{\gt}_j$ be a common upper bound of $q',p$ 
recalling that % 2021-02-10 06:10 and let 
$q' = q \restriction j$ where $q \in {\cI}$.  So
$r' \cup (q \restriction [j,i)) \in \bbP^{\gt}_i$ is a common upper bound of
$q,p$ as required.
\bigskip

\noindent
\underline{Case 2}:  $i = j+1,\cf(j) = \kappa$.

If $j \notin \Dom(p)$ it is trivial.  So \wilog \, for some 
$\beta < \delta,p(j)$ is a $\bbP^{\gt_\beta}_j$-name of a member 
of $\name{\bbQ}^{\gt_\beta}_j$; and \wilog \, $\alpha \le 
\beta < \delta$.  By the induction hypothesis 
$\bbP^{\gt_\beta}_j \lessdot \bbP^{\gt}_j$ hence there is $p' \in
\bbP^{\gt_\beta}_j$ such that $[p' \le p'' \in \bbP^{\gt_\beta}_j
\Rightarrow p'',p \restriction j$ are compatible in $\bbP^{\gt}_j]$. 

Let

\begin{equation*}
\begin{array}{clcr}
{\cJ} = \{q' \restriction j:&q' \in \bbP^{\gt_\beta}_i 
\text{ and } q' \text{ is above some member of } {\cI} \\
  &\text{and } q' \restriction j \Vdash_{\bbP^{\gt_\beta}_j} ``p(j) 
\le^{\name{\bbQ}^{\gt_\beta}_j} q'(j)"\}.
\end{array}
\end{equation*}

\mn
Now ${\cJ}$ is a dense subset of $\bbP^{\gt_\beta}_j$ (since if
$q \in \bbP^{\gt_\beta}_j \text{ then } q \cup \{\langle j,p(j) \rangle\}$
belongs to $\bbP^{\gt_\beta}_i$ hence is compatible with some member of
${\cI}$). 

Hence $p'$ is compatible with some $q'' \in {\cJ}$ 
(in $\bbP^{\gt_\beta}_j$),   % 2021-02-10 06:11 
so there is $r$ such that $p' \le r \in
\bbP^{\gt_\beta}_j,q'' \le r$.  As $q'' \in {\cJ}$ there is $q' \in
\bbP^{\gt_\beta}_i$ such that $q' \restriction j = q'',q'$ 
is above some $q^* \in \cI$ and $q' \rest j \Vdash ``p(j) 
\le^{\name{\bbQ}^{\gt_\beta}_j} q'(j)"$.

As $\bbP^{\gt_\beta}_j \models ``p' \le r \wedge q' \restriction j =
q'' \le r"$ and by the choice of $p'$ there is $p^* \in \bbP^{\gt}_j$ 
above $r$ (hence above $p'$ and above $q'' = q' \restriction j)$,
and above $p \restriction j$.  Now let $r^* = p^* \cup 
(q'' \restriction \{j\})$, clearly $r^* \in \bbP^{\gt}_i$ is above 
$p \restriction j$ and $r^* \restriction j$ forces that $r^*(j)$ is above
$p \restriction \{j\}$.  Clearly $r^* \restriction j$ is above $r$ and
$r^*$ is also  % 2021-02-10 06:12 also is
above $q^* \in {\cI}$ so we are done.
\bigskip

\noindent
\underline{Case 3}:  $i = j+1,j \in S $ % 2021-02-10 06:12 \cf(j) \ne \kappa$.

Use Subclaim \ref{au.C} above. 

So we have dealt with $\alpha < \delta \Rightarrow \bbP^{\gt_\alpha}_i 
\lessdot \bbP^{\gt}_i$.

Clearly we are done.  
\end{PROOF}

\begin{subc}
\label{au.E}
If ${\gt} \in {\gK}$ and $E$ is a $\kappa$-complete
non-principal ultrafilter on $\kappa$, \then \, we can find ${\gs}$
such that:
\mn
\begin{enumerate}
\item[$(i)$]   ${\gt} \le_{\gK} {\gs} \in {\gK}$
\sn
\item[$(ii)$]  there is $\langle \mathbf k_i,\mathbf j_i:i < \mu,\cf(i) 
\ne \kappa)$ such that:
\sn
\begin{enumerate}
\item[$(\alpha)$]  $\mathbf k_i$ is an isomorphism from
$(\bbP^{\gt}_i)^\kappa/E$ onto $\bbP^{\mathfrak s}_i$
\sn
\item[$(\beta)$]  $\mathbf j_i$ is the canonical embedding of
$\bbP^{\gt}_i$ into $(\bbP^{\gt}_i)^\kappa/E$
\sn
\item[$(\gamma)$]  $\mathbf k_i \circ \mathbf j_i =$ identity on
$\bbP^{\gt}_i$
\end{enumerate}
\sn
\item[$(iii)$]  $\name D^{\gs_i}$ is the image of
$(\name D_i)^\kappa/E$ under $\mathbf k_i$ and similarly
$\name \tau^{\gs_i}$ \If \, $i < \mu,\cf(i) \ne \kappa$
\sn
\item[$(iv)$]  if $i < \mu,\cf(i) = \kappa$, then
$\name \tau^{\gs_i}$ is defined such that, for $j < \kappa,\cf(j) 
\ne \kappa$ we have $\mathbf k_j$ is an isomorphism from
$(\bbP^{\gt}_i,\gamma',\tau^{\gt}_i)^\kappa/D$ onto 
$(\bbP^{\gs}_i,\gamma'',\tau^{\gs_i}_i)$ for some ordinals   % 2021-02-10 06:13 \mathfrak{s} , i
$\gamma',\gamma''$ (except that we do not require that the map from
$\gamma'$ to $\gamma''$ preserves order).
\end{enumerate}
\end{subc}

\begin{PROOF}{\ref{au.E}}
Straightforward.  

Note that if $\cf(i) = \kappa,i < \mu$ then $\name{\bbQ}^{\gs}_i$ is 
isomorphic to $\bbP^{\gs} _{i+1}/\bbP^{\gs}_i$ which 
is c.c.c. as by \L o\'s theorem for the logic $\bbL_{\kappa,\kappa}$ we have 
$\bigcup\limits_{j<i} (\bbP^{\gt}_j)^\kappa/E \lessdot
(\bbP^{\gt}_{i+1})^\kappa/E$, similarly for
$\name \tau_i$ which guarantees that the quotient 
is c.c.c., too (actually $\name \tau_i$ is not needed for
the c.c.c. here).  
\end{PROOF}

\begin{subc}
\label{au.F}
If ${\gt} \in {\gK}$ \then \, $\Vdash_{\bbP^{\gt}_\mu} ``{\gu} = {\gb} = 
{\gd} = \mu"$. 
\end{subc}

\begin{PROOF}{\ref{au.F}}
In $\mathbf V^{\bbP^{\gt}_\mu}$, the family ${\cD} = \{\Rang(\name \eta_i):
i < \mu$ and $\cf(i) \ne \kappa\} \cup \{[n,\omega):n < \omega\}$ generates a 
filter on 
  ${\cP}(\omega)^{\mathbf V[\bbP^{\gt}_\mu]}$, as
$\Rang(\name \eta_i) \in [\omega]^{\aleph_0}, 
i < j < \mu$ and $\cf(i) \ne \kappa$ and $\cf(j) \ne \kappa 
\Rightarrow \Rang(\name \eta_j) \subseteq^* \Rang(\name \eta_i)$.  

Also it is an ultrafilter as
${\cP}(\omega)^{\mathbf V[\bbP^{\gt}_\mu]} = \bigcup\limits_{i < \mu}
{\cP}(\omega)^{\mathbf V[\bbP^{\gt}_i]}$ and if $i < \mu$, then $\Rang
(\name \eta_{i+1})$ induces an ultrafilter on 
${\cP}(\omega)^{\mathbf V[\bbP^{\gt}_{i+1}]}$.  So ${\gu} \le \mu$.  
Also $({}^\omega \omega)^{\mathbf V[\bbP^{\gt}_\mu]} 
= \bigcup\limits_{i < \mu}({}^\omega \omega)^{\mathbf V[\bbP^{\gt}_i]},
({}^\omega \omega)^{\mathbf V[\bbP^{\gt}_i]}$ is increasing with
$i$ and if $\cf(i) \ne \kappa$ then $\name \eta_i \in {}^\omega \omega$ 
dominates
$({}^\omega \omega)^{\mathbf V[\bbP^{\gt}_i]}$ by
Subclaim \ref{au.B}, so ${\gb} = {\gd} = \mu$ as in previous cases.  

Lastly, always ${\gu} \ge {\gb}$ hence ${\gu} = \mu$.]
\end{PROOF}
\bigskip

\noindent
Now we define ${\gt}_\alpha \in {\gK}$ for $\alpha \le \lambda$ 
by induction on $\alpha$ satisfying $\langle {\gt}_\alpha:\alpha 
\le \lambda \rangle$ is $\le_{\gK}$-increasing continuous 
such that ${\gt}_{\alpha +1}$ is gotten from ${\gt}_\alpha$ as 
in Subclaim \ref{au.E}.  

Let $\bbP = \bbP^{\gt_\lambda}_\mu$, so $|\bbP| \le \lambda$ 
hence $(2^{\aleph_0})^{\mathbf V^{\bbP}} \le
(\lambda^{\aleph_0})^{\mathbf V}$ and easily equality holds. 
\bigskip   

\noindent
We finish by 
\begin{subc}
{\label{au.G}} 
We have\footnote{
 recall $ \lambda $ is regular; if we allow
 $ \lambda $ singular we have to use 
 $ \cf( \lambda)$.
       }   % 2021-02-10 06:14 sop hevarat julayim
$\Vdash_{\bbP^{g\mathfrak{t}  _ \alpha }  % 2021-02-10 06:15 
_\lambda} ``{\ga} \ge   % 2020-12-27 14:07 \cf(
\lambda"$.
\end{subc}
\begin{PROOF} {\ref{au.G}}
Why?  Assume toward a contradiction that $\theta <   % 2020-12-27 14:09 \cf(
\lambda$ and
$p \in \bbP$ and $p \Vdash_{\bbP} ``\name{\cA} = \{\name A_i:
i < \theta\}$ is a MAD family; i.e.
\mn
\begin{enumerate}
\item[$(i)$]  $A_i \in [\omega]^{\aleph_0}$
\sn
\item[$(ii)$]  $i \ne j \Rightarrow |\name A_i \cap \name A_j| < \aleph_0$
\sn
\item[$(iii)$]  under (i) + (ii), $\name{\cA}$ is maximal".
\end{enumerate}
\mn
Without loss of generality $\Vdash_{\bbP} ``\name A_i \in
[\omega]^{\aleph_0}"$.  As ${\ga} \ge {\gb} = \mu$ by Subclaim
\ref{au.F}, we have $\theta \ge \mu$.  For each $i < \theta$ and 
$m < \omega$ there is a maximal
antichain $\langle p_{i,m,n}:n < \omega \rangle$ of $\bbP$ and there is a
sequence $\langle \mathbf t_{i,m,n}:n < \omega \rangle$ of truth values such
that $p_{i,m,n} \Vdash ``(m \in \name A_i) \equiv \mathbf t_{i,m,n}"$.
We can find countable $w_i \subseteq \mu$ such that
$\bigcup\limits_{m,n < \omega} \Dom(p_{i,m,n}) \subseteq w_i$. Possibly
increasing $w_i$ retaining countability, we can find 
$\langle R_{i,\gamma}: \gamma \in w_i \rangle$ such that:
\mn
\begin{enumerate}
\item[$(\alpha)$]  $w_i$ has a maximal element and $\gamma \in w_i
\backslash \{\max(w_i)\} \Rightarrow \gamma +1 \in w_i$
\sn
\item[$(\beta)$]  $R_{i,\gamma}$ is a countable subset of
$\bbP^{\gt_\lambda}_\gamma$ and $q \in R_{i,\gamma} \Rightarrow
\Dom(q) \subseteq w_i \cap \gamma$
\sn
\item[$(\gamma)$]  for $\gamma_1 < \gamma_2$ in 
$w_i,q \in R_{i,\gamma_2} \Rightarrow q \restriction 
\gamma_1 \in R_{i,\gamma_1}$
\sn
\item[$(\delta)$]  for $\gamma_1 \in w_i,\gamma \in \gamma_1 \cap
w_i$ and $q \in R_{i,\gamma_1}$ the
$\bbP^{\gt}_\gamma$-name $q(\gamma)$ involves
$\aleph_0$ maximal antichains all included in $R_{i,\gamma}$
\sn
\item[$(\varepsilon)$]  $\{p_{i,m,n}:m,n <\omega\} \subseteq 
R_{i,\text{max}(w_i)}$.
\end{enumerate}
\mn
As $\cf(\lambda) > \aleph_0$ (as $\mu < \lambda = \cf(\lambda)$
by the assumption of Theorem \ref{au.1}) we have
$\bbP^{\gt}_\mu = \bigcup\limits_{\alpha < \lambda} 
\bbP^{\gt_\alpha}_\mu$.
Clearly for some $\alpha < \lambda$ 
we have $\cup \{R_{i,\gamma}:i < \theta,\gamma \in w_i\} \subseteq 
\bbP^{{\gt}_\alpha}_\mu$.   But $\bbP^{{\gt}_\alpha}_\mu \lessdot 
\bbP^{{\gt}_\lambda}_\mu$.  So $\Vdash_{{\bbP}^{{\gt}_\alpha}_\mu} 
``\name{\cA} = \{\name A_i:i < \theta\}$ is MAD".

Now, letting $\mathbf j$ be the canonical elementary embedding of 
$\mathbf V$ into $\mathbf V^\kappa/D$, we know:
\mn
\begin{enumerate}
\item[$(*)$]  in $\mathbf V^\kappa/D,\mathbf j(\name{\cA})$ is
a $\mathbf j(\bbP^{{\gt}_\alpha}_\mu)$-name of a MAD family.
\end{enumerate}
\mn
As $\mathbf V^\kappa/D$ is $\kappa$-closed, for 
c.c.c. forcing notions things are absolute enough but
$\{\mathbf j(i):i < \mu\}$ is not $\{i:\mathbf V^\kappa/D \models i < \mathbf j
(\mu)\}$, so in $\mathbf V$, it is forced for $\Vdash_{{\mathbf j}
(\bbP^{\mathfrak{ t}_\alpha }_\mu)}$,   % 2021-02-10 06:16 \alpha 
that $\{\mathbf j(\name A_i):i < \mu\}$ is not
MAD! 

Chasing arrows, clearly $\Vdash_{{\bbP}^{{\gt}_{\alpha +1}}_\mu} 
``\{\name A_i:i < \theta\}$ is not MAD" as required.
\end{PROOF} % 2020-12-27 14:10 \end{subc}

\begin{discussion}
\label{au.3}
1)  We can now look at other problems, like what can be the
order and equalities among ${\gd}, {\gb}, {\ga}, {\gu}$;  have not
considered it.  I have considered having ${\ga} = \mu$ but there
was a problem.

\noindent 
2)
 (2020)  In \ref{au.1} 	
We can add $ \mathfrak{p} = \mathfrak{t} = \mu$ % 2020-12-06 08:55 p=t= \mu
proving as in \ref{au.D}.  % 2020-12-06 08:55 4.9
Let me elaborate:
in Definition \ref{au.1b} % 2020-12-06 08:56 4.4 
(our forcing is $ \mathbb{P} _\mu$ 
for such $ \mathfrak{t} $),
we have an ultrafilter  generated by a sequence of subsets
of $\omega $  which is decreasing modulo finite;
see clauses (c) and (e).

\noindent 
3) % 2021-02-10 06:18 
So $ \mathbb{P} _\mu$ forces $ \mathfrak{s} $  is at least $ \mu$.
 But always $ \mathfrak{s} $  is at  most $ \mathfrak{u} $
 so in \ref{au.1} % 2020-12-06 08:59  4.1 
we can add  $ \mathfrak{u} =\mu$.  % 2021-02-10 06:19 \mathfrak{s} 
\end{discussion}
\newpage

%STOP HERE THE CONVERSION

\bibliographystyle{amsalpha}
\bibliography{shlhetal}

\end{document}